\documentstyle[amscd,amssymb,verbatim,diagrams]{amsart}
\newcommand{\Id}{\operatorname{Id}}

\newcommand{\ch}{\operatorname{ch}}

\renewcommand{\mod}{\operatorname{mod}}

\newcommand{\Tor}{\operatorname{Tor}}
\newcommand{\und}{\underline}

\newcommand{\OO}{{\cal O}}

\newcommand{\pr}{\operatorname{pr}}

\newcommand{\coker}{\operatorname{coker}}
\newcommand{\DD}{{\cal D}}
\newcommand{\NN}{{\cal N}}

\newcommand{\II}{{\cal I}}
\newcommand{\Gal}{\operatorname{Gal}}
\newcommand{\tot}{\operatorname{tot}}
\newcommand{\BB}{{\cal B}}

\newcommand{\G}{{\Bbb G}}

\newcommand{\hra}{\hookrightarrow}
\newcommand{\AV}{{\cal AV}}
\newcommand{\HHom}{{\cal H}om}

\newcommand{\Coh}{\operatorname{Coh}}
\newcommand{\GG}{{\cal G}}
\newcommand{\CC}{{\cal C}}

\newcommand{\Spec}{\operatorname{Spec}}

\newcommand{\si}{\sigma}

\newcommand{\Pic}{\operatorname{Pic}}

\newcommand{\ga}{\gamma}
\newcommand{\de}{\delta}

\renewcommand{\ker}{\operatorname{ker}}

\newcommand{\bff}{{\bf f}}

\newcommand{\bG}{{\bf G}}
\numberwithin{equation}{subsection}

\newtheorem{thm}{Theorem}[subsection]
\newtheorem{prop}[thm]{Proposition}
\newtheorem{lem}[thm]{Lemma}
\newtheorem{cor}[thm]{Corollary}
\newenvironment{rem}{\vspace{3mm}\noindent
{\bf Remark.}}{\vspace{3mm}}
\newenvironment{defi}{\vspace{3mm}\noindent
{\bf Definition.}}{\vspace{3mm}}

\newenvironment{ex}{\vspace{3mm}\noindent
{\bf Example.}}{\vspace{3mm}}
\newenvironment{exs}{\vspace{3mm}\noindent
{\bf Examples.}}{\vspace{3mm}}

\newcommand{\Pf}{\noindent {\it Proof}}
\newcommand{\id}{\operatorname{id}}

\newcommand{\ov}{\overline}

\newcommand{\ra}{\rightarrow}
\renewcommand{\AA}{{\cal A}}
\renewcommand{\AV}{{\cal AV}}
\newcommand{\For}{{\cal F}or}
\newcommand{\Aff}{{\cal A}ff}
\newcommand{\Ab}{{\cal A}b}
\newcommand{\Sh}{{\cal S}h}
\newcommand{\FF}{{\cal F}}

\newcommand{\MM}{{\cal M}}
\newcommand{\JJ}{{\cal J}}

\newcommand{\XX}{{\cal X}}

\newcommand{\PP}{{\cal P}}
\newcommand{\QQ}{{\cal Q}}
\newcommand{\VV}{{\cal V}}
\newcommand{\SS}{{\cal S}}
\newcommand{\UU}{{\frak U}}
\newcommand{\LL}{{\cal L}}
\newcommand{\Du}{{\Bbb D}}

\newcommand{\Hom}{\operatorname{Hom}}
\newcommand{\Biext}{\operatorname{Biext}}
\newcommand{\Ext}{\operatorname{Ext}}

\newcommand{\PPic}{{\cal P}ic}
\newcommand{\Eext}{EXT}
\newcommand{\EExt}{EXT}
\renewcommand{\a}{\alpha}
\renewcommand{\b}{\beta}
\newcommand{\om}{\omega}
\newcommand{\De}{\Delta}

\newcommand{\th}{\theta}

\newcommand{\Z}{{\Bbb Z}}

\renewcommand{\L}{{\Bbb L}}
\newcommand{\La}{\Lambda}
\newcommand{\Ga}{\Gamma}

\newcommand{\wt}{\widetilde}
\newcommand{\ot}{\otimes}

\newcommand{\sub}{\subset}
\newcommand{\ed}{\qed\vspace{3mm}}

\newcommand{\Qcoh}{\operatorname{Qcoh}}
\newcommand{\clim}{\overrightarrow{c}}

\newcommand{\Spf}{\operatorname{Spf}}
\newcommand{\TTor}{{\cal T}or}
\newcommand{\Auteq}{\operatorname{Auteq}}
\newcommand{\Aut}{\operatorname{Aut}}

\title{Kernel algebras and generalized Fourier-Mukai transforms}
\author{A. Polishchuk}
\address{Department of Mathematics, University of Oregon, Eugene, OR 97405}
\email{apolish@@math.uoregon.edu}
\thanks{Supported in part by NSF grant}

\begin{document}
\begin{abstract} We introduce and study {\it kernel algebras}, i.e., algebras in the category of
sheaves on a square of a scheme, where the latter category is equipped with a monoidal structure
via a natural convolution operation. We show that many interesting categories, such as 
$D$-modules, equivariant sheaves and their twisted versions, arise
as categories of modules over kernel algebras.
We develop the techniques of constructing derived equivalences between these module categories.
As one application we generalize the results of \cite{PR} 
concerning modules over algebras of twisted differential operators 
on abelian varieties. As another application we recover and generalize the results of 
Laumon~\cite{Lau} concerning an analog of the Fourier transform
for derived categories of quasicoherent sheaves on a dual pair of generalized 1-motives. 
\end{abstract}
\maketitle

\bigskip

\centerline{\sc Introduction}

The classical Morita theory gives a way to construct all equivalences between the categories
of modules over rings $A$ and $B$ in terms of $A$-$B$-bimodules. This theory was generalized
to bounded derived categories of module categories by Ricard~\cite{Ric}. On the other hand,
starting from the pioneering work of Mukai~\cite{Mukai} algebraic geometers got interested in the study
of the functors between bounded derived categories of 
coherent sheaves on projective varieties $X$ and $Y$
associated with complexes of coherent sheaves ({\it kernels}) on $X\times Y$ (we refer to such
functors as {\it integral transforms}). The result of Orlov~\cite{Orlov} states that in the case when $X$ and $Y$ are smooth projective varieties all exact equivalences
between these derived categories are of this form. Note that the latter theory belongs to the commutative world, although there exist generalizations to stacks (see \cite{Kaw}, \cite{BZFN}) and
to twisted sheaves (see \cite{CanSt}). There are interesting ``noncommutative" categories of geometric origin that are left out from this picture, such as categories of $D$-modules. It seems that there is no
straightforward generalization of the above theory to this case. One way to include both the
derived Morita theory and the Fourier-Mukai transforms into one framework is by working with
dg-categories, as outlined in \cite{Toen}. In the present paper we propose a more specialized extension
of the techniques of integral transforms to a (partially) noncommutative world that does not require passing to dg-categories. Our immediate goal was to understand and generalize two
concrete examples, namely, the Fourier duality for generalized $1$-motives (see \cite{Lau}) 
and for modules over algebras of twisted differential operators on abelian varieties (see \cite{PR}). 
The resulting framework, although much more limited in applicability than the dg-techniques,
covers both these examples and their generalizations.

The main idea is quite simple. Composition of integral transforms between derived categories of sheaves corresponds to taking convolution of kernels (see section \ref{circle-sec}). 
In particular, we get a monoidal structure on the category of kernels on $X\times X$ with 
the unit object $\De_*\OO_X$, where $\De:X\to X\times X$ is the diagonal. 
Thus, in the context of derived equivalences
it is natural to work with algebra objects with respect to this monoidal structure.
This is what we call {\it kernel algebras}. However, if one wants to define analogs of module
categories one needs to impose a kind of flatness condition: the transform from the derived category
of sheaves on $X$ to itself defined by our kernel algebra
should preserve the abelian subcategory of sheaves. We call kernel algebras with this property {\it pure}. 
Generalizing Theorem 6.5 of \cite{PR}, we prove that an equivalence of derived 
categories of sheaves on $X$ and $Y$ extends to a correspondence between kernel algebras
over $X$ and $Y$ such that the corresponding derived categories of modules are equivalent,
provided both kernel algebras are pure (and other technical assumptions are satisfied).
Note that any D-algebra in the sense of \cite{BB} can be
viewed as a pure kernel algebra, so the context of \cite{PR} is embedded into our framework. 
We also show that one can associate a kernel algebra with
an action of a finite group scheme $G$ (or a formal group) on the derived category $D_{qc}(X)$ of 
quasicoherent sheaves on $X$ (provided this action is given by integral transforms). 
In the case when such an action is geometric (i.e., induced by the action of
$G$ on $X$) the category of modules over the corresponding pure kernel algebra is equivalent to
the category of $G$-equivariant sheaves on $X$. More generally, one can associate a pure kernel
algebra with an appropriate class of groupoids (resp., formal groupoids). We also consider
in detail the situation when $G$ acts on the category of sheaves on $X$ 
by autoequivalences induced by automorphisms of $X$ combined with tensoring by line bundles.
Modules over the corresponding $G$-algebras can be interpreted as twisted $G$-equivariant
sheaves on $X$.

Let us explain how kernel algebras can be used to establish an analog of the Fourier-Mukai duality
for the derived categories of sheaves on generalized $1$-motives considered in \cite{Lau}. 
We consider a slightly broader class of generalized $1$-motives than in \cite{Lau}, namely, 
complexes $[G\to E]$ concentrated in degrees $-1$ and $0$, 
where $G$ is a commutative formal group (satisfying
some finiteness assumptions) and $E$ is a commutative algebraic group.\footnote{Thus, unlike in \cite{Lau}, we allow $G$ to have torsion and $E$ to be nonconnected.} Quasicoherent
sheaves on $[G\to E]$ are simply $G$-equivariant quasicoherent sheaves on $E$,
so they can be described by an appropriate pure kernel algebra on $E$. Note that $E$ is an extension
of an abelian variety $A$ by an affine algebraic group. Using an appropriate notion
of the push-forward for kernel algebras we can equivalently describe sheaves on $K=[G\to E]$ as
modules over some pure kernel algebra $\AA(K)$ on $A$. The dual $1$-motive
$\Du(K)$ is presented by a complex $[G'\to E']$, where $E'$ is an extension of the dual abelian
variety $\hat{A}$ by an affine group. Applying the same procedure to $\Du(K)$ we get
a kernel algebra $\AA(\Du(K))$ on $\hat{A}$ responsible for sheaves on $\Du(K)$.
Furthermore, we interpret $\AA(K)$ (resp., $\AA(\Du(K))$) as the kernel algebra associated with
a homomorphism from $G\times G'$ to $\Auteq(D_{qc}(A))$ (resp., $\Auteq(D_{qc}(\hat{A}))$). Thus, sheaves
on $K$ (resp., $\Du(K)$) can be viewed 
as twisted $G\times G'$-equivariant sheaves on $A$ (resp., $\hat{A}$). Using this point of view
we check that $\AA(K)$ and $\AA(\Du(K))$
correspond to each other under the standard Fourier-Mukai transform that
gives an equivalence between the derived categories of sheaves on $A$ and on $\hat{A}$.
This immediately leads to the desired equivalence of derived categories of sheaves on $K$ and
on $\Du(K)$. 

The paper is organized as follows. In section \ref{prelim-sec} we collect
some results on quasicoherent sheaves over schemes and formal schemes, as well as few facts
about biextensions. Most of these results are well known. However, in some cases we had to prove
simple extensions of the known theorems to new situations (for example, we give three versions
of the base change formula for sheaves on formal schemes). Also, we prove here some
auxiliary statements involving duality for sheaves on formal schemes.
Section \ref{keralg-sec} develops the theory of kernel algebras. Here we define a convolution
operation for modules over pure kernel algebras satisfying an additional technical assumption. The main result of this section is Theorem \ref{PR-thm} that extends equivalences
between derived categories of sheaves on schemes $X$ and $Y$ to derived categories of modules
over kernel algebras over $X$ and $Y$.
In section \ref{kerrep-sec} we define the notion of a {\it kernel representation} of a group scheme
(or a formal group) $G$ over a scheme $X$, which is a scheme-theoretic version of a homomorphism
from a discrete group to $\Auteq(D_{qc}(X))$.  
In the case when $G$ is a finite flat group scheme (resp., formal group) we construct a kernel algebra
over $X$ associated with a kernel representation. We study in particular such kernel algebras corresponding to homomorphisms from $G$ to $\Aut(X)\ltimes\Pic(X)$, or equivalently,
with $1$-cocycles of $G$ with values in $\Pic(X)$. In the case when $G$ is commutative
the interplay between such algebras and the Cartier duality plays a crucial rule in applications to
$1$-motives. On the other hand, Theorem \ref{equiv-sh-thm} establishes an equivalence between the category of sheaves equivariant
with respect to an action of a groupoid (resp., formal groupoid) on $X$ and modules over the corresponding kernel algebra.
We also construct an equivalence of derived categories for twisted $G$-equivariant sheaves on a pair of dual abelian schemes (see Theorem \ref{twist-eq-thm}) generalizing equivalences
for modules over twisted differential operators considered in \cite{PR}.  
Finally, in section \ref{FM-sec} we apply the developed techniques to constructing
an analog of Fourier-Mukai transform for generalized $1$-motives and to generalizing Theorem
\ref{twist-eq-thm} to this context.

\medskip

\noindent
{\it Acknowledgment}. I am grateful to Alexander Beilinson, Oren Ben-Bassat,
Joseph Lipman and Daniel Hernandez Ruiperez for helpful discussions.

\medskip

\noindent{\it Conventions and notation}. 
All schemes in this paper are assumed to be noetherian, and all
formal schemes are assumed to be locally noetherian. 
For our conventions on formal groups see the beginning of section \ref{kerrep-sec}.
By an {\it algebraic group over a field $k$} we mean a group scheme of finite type over $k$.
 
Starting from section \ref{keralg-sec}, all the functors are assumed to be derived. Thus,
for a morphism $f$ between schemes 
$f_*$ (resp., $f^*$) will denote the derived push-forward (resp., pull-back) of complexes of sheaves
of $\OO$-modules. Similarly, $\ot$ will denote the derived tensor product (also starting from section
\ref{keralg-sec}).
Recall that these functors are defined on the unbounded derived categories 
of $\OO$-modules (see \cite{Spalt} and ch.~2 of \cite{Lipman}). Similar conventions are applied for
formal schemes.  


All inductive (resp., projective) systems are assumed to be small and filtered.
We denote by $\lim$ (resp., $\projlim$) the inductive (resp., projective) limit of an inductive
(resp., projective) system.


\section{Preliminaries}\label{prelim-sec}

\subsection{Quasicoherent sheaves on schemes}

Our primary source 
regarding derived categories of quasicoherent sheaves on schemes
are the notes \cite{Lipman} (although we are not aiming at the same level of generality).
Recall that we assume all schemes to be noetherian, so every morphism between
them is quasi-compact and quasi-separated (see \cite{EGAI}, (6.1.13)).
We denote by $\Qcoh(X)$ (resp., $\Coh(X)$)
the category of quasicoherent (resp., coherent) sheaves of $\OO$-modules on a scheme $X$.
We denote by $D_{qc}(X)$ (resp., $D_c(X)$) the full subcategory in the derived category of 
$\OO$-modules consisting of complexes with quasicoherent (resp., coherent) cohomology.

First, we observe that the standard (derived) functors with complexes of $\OO$-modules
preserve the subcategories $D_{qc}$.
For the pull-back functors and for the tensor product this requires no assumptions on the schemes
involved (see \cite{Lipman}(3.9.1) and (2.5.8)). For the push-forward functor this is true 
because we consider only noetherian schemes (\cite{Lipman} (3.9.2)). For the same reason
for a map $f:X\to Y$ the functor 
$Rf_*:D_{qc}(X)\to D_{qc}(Y)$ is bounded above (\cite{Lipman} (3.9.2)), 
i.e., it sends $D_{qc}^{\le 0}(X)$ to $D_{qc}^{\le n}(Y)$ for some $n$ (it is also bounded below
for trivial reasons).
In the case when $f:X\to Y$ is an affine morphism the functor $f_*$ is exact
and the induced functor $D_{qc}(X)\to D_{qc}(Y)$ is
conservative (this follows from the fact that for an affine scheme $X=\Spec(A)$ the functor of global sections induces an equivalence of $D_{qc}(X)$ with $D(A-\mod)$). 

The following result follows immediately from Lemma (3.9.3.1) in \cite{Lipman}. 
Note that the crucial case of $R^0f_*$ traces back to the proof of Theorem 3.10 in \cite{Godement}. 

\begin{lem}\label{lim-push-lem} 
Let $f:X\to Y$ be a morphism of schemes. Then
$R^if_*$ commutes with inductive limits. 
\end{lem} 

We will use the following versions of the projection formula and the base change formula.

\begin{thm}\label{proj-formula-thm} (\cite{Lipman}(3.9.4))
Let $f:X\to Y$ be a morphism of schemes.
Then the natural map
$$Rf_*F\otimes^{\L} G\to Rf_*(F\otimes^{\L} Lf^*G)$$
is an isomorphism for $F\in D_{qc}(X)$, $G\in D_{qc}(Y)$.
\end{thm}

\begin{thm}\label{base-change-thm} (\cite{Lipman}(3.9.5) and (3.10.3); \cite{BBH}, Prop. A.85)
Assume that we have a  cartesian diagram 
\begin{diagram}
X' & \rTo{v} & X\\
\dTo{f'} & &\dTo{f}\\
Y' & \rTo{u} & Y
\end{diagram}
where either $u$ or $f$ is flat.
Then for $F\in D_{qc}(X)$ the base change map $u^*Rf_*F\to Rf'_*v^*F$
is an isomorphism.
\end{thm}

A scheme $X$ is called {\it semi-separated} 
if the diagonal morphism $X\to X\times X$ is affine. 
It is easy to see that any morphism from an affine scheme to a semi-separated
scheme is affine. For a quasicoherent sheaf on a semi-separated scheme 
one can construct 
a quasi-isomorphic complex of flat quasicoherent sheaves (using Cech resolution, 
see \cite{AJL-loc}, Prop.1.1). 
This easily implies the following statement (see also \cite{Murfet}, Prop. 16).

\begin{lem}\label{flat-cover-lem} 
Let $X$ be a semi-separated scheme. Then for
every quasicoherent sheaf $F$ on $X$ there exists a surjection $P\to F$, where
$P$ is flat quasicoherent. Moreover, this surjection can be constructed functorially in $F$.
\end{lem}

Recall that for a finite morphism $f:X\to Y$ the functor
$f_*:D_{qc}(X)\to D_{qc}(Y)$ admits the right adjoint
$f^!:D_{qc}(Y)\to D_{qc}(X)$. We need the following simple fact about this functor.

\begin{prop}\label{finite-duality-prop}
Let $f:X\to Y$ be a finite flat morphism. Then $f^!\OO_Y$ is a line bundle on $X$ and for 
$F\in D_{qc}(Y)$ the natural morphism
$$f^!\OO_X\otimes f^*F\to f^!F$$
in $D_{qc}(X)$ is an isomorphism. Given any morphism $u:U\to Y$ let us consider
the cartesian square
\begin{diagram}
V &\rTo{v}& X\\
\dTo{g} &&\dTo{f}\\
U &\rTo{u}&Y
\end{diagram}
Then one has a natural isomorphism of functors $Lv^*f^!\simeq g^!Lu^*$. 
\end{prop}

\Pf . 
It is enough to check that our map becomes an isomorphism after applying $f_*$.
By the local duality isomorphism we have
$$f_*f^!F\simeq \HHom(f_*\OO_X,F).$$
Note that $f_*\OO_X$ is a vector bundle, so the natural map
$$\HHom(f_*\OO_X,\OO_Y)\otimes F\to \HHom(f_*\OO_X,F)$$
is an isomorphism. It remains to observe that the target of this map is isomorphic
to $f_*(f^!\OO_X\otimes f^*F)$ (by the projection formula).

Given the base change diagram as above, 
we have a natural morphism $v^*f^!\OO_Y\to g^!\OO_U$
corresponding by adjunction to 
$$g_*v^*f^!\OO_Y\simeq u^*f_*f^!\OO_Y\to \OO_U.$$
To check that it is an isomorphism it is enough to compare the push-forwards
by $g$. It remains to use an isomorphism
$$\HHom(g_*\OO_V,\OO_U)\simeq u^*\HHom(f_*\OO_X,\OO_Y)$$
that follows from $g_*\OO_V\simeq u^*f_*\OO_X$.
\ed

\subsection{Sheaves on formal schemes}

We use the definitions of \cite{EGAI} concerning formal schemes. 
{\it All formal schemes considered in
the present paper are assumed to be locally noetherian}. Thus, our formal schemes are locally of
the form $\Spf A$, where $A$ is an adic noetherian ring.
Following \cite{AJL} we say that a morphism $f:X\to Y$ between formal schemes is of
{\it pseudo-finite type} (resp., {\it pseudo-finite}) if there exist ideals of definition $\JJ_X\sub\OO_X$,
$\JJ_Y\sub\OO_Y$, such that $\JJ_Y\OO_X\sub\JJ_X$ and the induced map of schemes
$(X,\OO_X/\JJ_X)\to(Y,\OO_Y/\JJ_Y)$ is of finite type (resp., finite).

By a sheaf on $X$ we always mean a sheaf of $\OO_X$-modules.
We denote by $A(X)$ the abelian category of sheaves on $X$ and by $D(X)$ its (unbounded) derived
category.

Following \cite{AJL} we work with several full subcategories of $A(X)$ and $D(X)$.
First, we consider the categories $A_{qc}(X)\sub A(X)$ (resp., $A_c(X)$) of quasicoherent 
(resp., coherent) sheaves. 
Next, let us recall the definition of torsion sheaves on $X$.
For a coherent ideal sheaf $\JJ\sub\OO_X$ one defines the subfunctor of the identity functor
$$\Ga_{\JJ}\FF=\lim_n\HHom(\OO_X/\JJ^n,\FF).$$
This functor is left exact and depends only on the topology on $\OO_X$ defined by $\JJ$. By
definition, $\JJ$-torsion sheaves are sheaves $\FF$ satisfying $\Ga_{\JJ}\FF=\FF$.
For a formal scheme $X$ {\it torsion sheaves} on $X$
are defined as $\JJ$-torsion sheaves, where $\JJ$ is an ideal of definition for $X$
(recall that on locally noetherian formal schemes ideals of definition exist globally, see 
\cite{EGAI}, Prop.~(10.5.4)).
We denote by $A_{qct}(X)\sub A_{qc}(X)$ the subcategory of torsion quasicoherent sheaves.

Although we are mainly interested in the categories $A_c(X)$ and $A_{qct}(X)$,
it is sometimes convenient to work
with the subcategory $A_{\clim}(X)\sub A_{qc}(X)$ of sheaves that 
can be be presented as direct limits of coherent sheaves. 
Note that $A_{\clim}(X)$ contains both $A_c(X)$ and $A_{qct}(X)$ (see Lemma 5.1.4 of \cite{AJL}).
Also, locally every quasicoherent sheaf is in $A_{\clim}$ (see Cor. 3.1.4 of \cite{AJL}).

For a pair of quasicoherent sheaves $F$ and $G$ on a formal scheme $X$ we set
\begin{equation}\label{hat-Hom-eq}
\widehat{\Hom}(F,G)=\projlim_n\Hom(F,G\otimes\OO_X/\JJ^n),
\end{equation}
where $\JJ$ is an ideal of definition for $X$. It is easy to see that these groups do not depend
on a choice of an ideal of definition. If $i_n:X_n\hra X$ is a closed embedding of the subscheme
corresponding to the ideal sheaf $\JJ^n$ then we have
$$\widehat{\Hom}(F,G)=\projlim_n\Hom(i_n^*F,i_n^*G).$$
Hence, we have a natural composition law for $\widehat{\Hom}$.
Given an element $\a\in\widehat{\Hom}(F,G)$ then for
every morphism $f:Y\to X$ of formal schemes we have the induced element
$f^*\a:\widehat{\Hom}(F,G)$ (use ideals of definition $\JJ_X\in\OO_X$ and $\JJ_Y\in \OO_Y$
such that $\JJ_X\OO_Y\sub\JJ_Y$). 


For $*=c, qc, \clim, qct$ we denote by $D_*(X)$ the full subcategory of complexes with cohomology in
$A_*(X)$. 

Let $R\Ga_{\JJ}$ denote the derived functor to $\Ga_{\JJ}$.
For a complex of sheaves with $\JJ$-torsion cohomology
the natural map $R\Ga_{\JJ}F\to F$ is an isomorphism (see the proof of Proposition 5.2.1(a) of \cite{AJL}). 
In the case when $\JJ$ is an ideal of definition of $X$ we set $R\Ga'_X=R\Ga_{\JJ}$.
For $F\in D_{qc}(X)$ we have $R\Ga'_X(F)\in D_{qct}(X)$. Also, for such $F$ the natural map
$R\Ga'_X(F)\to F$ is an isomorphism iff $F\in D_{qct}(X)$ (see Prop. 5.2.1 of \cite{AJL}).

Since we consider only locally noetherian formal schemes, by Cor. 3.1.2 of \cite{AJL-loc},
the natural map
\begin{equation}\label{torsion-isom-eq}
E\otimes^{\L} R\Ga'_X(F)\to R\Ga'_X(E\otimes^{\L} F)
\end{equation}
is an isomorphism for all $E,F\in D(X)$. 

\begin{lem}\label{tensor-lem} 
Let $X$ be a formal scheme
The subcategories $D_c(X)$, $D_{\clim}(X)$ and $D_{qc}$ are all stable under
the derived tensor product functor $\otimes^{\L}: D(X)\times D(X)\to D(X)$.
Also, if $A\in D(X)$ has $\JJ$-torsion cohomology, where $\JJ\sub\OO_X$ is a coherent ideal sheaf
then the same is true for $A\otimes^{\L} B$ for any $B\in D(X)$.
\end{lem}

\Pf . Using the formula $A\otimes^{\L} B=\lim_{m,n}\tau_{\le m}A\otimes^{\L}\tau_{\le n}B$ and
the spectral sequence we can reduce the first assertion to showing that for $F,G\in A_*(X)$
(where $*=c,\clim,qc$) one has $\TTor_n(F,G)\in A_*(X)$.
In the case $F$ and $G$ are coherent sheaves this is clear by reducing to the affine case,
where the tensor product simply corresponds to the tensor product over the corresponding ring
(see \cite{EGAI}, Prop. (10.10.2.3)). The case $F,G\in A_{\clim}(X)$ follows because 
$\TTor_n$ commutes with direct limits (see 2.5 of \cite{Lipman}). Since locally any quasicoherent
sheaf is in $A_{\clim}(X)$, this also implies the case of $*=qc$. 

As above, the last statement can be reduced to the case when $A$ is a $\JJ$-torsion sheaf.
Now the assertion follows by choosing a $q$-flat resolution of $B$ and using the similar
statement for the underived tensor product of sheaves.
\ed


For a morphism $f:X\to Y$ of formal schemes one has the derived functors of push-forward
$Rf_*:D(X)\to D(Y)$ and pull-back $Lf^*:D(Y)\to D(X)$. In the following proposition we collect
some of their properties (mostly proved in \cite{AJL}).

\begin{prop}\label{formal-pull-push-prop} 
(i) One has $Lf^*D_{\clim}(Y)\sub D_{qc}(X)$. If in addition, $Y$ is affine
then $Lf^*D_{\clim}(Y)\sub D_{\clim}(X)$. On the other hand, if $f$ is flat then 
$f^*D_c(Y)\sub D_c(X)$ and $f^*D_{\clim}(Y)\sub D_{\clim}(X)$.

\noindent
(ii) Let $\JJ\sub\OO_Y$ be a coherent sheaf of ideals, and
let $D_{\JJ}(Y)\sub D(Y)$ (resp., $D_{\JJ\OO_X}(X)\sub D(X)$) denote the subcategory
of complexes with $\JJ$-torsion (resp., $\JJ\OO_X$-torsion) cohomology.
There is an isomorphism of functors
$$Lf^*R\circ \Ga_{\JJ}simeq R\Ga_{\JJ\OO_X}\circ Lf^*.$$ 
Hence, $Lf^*D_{\JJ}(Y)\sub D_{\JJ\OO_X}(X)$.
If $X$ is noetherian then we also have
$$R\Ga_{\JJ}\circ Rf_*\simeq Rf_*R\circ\Ga_{\JJ\OO_Y},$$ 
so in this case 
$Rf_*D_{\JJ\OO_X}(X)\sub D_{\JJ}(Y)$.

\noindent
(iii) Assume that $X$ and $Y$ are noetherian. Then the functor $Rf_*$ is
bounded on $D_{\clim}(X)$.

\noindent
(iv) Assume that $f$ is quasicompact. Then
$Rf_*D_{qct}(X)\sub D_{qct}(Y)$. 

\noindent
(v) Assume that $f$ is quasicompact. Then the functor 
$Rf_*$ commutes with small direct sums in $D_{\clim}(X)$. Also, in this case the functors 
$R^if_*$ commute with (small filtered) direct limits. 

\noindent
(vi) If $f$ is proper then $Rf_*D_c(X)\sub D_c(Y)$. If in addition $X$ and $Y$ are noetherian then
$Rf_*D_{\clim}(X)\sub D_{\clim}(Y)$.

\noindent
(vii) If $f$ is affine then $f_*$ is exact on $D_{\clim}(X)$ and $f_*D_{\clim}(X)\sub D_{qc}(Y)$.
\end{prop}

\Pf . (i) The first two assertions follow from Prop. 3.3.5 of \cite{AJL}. For the last assertion we use
the well-known facts that $f^*A_c(Y)\sub A_c(X)$ and $f^*$ commutes with
direct limits.

\noindent
(ii) See Prop. 5.2.8(b),(d) of \cite{AJL}.

\noindent
(iii) See Prop. 3.4.3(b) of \cite{AJL}. 

\noindent
(iv) Since the assertion is local in $Y$, we can assume that $X$ and $Y$ are noetherian.
Then we can use Prop. 5.2.6 of \cite{AJL}. 

\noindent
(v) Since the assertion is local in $Y$, we can assume that $X$ and $Y$ are noetherian.
Then the first assertion follows from Prop. 3.5.2 of \cite{AJL}. 
The second assertion is checked exactly as in the case of schemes
(see Lemma \ref{lim-push-lem}). 

\noindent
(vi) Since the first assertion is local in $Y$, it is enough to consider the case when $X$ and $Y$ are noetherian, so we can use Prop. 3.5.1 of \cite{AJL}.

\noindent
(vii) We can assume that both $X$ and $Y$ are affine: $X=\Spf(A)$ and $Y=\Spf(B)$. 
Exactness of $f_*$ follows from Lemma 3.4.2 of \cite{AJL}.
By part (v), it is enough to check that $f_*A_c(X)\sub A_{qc}(Y)$.
Every coherent sheaf on $X$ is of the form $M^{\Delta}$ for some finitely generated $A$-module $M$.
It is easy to see that $f_*(M^{\Delta})$ 
is the $\OO_Y$-module associated with $M$ viewed as a $B$-module
via the construction of Prop. 3.2 of \cite{Yekut}. Hence, $f_*M\in A_{\clim}(Y)\sub A_{qc}(Y)$.

We will also use push-forwards with quasi-compact support in a very special case. Instead of developing
the general theory we will give an ad hoc definition for this case.

\begin{defi} (i) Let ${\bf P}$ be some property of morphisms between formal schemes, local over base,
and let $f:X\to Y$ be a morphism of formal schemes. 
We say that $f$ is $ldu-{\bf P}$ ({\it ldu}=locally disjoint union) if 
there exists an open covering $(U_i)$ of $Y$ such that for every $i$,
$f^{-1}(U_i)$ is a disjoint union of formal schemes $V_{i,j}$
such that all the maps $f|_{V_{i,j}}:V_{i,j}\to U_i$ have property ${\bf P}$.
Note that if the property ${\bf P}$ is stable under base changes then the same is true
for the property ``$ldu-{\bf P}$".

\noindent
(ii) Let $f:X\to Y$ be an {\it ldu}-quasi-compact morphism between formal schemes. 
Let us consider the subfunctor
$$f_{\sigma*}:A(X)\to A(Y)$$
in $f_*$, where $f_{\sigma*}(F)\sub f_*(F)$ is the subsheaf of sections with quasicompact support
over $Y$. This functor is left exact and we denote by $Rf_{\sigma*}:D(X)\to D(Y)$ its derived functor.
Note that $Y$ can be covered by open subsets $U\sub Y$
such that $f^{-1}(U)=\sqcup V_j$ where each $f_j=f|_{V_j}:V_j\to U$ is quasicompact, and we have
$$Rf_{\sigma*}(F)|_U\simeq\oplus_j Rf_{j*}(F|_{V_j}).$$
Observe also that the definition of $f_{\sigma*}$ is local in $Y$.
Hence, Proposition \ref{formal-pull-push-prop}(iv) implies that
$Rf_{\sigma*}D_{qct}(X)\sub D_{qct}(Y)$. 
If $f$ is quasicompact then $f_{\sigma*}=f_*$.
In the case when $f$ is {\it ldu}-pseudo-finite we will denote $f_{\sigma*}$ by $f_!$
(note that this functor is exact).
\end{defi}

One can immediately generalize many properties of the push-forward functors for formal
schemes to the functors $Rf_{\sigma*}$ for {\it ldu}-quasicompact morphisms. 
For example, for {\it ldu}-quasicompact maps of formal schemes $f:X\to Y$ and $g:Z\to X$
one has $R(f\circ g)_{\sigma*}\simeq Rf_{\sigma*}\circ Rg_{\sigma*}$. 

We now turn to versions of the projection and base change formulae for formal schemes.

\begin{thm}\label{formal-proj-formula}
Let $f:X\to Y$ be a quasicompact map of formal schemes. 
Then for $F\in D_{\clim}(X)$ and $G\in D_{\clim}(Y)$ the natural map
$$G\otimes^{\L} Rf_*F\to Rf_*(Lf^*G\otimes^{\L} F)$$
is an  isomorphism. The similar result holds if $f$ is {\it ldu}-quasicompact and
$Rf_*$ is replaced with $Rf_{\sigma*}$. 
\end{thm}

\Pf .
The question is local in $Y$, so we can assume that $Y$ is affine and $X$ is noetherian
(replacing $X$ by one of its connected component in the second case).
Note that $Lf^*G\in D_{\clim}(X)$ by Proposition \ref{formal-pull-push-prop}(i). 
Hence, $Lf^*G\otimes^{\L} F\in D_{\clim}(X)$ by Lemma \ref{tensor-lem}.
First, assume that both $F$ and $G$ are bounded above.
Fix $F$. Then both sides respect coproducts as functors in $G$ (use \cite{AJL}, Prop. 3.5.2). 
If $G=\OO_Y$ then the statement is clear. Now use the fact that every object of $A_{\clim}(Y)$
is a quotient of a free module plus boundedness above of both sides (by the way out argument---see 
\cite{Lipman}, (1.11.3.1)). 
The case when $F$ and $G$ are unbounded can be deduced from this as in \cite{Lipman}, (3.9.4),
with $D_{qc}$ replaced by $D_{\clim}$.
\ed

\begin{rem} In the above theorem it is enough to assume that there exists an open covering
$(U_i)$ of $Y$ such that $G|_{U_i}\in D_{\clim}(U_i)$ and $F\in D_{\clim}(f^{-1}(U_i))$
(since the question is local in $Y$). For example, if $f$ is a closed embedding then the projection
formula holds for $F\in D^b_{qc}(X)$ and $G\in D^b_{qc}(Y)$.
\end{rem}




We will use two versions of the flat base change formula for sheaves on formal schemes based on
Proposition 7.2 of \cite{AJL}. In addition we prove a base change formula for flat
adic morphisms (part (i)(b) of the Theorem below).

\begin{thm}\label{formal-base-change}
(i) Let
\begin{diagram}
V &\rTo{v}& X\\
\dTo{g} &&\dTo{f}\\
U &\rTo{u}&Y
\end{diagram}
be a cartesian diagram of formal schemes. 

(a) Assume that $f$ is proper and $u$ is flat.
Then for every $F\in D_{\clim}(X)$ the natural map
$$\theta_F:Lu^*Rf_*F\to Rg_*Lv^*F$$
is an isomorphism.

(b) Now assume that $f$ is flat. In addition assume that $f$ is adic and quasicompact (resp., proper). 
Then the map $\th_F$ is an isomorphism
for every $F\in D_{qct}(X)$ (resp., $F\in D_{\clim}(X)$).

\noindent
(ii) Consider a cartesian diagram as above, where
$f$ is of pseudo-finite type and $u$ is flat.
Assume also that $\II$ is a coherent sheaf of ideals on $U$, such that if $\JJ$ is an ideal of definition
on $Y$ then $\JJ\OO_U+\II$ is an ideal of definition on $U$.
Then for every $F\in D_{qct}(X)$ and every $G\in D_{qc}(U)$ such that $G$ has $\II$-torsion cohomology sheaves, the map induced by the base change map
\begin{equation}\label{base-ch-map}
u^*Rf_*F\otimes^{\L} G\to Rg_*v^*F\otimes^{\L} G
\end{equation}
is an isomorphism.

\noindent
(iii) The assertions of (i) and (ii) also hold if 
$Rf_*$ (resp., $Rg_*$) is replaced with $Rf_{\sigma*}$ (resp., $Rg_{\sigma*}$),
and every property ${\bf P}$ of $f$ is replaced with ``$ldu-{\bf P}$".
\end{thm}

\Pf . Let us observe that all the assertions are local in $Y$ and $U$, so in the proof
we can (and will) assume them to be noetherian.

\noindent
(i) (a)
Using Proposition 7.2(c) of \cite{AJL} we see that $R\Gamma'_U\th_F$ is an isomorphism.
By Prop. 6.2.1 of \cite{AJL} this implies that $\th_F$ itself is an isomorphism when $F$ is a
coherent sheaf. Indeed, in this case the source and the target of $\th_F$ are in $D_c(Y')$ by
Proposition \ref{formal-pull-push-prop}(vi). 
It follows that $\th_F$ is an isomorphism for $F\in A_{\clim}(X)$ (using Proposition 
\ref{formal-pull-push-prop}(v)).
By Proposition \ref{formal-pull-push-prop}(iii), the source and the target of $\th_F$
are bounded functors, so we can finish the proof by the way out argument 
(see \cite{Lipman}, (1.11.3)).

(b) Since the assertion is local in $U$, we can assume it to be affine. 
Then it is enough to check that the map $u_*\theta_F$ on $Y$ is an isomorphism.
Note that $Rf_*F\in D_{qct}(Y)$ by Proposition \ref{formal-pull-push-prop}(iv)
(resp., $Rf_*F\in D_{\clim}(Y)$ by Proposition \ref{formal-pull-push-prop}(vi)).
Therefore, using the projection formula (see Theorem \ref{formal-proj-formula}) 
we can identify the source of $u_*\theta_F$ with
$u_*\OO_U\otimes^{\L} Rf_*F$ and the target of $u_*\theta_F$ with
$Rf_*(v_*\OO_V\otimes^{\L} F)$. It remains to use the isomorphism $v_*\OO_V\simeq f^*u_*\OO_U$
(recall that $f$ is adic) and the projection formula again.

\vspace{2mm}

\noindent
(ii) By Proposition 7.2(b) of \cite{AJL}, the map 
$$R\Gamma'_Uu^*Rf_*F\to R\Gamma'_URg_*v^*F$$
is an isomorphism. Using the isomorphism \eqref{torsion-isom-eq} we see that 
the map \eqref{base-ch-map} becomes an isomorphism after
applying $R\Gamma'_U$. It remains to check that the target and the source of \eqref{base-ch-map}
have torsion cohomology sheaves. Since $G$ has $\II$-torsion cohomology sheaves, it is enough
to check that the cohomology of $u^*Rf_*F$ and $Rg_*v^*F$ are
$\JJ\OO_U$-torsion sheaves. For $u^*Rf_*F$ this follows from the fact that $Rf_*F$ has $\JJ$-torsion
cohomology (see Proposition \ref{formal-pull-push-prop}(ii)).
On the other hand, by the same Proposition, 
to check that the cohomology of $Rg_*v^*F$ are $\JJ\OO_U$-torsion sheaves it is enough to check that $v^*F$ has $\JJ\OO_V$-torsion cohomology. But this immediately follows from
the assumption that $F$ has torsion cohomology on $X$.

\vspace{2mm}

\noindent
(iii) Shrinking $Y$ we can assume that all the connected components $X_i$ of $X$ are noetherian.
Then replacing $X$ by $X_i$ and $V$ by $V_i=U\times_Y X_i$.
we reduce ourselves to the situation considered in (i) and (ii).
\ed

\begin{rem}
One situation where the case (ii) of the above theorem applies is this.
Assume that $X$, $Y$ and $Z$ are formal schemes over a usual scheme $S$,
where $Z$ is flat over $S$. Then for a map $f:X\to Y$ of pseudo-finite type we can consider the
cartesian diagram
\begin{diagram}
X'\times_S Z &\rTo{p_1}& X\\
\dTo{f\times\id_Z}&&\dTo{f}\\
Y\times_S Z &\rTo{p_1}&Y
\end{diagram}
Then for every $F\in D_{qct}(X)$ and $G\in D_{qct}(Z)$ the natural map
$$p_1^*Rf_*F\otimes^{\L} Lp_2^*G\to R(f\times\id_Z)_*p_1^*F\otimes^{\L} Lp_2^*G$$
is an isomorphism. 
Similar assertion holds for $Rf_{\sigma*}$ instead of $Rf_*$ provided
$f$ is of {\it ldu}-pseudo-finite type.
\end{rem}


Finally, we need a version of Proposition \ref{finite-duality-prop}
for a certain class of morphisms between formal schemes. 
For the kind of morphisms we are interested in
the situation is much simpler than in the general duality setup considered in \cite{AJL}
(in particular, in our case the functor $f^!$ can be easily constructed at the level of abelian categories
of sheaves). 

\begin{defi} Let $f:X\to Y$ be a flat morphism of formal schemes.
We say that $f$ is {\it nicely ind-finite} if each connected component
of $X$ is affine over $Y$ and can be presented as the inductive limit of a system of closed formal
subschemes $X_0\sub X_1\sub\ldots\sub X$, where each $X_n$ is flat and finite
over $Y$. We say that $f$ is {\it locally nicely ind-finite} if the above condition holds locally in $Y$.
We say that a formal scheme $X$ over a usual scheme $S$ is 
{\it (locally) nicely ind-finite flat over $S$} if
the morphism $\pi:X\to S$ is such. 
\end{defi}

It is easy to see 
that the property of a flat morphism to be nicely ind-finite (resp., locally nicely ind-finite)
is preserved under arbitrary base changes. Also, if $f$ is locally nicely ind-finite morphism
then it is {\it ldu}-pseudo-finite, so we have the exact functor $f_!$.

\begin{prop}\label{formal-finite-duality-prop} 
(i) Let $f:X\to Y$ be a finite flat morphism of formal schemes.
There exists a natural line bundle $f^!\OO_Y$ on $X$, equipped with a morphism
\begin{equation}\label{f!-O-map}
f_*f^!\OO_Y\to \OO_Y,
\end{equation}
such that for every map from a
usual scheme $Y'$ to $Y$ the pull-back of $f^!\OO_Y$ to $X'=Y'\times_Y X$ is isomorphic
to $(f')^!\OO_{Y'}$, where $f':X'\to Y'$ is the induced finite flat morphism of usual schemes.
Under this isomorphism the pull-back to $Y'$ of
the map \eqref{f!-O-map} gets identified with the canonical map
$f'_*(f')^!\OO_{Y'}\to\OO_{Y'}$. 
The formation of $f^!\OO_Y$ and of the map \eqref{f!-O-map} is also compatible with arbitrary base changes $Y'\to Y$ of formal schemes.

If we define the functor
$$f^!:A_{qc}(Y)\to A_{qc}(X):G\mapsto f^!G=f^!\OO_Y\otimes f^*G$$
then for $F\in A_{qct}(X)$ and $G\in A_{qc}(Y)$ the composed map
$$\Hom(F,f^!G)\to \Hom(f_*F,f_*f^!G)\to \Hom(f_*F,G),$$
where the second arrow is induced by \eqref{f!-O-map}, is an isomorphism.

If $g:Z\to X$ is another finite flat morphism then
$(fg)^!\simeq g^!\circ f^!$ on $A_{qc}(Y)$.

If $i:X'\hra X$ is a closed embedding of formal schemes such that $f'=f\circ i:X'\to Y$
is still flat then we have a natural map
\begin{equation}\label{if!-O-map}
i_*(f')^!\OO_Y\to f^!\OO_Y
\end{equation}
inducing similar canonical maps after an arbitrary base change $Y'\to Y$ with $Y$' a usual scheme.
The push-forward of this map by $f_!$ is compatible with the maps \eqref{f!-O-map} for $f$ and $f'$.
For $G\in A_{qc}(Y)$ the induced map
\begin{equation}\label{closed-dual-eq}
i_*(f')^!G\to \HHom(i_*\OO_{X'},f^!G)
\end{equation}
is an isomorphism.

\noindent
(ii) Let $f:X\to Y$ be a locally nicely ind-finite flat morphism of formal schemes.
Then the functor $f_!:A_{qct}(X)\to A_{qct}(Y)$ admits an exact right adjoint functor
$f^!:A_{qct}(Y)\to A_{qct}(X)$.

In the case when $Y$ is a usual scheme, so that $\OO_Y\in A_{qct}(Y)=\Qcoh(Y)$, the
natural map $f_!f^!\OO_Y\to \OO_Y$ induces an isomorphism
\begin{equation}\label{dual-formal-eq}
f^!\OO_Y\otimes f^*G\wt{\to} f^!G
\end{equation}
for $G\in \Qcoh(Y)$.

\noindent
(iii) Let $f:X\to Y$ be a (globally) nicely ind-finite morphism of formal schemes. Then 
there exists a sheaf $f^!\OO_Y\in A_{\clim}(X)$ (possibly depending on a nicely ind-finite structure
$(X_n)$), 
flat over $Y$, equipped with a morphism $c_f:f_!(f^!\OO_Y)\to\OO_Y$, 
such that we have the induced isomorphism \eqref{dual-formal-eq} for $G\in A_{qct}(Y)$.
The formation of $(f^!\OO_Y,c_f)$ is compatible with arbitrary base changes (where the presentation
as a limit is obtained by the base change).

Let $g:Z\to X$ be another nicely ind-finite flat morphism. Then $fg:Z\to Y$ has a natural nicely ind-finite 
structure and there is an isomorphism
\begin{equation}\label{fg-dual-O-eq}
g^*(f^!\OO_Y)\otimes g^!\OO_X\wt{\to}(fg)^!\OO_Y
\end{equation}
and an isomorphism of functors
\begin{equation}\label{fg-dual-eq}
(fg)^!\simeq g^!\circ f^!:A_{qct}(Y)\to A_{qct}(X).
\end{equation}
Furthermore, for $F\in A_{qct}(Y)$ the following diagram is commutative
\begin{diagram}
g^*f^*F\ot g^*f^!\OO_Y\ot g^!\OO_X&\rTo{}&(fg)^*F\ot g^*f^!\OO_Y\ot g^!\OO_X&\rTo{}&
(fg)^*F\ot (fg)^!\OO_Y\\
\dTo{}&&&&\dTo{}\\
g^*f^!F\ot g^!\OO_X&\rTo{}&g^!f^!F&\rTo{}&(fg)^!F
\end{diagram}
with arrows induced by \eqref{dual-formal-eq}, \eqref{fg-dual-O-eq} and \eqref{fg-dual-eq}.

\noindent
(iv) Let $f:X\to Y$ be a locally nicely ind-finite morphism from a formal scheme 
$X$ to a usual scheme $Y$.
Then for quasicoherent sheaves $F$ on $X$ and $G$ on $Y$ one has a natural isomorphism 
\begin{equation}\label{f!-hom-formal-eq}
\Hom(F\ot f^!\OO_Y,f^*G\ot f^!\OO_Y)\simeq \widehat{\Hom}(F,f^*G),
\end{equation}
local in $Y$ (see \eqref{hat-Hom-eq}).
\end{prop}

\Pf . (i) Let us represent $Y$ as the limit of closed subschemes $i_n:Y_n\hra Y$, and set
$j_n:X_n=f^{-1}(Y_n)\hra X$, $f_n=f|_{X_n}:X_n\to Y_n$. Note that $X_n$ are usual schemes
and $X$ is the inductive limit of $X_n$.
Then every $F\in A_{qct}(X)$ is the inductive limit of the subsheaves $j_{n*}F_n$, where
$F_n=j_n^*\HHom(j_{n*}\OO_{X_n},F)$. Note that each $F_n$ is quasicoherent (see \cite{AJL}, 
Cor. 3.1.6(d)).
Also, by Proposition \ref{formal-pull-push-prop}(v), we have $f_*F=\lim_n i_{n*}f_{n*}F_n$.
Similarly, for $G\in A_{qc}(Y)$ let us denote $G_n=i_n^*\HHom(i_{n*}\OO_{Y_n},G)\in A_{qc}(Y_n)$. 
Then by Proposition \ref{finite-duality-prop}, we have
$$\Hom(f_*F,G)=\projlim\Hom(f_{n*}F_n,G_n)\simeq\projlim\Hom(F_n,f_n^!G_n)\simeq
\projlim\Hom(F_n,f_n^!\OO_{Y_n}\otimes f_n^*G_n).$$
Using the compatibility of $f_{n+1}^!$ with the base change $Y_n\to Y_{n+1}$
(Proposition \ref{finite-duality-prop}) 
we obtain 
$$f_{n+1}^!\OO_{Y_{n+1}}|_{X_n}\simeq f_n^!\OO_{Y_n}.$$
Therefore, 
$$f^!\OO_Y:=\projlim_n j_{n*}f_n^!\OO_{Y_n}$$
is a line bundle on $X$ equipped with isomorphisms $j_n^*f^!\OO_Y\simeq f_n^!\OO_{Y_n}$.
Hence,
\begin{align*}
&\Hom(F,f^!\OO_Y\otimes G)\simeq\projlim\Hom(j_{n*}F_n,f^!\OO_Y\otimes f^*G)\simeq
\projlim\Hom(F_n,j_n^*(f^!\OO_Y)\otimes f_n^*G_n)\\
&\simeq\projlim\Hom(F_n,f_n^!\OO_{Y_n}\otimes f_n^*G_n)
\end{align*}
which gives a natural isomorphism
$$\Hom(f_*F,G)\simeq\Hom(F,f^!G).$$ 
Furthermore, we have 
$$f_*f^!OO_Y\simeq\projlim f_{n*}f_n^!\OO_{Y_n},$$
so the canonical maps $f_{n*}f_n^!\OO_{Y_n}\to \OO_{Y_n}$ give rise to a map
\eqref{f!-O-map}. Is is easy to see that the above isomorphism is induced by \eqref{f!-O-map}
(one can replace $F$ by $j_{n*}F_n$, so the statement reduces to the similar statement for
$f_n$).

Now let us consider a base change diagram
\begin{diagram}
V &\rTo{v}& X\\
\dTo{g} &&\dTo{f}\\
U &\rTo{u}&Y
\end{diagram}
Assume first that $U$ is a usual scheme.
Then using the adjoint pair $(g_*,g^!)$ and the morphism $f_*f^!\OO_Y\to\OO_U$,
as in the proof of Proposition \ref{finite-duality-prop},
we construct a natural morphism 
\begin{equation}\label{v*f!-eq}
v^*f^!\OO_Y\to g^!\OO_U.
\end{equation}
Since the map $U\to Y$ locally factors through some subscheme $Y_n$ it follows easily from our
construction that \eqref{v*f!-eq} is an isomorphism. Next, if $U$ is a formal scheme we can
represent it as the limit $\lim_n U_n$ of closed subschemes and consider the corresponding
closed subschemes $V_n\sub V$, so that $V=\lim_n V_n$, and let $g_n:V_n\to U_n$ be
the induced morphisms. Let also $v_n:V_n\to X$ be the morphisms induced by $v:V\to X$. 
Then the above argument gives a natural isomorphism $v_n^*f^!\OO_Y\simeq g_n^!\OO_{U_n}$
for each $n$. Passing to inverse limits we get that \eqref{v*f!-eq} is an isomorphism
(to represent the left-hand side as a limit we use the fact that $f^!\OO_Y$ is a line bundle on $X$).

If $g:Z\to X$ is another finite flat morphism then to construct an isomorphism of functors
$(f\circ g)^!\simeq g^!\circ f^!$ is equivalent to constructing an isomorphism
$$(f\circ g)^!\OO_Y\simeq g^!\OO_Y\otimes g^*f^!\OO_Y$$
of line bundles on $Z$. By choosing a representation $Y=\lim_n Y_n$ as above and using
the definitions this immediately reduces to the case of schemes,
where we can use Proposition \ref{finite-duality-prop}.

Let $i:X'\to X$ be a closed embedding such that $f'=f\circ i:X'\to Y$ is still flat. The map
\eqref{if!-O-map} is obtained by passing to limit from the similar maps associated with the closed
embeddings $X'\cap X_n\to X_n$ and the finite flat maps $f_n:X_n\to Y_n$.
The proof of the fact that \eqref{closed-dual-eq} is an isomorphism easily reduces
to the case of schemes.

\vspace{2mm}

\noindent
(ii) First, let us replace $Y$ by its affine open subset, and $X$ by its connected component, so that $X$ 
is affine and is the inductive limit of a sequence of closed formal
subschemes $k_n:X_n\hra X$, where each $X_n$ is finite and flat over $Y$.  
Note that in this situation we have $f_!=f_*$.
Set $f_n=f\circ k_n:X_n\to Y$. 
Since $X$ and $X_n$ are affine, we have $X=\Spf(A)$,
$X_n=\Spf(A/I_n)$, where $A=\projlim_n A/I_n$ as a topological ring.
Therefore, any ideal of definition $J\sub A$ contains $I_n$ for some $n$. 
It follows that for every $F\in A_{qct}(X)$
we have $F=\lim_n \HHom(\OO_{X}/\II_n,F)$,
where $\II_n\sub\OO_X$ is the ideal sheaf of $X_n$. 
Set $F_n=k_n^*\HHom(\OO_{X}/\II_n,F)\in A_{qct}(X_n)$.
As in part (i), we have 
$$\Hom(f_*F,G)\simeq\projlim_n\Hom(F_n,f_n^!G).$$
Next, by part (i), for $G\in A_{qct}(Y)$ we have isomorphisms
$$(k_{n,n+1})_*(f_n^!G)\simeq\HHom((k_{n,n+1})_*\OO_{X_n},f_{n+1}^!G),$$
where $k_{n,n+1}:X_n\hra X_{n+1}$ is the natural embedding.
Therefore, we have natural maps $k_{n*}f_n^!G\to (k_{n+1})_*f_{n+1}^!G$,
so that setting 
\begin{equation}\label{f!-G-eq}
f^!G:=\lim_n k_{n*}f_n^!G
\end{equation}
we have $k_{n*}f_n^!G\simeq \HHom(\OO_X/\II_n,f^!G)$.
It follows that
$$\projlim_n\Hom(F_n,f_n^!G)\simeq\projlim_n\Hom(k_{n*}F_n,f^!G)\simeq\Hom(F,f^!G),$$
so $G\to f^!G$ is right adjoint to $f_*$.

Next, assume that $Y$ is covered by open subsets $Y_i$ such that the morphisms 
$f_{Y_i}:X_i=f^{-1}(Y_i)\to Y_i$ have the above structure. Then for each $i$ we have
the exact functor $f_{Y_i}^!:A_{qct}(Y_i)\to A_{qct}(X_i)$, right adjoint to $(f_{Y_i})_!$. 
Furthermore, these functors
are compatible with restrictions to open subsets, so we can glue them into an exact functor
$f^!:A_{qct}(Y)\to A_{qct}(X)$, right adjoint to $f_!$.

The isomorphism \eqref{dual-formal-eq} in the case when $Y$ is a usual scheme will follow from
part (iii).

\vspace{2mm}

\noindent
(iii) Let us again consider the situation from the beginning of the proof of (ii). We have natural maps
$$(k_{n,n+1})_*(f_n^!\OO_Y)\to\HHom((k_{n,n+1})_*\OO_{X_n},f_{n+1}^!\OO_Y)\to f_{n+1}^!\OO_Y,$$
so the sheaves $(k_{n*}f_n^!\OO_Y)$ form an inductive system in $A_c(X)$, and we can set 
\begin{equation}\label{f!-O-eq}
f^!\OO_Y=\lim_n k_{n*}f_n^!\OO_Y.
\end{equation}
Note that in the case when $Y$ is a usual scheme we have $\OO_Y\in A_{qct}(Y)$, 
and this definition agrees with \eqref{f!-G-eq}.
Each $k_{n*}f_n^!\OO_Y$ is flat over $Y$, hence so is $f^!\OO_Y$. 
Now using part (i) and \eqref{f!-G-eq} we obtain for $G\in A_{qct}(Y)$
$$f^!\OO_Y\otimes f^*G\simeq\lim_n k_{n*}(f_n^!\OO_Y\otimes f_n^*G)\simeq\lim_n k_{n*}f_n^!G=f^!G$$
as claimed. Note also that
$f_*f^!\OO_Y=\lim_n f_{n*}f_n^!\OO_Y$ and we have a system of compatible 
morphisms $f_{n*}f_n^!\OO_Y\to\OO_Y$ (see part (i)), 
so we get a morphism $f_*f^!\OO_Y\to\OO_Y$. Hence, by the projection formula,
we have a natural map
$$f_*(f^!\OO_Y\ot f^*G)\simeq f_*f^!\OO_Y\ot G\to G,$$
hence, by adjunction, a map \eqref{dual-formal-eq}. One can easily check that it coincides with
the above isomorphism by reducing to the case of a finite morphism.

The fact that the formation of $f^!\OO_Y$ is compatible with any base change $Y'\to Y$
follows from a similar result for finite morphisms (see part (i)).

If $g:Z\to Y$ is another nicely ind-finite morphism then replacing $Z$ by its connected component
we can assume that $Z=\lim_m Z_m$, where $g_m=g|_{Z_m}:Z_m\to X$ are finite and flat.
Let us set $Z_{m,n}=g_m^{-1}(X_n)\sub Z_m$. Then the induced map $g_{m,n}:Z_{m,n}\to X_n$
is finite and flat, hence, $Z_{m,n}$ is finite and flat over $Y$. Thus, $Z=\lim_{m,n}Z_{m,n}$ 
is a nicely ind-finite structure on $fg:Z\to Y$.
The isomorphism $(f g)^!\simeq g^!\circ f^!$ follows by passing to adjoint functors from
$(f g)_!\simeq f_!\circ g_!$.

To construct a map \eqref{fg-dual-O-eq}, we note that it should correspond by adjunction
to a map
$$(fg)_!(g^*f^!\OO_Y\ot g^!\OO_X)\simeq f_!(f^!\OO_Y\ot g_!g^!\OO_X)\to\OO_Y,$$
where we used the projection formula. We clearly have such a map induced by the maps
$g_!g^!\OO_X\to\OO_X$ and $f_!f^!\OO_Y\to\OO_Y$. To see that it is an isomorphism
we apply the definitions \eqref{f!-O-eq} using the nicely ind-finite structures on $f$, $g$ and $fg$.
Namely, let $l_m:Z_m\to Z$ and $l_{m,n}:Z_{m,n}\to Z$ denote the natural closed embeddings. 
Then we have 
\begin{align*}
&g^*f^!\OO_Y\ot g^!\OO_X\stackrel{(1)}{\simeq} g^*f^!\OO_Y\ot \lim_m l_{m*}g_m^!\OO_X
\stackrel{(2)}{\simeq}\lim_m l_{m*}(g_m^*f^!\OO_Y\ot g_m^!\OO_X)\stackrel{(3)}{\simeq}\\
&\lim_{m,n}l_{m,n*}(g_{m,n}^*f_n^!\OO_Y\ot g_{m,n}^!\OO_X)\stackrel{(4)}{\simeq}
\lim_{m,n}l_{m,n*}(f_ng_{m,n})^!\OO_Y\stackrel{(5)}{\simeq} (fg)^!\OO_Y,
\end{align*} 
where (1) and (5) are given by \eqref{f!-O-eq}, (2) follows from the projection formula,
(3) uses \eqref{f!-O-eq} and the compatibility of $g_m^!\OO_X$ with the base change (see part (i)),
and (4) uses an isomorphism \eqref{fg-dual-O-eq} for finite morphisms that follows from part (i).

The proof of the commutativity of the diagram is tedious but straightforward: since the target vertex
is $(fg)^!F$, we can use adjunction to rewrite the commutative diagram applying $(fg)_!=f_!g_!$ to
other vertices.
\vspace{2mm}

\noindent
(iv) It is enough to consider the situation when $X=\lim X_n$, where each $X_n$ is flat and finite
over $Y$. Let us use the same notation as in part (ii).
The natural map is obtained as follows: an element $\a\in\widehat{\Hom}(F,f^*G)$
induces a map $k_n^*F\to f_n^*G$ for each $n$, hence, we get a compatible system of
maps 
$$F\ot k_{n*}f_n^!\OO_Y\to f^*G\ot k_{n*}f_n^!\OO_Y.$$
Taking the inductive limits and using the formula \eqref{f!-O-eq}
we get a map $F\ot f^!\OO_Y\to f^*G\ot f^!\OO_Y$.
To see that it is an isomorphism we use \eqref{dual-formal-eq}:
$$\Hom(F\ot f^!\OO_Y, f^*G\ot f^!\OO_Y)\simeq \projlim(k_n^*(k_n^*F\ot f_n^!\OO_Y, f^!G)
\simeq\projlim(k_n^*F\ot f_n^!\OO_Y, f_n^!G).$$
By Proposition \ref{finite-duality-prop}, this is isomorphic to 
$$\projlim(k_n^*F,f_n^*G)\simeq\widehat{\Hom}(F,f^*G).$$
\ed

\begin{ex}
Let $k$ be a field. A formal scheme $X$ is nicely ind-finite over $\Spec(k)$ 
if and only it is {\it ldu}-pseudo-finite over $k$, i.e., every connected component of $X$ is of the form 
$\Spf A$, where $A$ is a noetherian adic
$k$-algebra, such that $A/J$ is finite-dimesional over $k$ for an ideal of definition $J$.
Thus, if $G$ is a formal $k$-group (see the beginning of section \ref{kerrep-sec} for our
conventions on formal $k$-groups) then $G$ is nicely ind-finite over $\Spec(k)$.
More generally, if $X$ is a formal scheme over $k$, and $f:E\to X$ is a $G$-torsor
then $f$ is locally nicely ind-finite (this immediately reduces to the case of a trivial $G$-torsor).
\end{ex}

Let us point out the following corollary from Proposition \ref{formal-finite-duality-prop}.
For a formal scheme $X$, locally nicely ind-finite flat over a usual scheme $S$, we set
$\om_{X/S}=\om_{\pi}:=\pi^!\OO_S\in A_{qct}(X)$, where $\pi:X\to S$ is the structure morphism.

\begin{cor}\label{formal-duality-cor} 
Let $S$ be a scheme, $X$ and $Y$ formal schemes over $S$, where $Y$ is 
locally nicely ind-finite flat over $S$. Then for $F\in D_{qct}(X)$ there is a natural isomorphism on
$X\times_S Y$
\begin{equation}\label{dual-XY-eq}
\a_{X,Y}:p_X^*F\ot p_Y^*\om_{Y/S}\wt{\to} p_X^!F,
\end{equation}
where $p_X$ and $p_Y$ are the projections from $X\times_S Y$ to $X$ and $Y$, respectively.
Assume that $T$ is another scheme, $Y\to T$ is a morphism, and 
$Z$ is a formal scheme, nicely ind-finite flat over $T$. Let us form the fibered product
$X\times_S Y\times_T Z$, and let $p_X$, $p_Y$, $p_{XY}$, etc., denote the projections
from this product to the partial products of factors. Then $Y\times_T Z$ is locally nicely ind-finite flat
over $T$, and the following diagram is commutative
\begin{diagram}
p_X^*F \ot p_Y^*\om_{Y/S}\ot p_Z^*\om_{Z/T}&\rTo{p_{YZ}^*\a_{Y,Z}}&
p_X^*F\ot p_{YZ}^*(p^{YZ}_Y)^!\om_{Y/S}&\rTo{}&p_X^*F \ot p_{YZ}^*\om_{Y\times_T Z/S}\\
\dTo{\a_{X,Y}}&&&&\dTo{\a_{X,YZ}}\\
p_{XY}^*((p^{XY}_X)^!F) \ot p_Z^*\om_{Z/T}&\rTo{\a_{XY,Z}}&
p_{XY}^!(p^{XY}_X)^!F&\rTo{}&p_X^!F
\end{diagram}
where $p^{XY}_X:X\times_S Y\to X$ (resp., $p^{YZ}_Y:Y\times_T Z\to Y$) are the natural projections,
the marked arrows are induced by \eqref{dual-XY-eq}, the remaining arrows are obtained from
the natural isomorphism $(fg)^!\simeq g^!\circ f^!$. 
\end{cor}

\Pf . By adjunction, a map \eqref{dual-XY-eq} should correspond to
morphism on $X$
$$p_{X!}(p_X^*F\ot p_Y^*\om_{Y/S})\simeq F\ot p_{X!}p_Y^*\om_{Y/S}\to F$$
(where we used the projection formula).
Applying Theorem \ref{formal-base-change}(ii),(iii) to the cartesian diagram
\begin{diagram}
X\times_S Y &\rTo{p_Y}& Y\\
\dTo{p_X}&&\dTo{\pi_Y}\\
X&\rTo{\pi_X}& S
\end{diagram}
we obtain an isomorphism $F\ot p_{X!}p_Y^*\om_{Y/S}\simeq F\ot \pi_X^*(\pi_{Y!})\om_{Y/S}$.
Now the desired morphism is induced by the canonical
map $\pi_{Y!}\om_{Y/S}=\pi_{Y!}\pi_Y^!\OO_S\to \OO_S$.
To check that \eqref{dual-XY-eq} is an isomorphism we can argue locally in $S$, so we can
assume that $\pi_Y:Y\to S$ is (globally) nicely ind-finite. Then we can apply the compatibility of the formation of $\om_{Y/S}=\pi_Y^!\OO_S$ with the base change, which gives an isomorphism
$p_Y^*\om_{Y/S}\simeq p_X^!\OO_X$. Now \eqref{dual-XY-eq} follows from \eqref{dual-formal-eq}.

In the case when $Z$ is nicely ind-finite over $T$, then the projection $Y\times_T Z\to S$
is the composition of $p^{YZ}_Y$ with $\pi_Y$, so it is nicely ind-finite. 
The desired commutativity of the diagram follows from the commutativity of the diagram
in Proposition \ref{formal-finite-duality-prop}(iii) applied to the composition of $p_X^{XY}$ with $p_{XY}$
(note that both these projections are nicely ind-finite).
\ed

\subsection{Biextensions}\label{biext-sec}

Let $\SS$ be a site, and let $\Sh(\SS)$ denote the category of sheaves of abelian groups over $\SS$.
 We refer to sec.~1.4 of \cite{SGA4}  for basic definitions concerning 
(strictly commutative) Picard stacks over $\SS$. Roughly speaking, these are sheaves
of categories $\PP$ over $\SS$ equipped with a commutative group law $\PP\times\PP\to \PP$
satisfying appropriate axioms. If $K=[K^{-1}\stackrel{d}{\to} K^0]$ 
is a complex over $\Sh(\SS)$, concentrated in degrees $-1$ and $0$, 
then one has the corresponding Picard stack $\ch(K)$ with objects given
by sections of $K^0$, where a morphism between $x$ and $y$ in $K^0(U)$ is an element
$f\in K^{-1}(U)$ such that $df=y-x$. In fact, every Picard stack can be represented in this way
(see \cite{SGA4}, 1.4.13).
For example, if $G$ is a sheaf of abelian groups then the Picard stack of 
$G$-torsors is equivalent to $\ch(G[1])$
(see \cite{SGA4}, 1.4.21).

For a pair of Picard stacks $\PP_1,\PP_2$ there is a natural Picard stack $HOM(\PP_1,\PP_2)$
of {\it homofunctors} $\PP_1\to\PP_2$ (the corresponding term in \cite{SGA4} ,1.4 is ``foncteur  additif"). 
For complexes of sheaves $K_1$ and $K_2$, concentrated
in degrees $-1$ and $0$, one has
\begin{equation}\label{Pic-Hom-eq}
HOM(\ch(K_1),\ch(K_2))\simeq\ch\tau_{\le 0}R\und{\Hom}(K_1,K_2),
\end{equation}
where in the right-hand side we view $K_1$ and $K_2$ as objects of the derived category
$D^b(\Sh(\SS))$ (see \cite{SGA4}(1.4.18.1)).

\begin{lem}\label{Pic-ext-lem} (\cite{SGA4}, 1.4.23)
Let $K=[K^{-1}\to K^0]$ be a complex over $\Sh(\SS)$, and let $G\in\Sh(\SS)$. Then
the Picard stack of extensions
$\EExt(K,G)$ (i.e., of extensions of $K^0$ by $G$, trivialized over $K^{-1}$)
is equivalent to $HOM(\ch(K),\ch(G[1]))$. 
\end{lem}

For a triple $\PP_1,\PP_2,\PP_3$ of Picard stacks there is a Picard stack
$HOM(\PP_1,\PP_2;\PP_3)$ of biadditive functors $\PP_1\times\PP_2\to \PP_3$ 
(see \cite{SGA4}, 1.4.8). In terms of complexes of sheaves $K_1,K_2,K_3$ as above, we have
\begin{equation}\label{bi-Hom-eq}
HOM(\ch(K_1),\ch(K_2);\ch(K_3))\simeq\ch\tau_{\le 0}R\und{\Hom}(\tau_{\ge -1}(K_1\ot^{\L} K_2),K_3)
\simeq\ch\tau_{\le 0}R\und{\Hom}(K_1\ot^{\L} K_2,K_3)
\end{equation}
(using \cite{SGA4}, (1.4.8.1), 1.4.20 and \eqref{Pic-Hom-eq}).

\begin{lem}
One has a natural equivalence
\begin{equation}\label{biadd-equiv}
HOM(\PP_1,\PP_2;\PP_3)\simeq HOM(\PP_1,HOM(\PP_2,\PP_3)).
\end{equation}
\end{lem}

\Pf . This can be deduced directly from the definitions. Alternatively, 
we can realize each Picard stack $\PP_i$ as $\ch(K_i)$ for some complexes $K_i$ ($i=1,2,3$)
concentrated in degrees $-1$ and $0$. Using 
\eqref{Pic-Hom-eq} and \eqref{bi-Hom-eq}
we obtain similar representations for both sides of \eqref{biadd-equiv}. It remains to use
the isomorphism
$$R\underline{\Hom}(K_1\otimes^{\L} K_2,K_3)\simeq
R\underline{\Hom}(K_1,R\underline{\Hom}(K_2,K_3))\simeq
R\underline{\Hom}(K_1,\tau_{\le 0}R\underline{\Hom}(K_2,K_3)).$$
\ed

\begin{prop}\label{biext-prop1} 
Let $P$, $Q$ and $G$ be three sheaves of commutative groups on some site.
Then the category of biextensions of $P\times Q$ by $G$ is equivalent to the category of
homofunctors $P\to\EExt(Q,G)$.
\end{prop}

\Pf . Let $BIEXT(P,Q;G)$ denote the Picard stack of biextensions of $P\times Q$ by $G$.
Then we have a natural functor
\begin{equation}\label{biext-hom-eq}
BIEXT(P,Q;G)\to HOM(\ch(P),\ch(Q);\ch(G[1]))
\end{equation}
that associates with a biextension the corresponding map from pairs of sections $(p,q)$ of $P\times Q$ to $G$-torsors. Using \eqref{biadd-equiv} and Lemma \ref{Pic-ext-lem}
we see that the right-hand side is equivalent
to the Picard stack of homofunctors from $P$ to $\EExt(Q,G)$.
Thus, it suffices to prove that \eqref{biext-hom-eq} is an equivalence. 
This can be checked directly using the definition of the biextension
as in \cite{SGA7-VII}, 2.0. Here is a different way to check this.
By \eqref{bi-Hom-eq}, we have 
$$HOM(\ch(P),\ch(Q);\ch(G[1]))\simeq\ch\tau_{\le 0}R\underline{\Hom}(P\otimes^{\L} Q,G[1]).$$
Therefore, the sheaf of automorphisms of the neutral object of this Picard stack is 
$\underline{\Hom}(P\otimes^{\L} Q, G)$ while the sheaf associated with the isomorphism classes of
objects is $\underline{\Ext}^1(P\otimes^{\L} Q,G)$. But these are exactly the same sheaves as one
gets for the Picard stack of biextensions (see \cite{SGA7-VII}, 2.5.4 and 3.6.5).
\ed

The following result follows immediately from \cite{SGA7-VIII} (1.1.6) and 1.5 
(it can also be derived from
the above Proposition by considering the subcategory of locally trivial extensions in $\EExt(Q,G)$).

\begin{prop}\label{biext-prop2} 
For $P,Q,G$ as above there is a fully faithful functor
$$\EExt(P,\und{\Hom}(Q,G))\to BIEXT(P,Q;G).$$
If $\underline{\Ext}^1(Q,G)=0$ then this functor is an equivalence.
\end{prop}

The condition $\underline{\Ext}^1(Q,G)=0$ is satisfied in the following important case.

\begin{lem}\label{finite-ext-van-lem} (\cite{SGA7-VIII}, 3.3.1)
Let $G$ be a finite flat group scheme over $S$
Then $\und{\Ext}^1(G,\G_m)=0$.
\end{lem}

Now let $k$ be a field of characteristic zero. Then similar vanishing holds for 
affine commutative algebraic groups over $k$ and for formal $k$-groups 
(for our assumptions on formal $k$-groups see the beginning of section \ref{kerrep-sec}).

\begin{lem}\label{formal-ext-van-lem}
(i) Let $G$ be an affine commutative algebraic group over $k$ or a formal $k$-group.
Then
$\und{\Ext}^1(G,\G_m)=0$.

\noindent
(ii) If $0\to G_1\to G_2\to G_3\to 0$ is an exact sequence of formal $k$-groups then
there exists a finite field extension $k\sub k'$ and a section $(G_3)_{k'}\to (G_2)_{k'}$ (not
required to be a group homomorphism).
Hence, $G_2\to G_3$ is a $G_1$-torsor in flat topology. 
\end{lem}

\Pf . (i) For affine commutative algebraic groups this follows immediately from the Cartier
duality since $\G_m$ is dual to $\Z$. Now assume $G$ is a formal $k$-group.
Let $G^0$ be the component of identity in $G$. Passing to a finite extension of $k$ we can
assume that $G/G^0$ is a finitely generated discrete abelian group, so $\und{\Ext}^1(G/G^0,\G_m)=0$.
Now If we have a group extension 
$$1\to G_m\to H\to G^0\to 1$$
then we can pass to the corresponding extension of Lie algebras
$$0\to k\to LH\to LG^0\to 0.$$
Since the characteristic of $k$ is zero and all the Lie brackets here are trivial,
we have a splitting $LG^0\to LH$ that induces a splitting $G^0\to H$

\noindent
(ii) Since all our formal $k$-groups are assumed to be commutative,
the extension splits over $G_3^0$ (by considering the extension of Lie algebras). Thus, it is enough
to consider the case when $G_3$ is \'etale.
Passing to a finite extension of $k$ we can assume that $G_3$ is a finitely
generated discrete abelian group, in which case the statement is clear. 
\ed

\begin{rem} The vanishing in Lemma \ref{formal-ext-van-lem}(i)
also holds for formal groups over fields of finite characteristic---one can reduce this
to the case of finite group schemes.
\end{rem}

\section{Kernel algebras}
\label{keralg-sec}

In this section we present a formalism generalizing the techniques of Fourier-Mukai transforms for
D-algebras developed in \cite{PR}. We consider more general ``nonlocal" algebras and
generalize the main construction of the ``circle product" from \cite{PR} by making use
of appropriate resolutions. The main result of this section
is Theorem \ref{PR-thm} generalizing Thm.~6.5 of \cite{PR}.
Throughout this section we work over a fixed base scheme $S$, and all our schemes
are assumed to be $S$-schemes. 
We denote simply by $X\times Y$ the fibered product over $S$.
{\it All schemes in this section are assumed to be semi-separated}.

\subsection{Convolution for quasicoherent sheaves}\label{circle-sec}

Let us start by recalling the basic properties of the convolution operation for quasicoherent sheaves
(called the {\it circle product} in \cite{PR}). 
Given $F\in D_{qc}(X\times Y)$ and $G\in D_{qc}(Y\times Z)$ we set
$$F\circ _Y G=p_{13*}(p_{12}^*F\ot p_{23}^*G),$$
where $p_{ij}$ are the projections from $X\times Y\times Z$ to the double products of factors.

\begin{lem}\label{circle-prop-lem}
(i) Assume that in each pair of $S$-schemes
$(X,Z)$ and $(Y,T)$ at least one is flat over $S$. Then for $F\in D_{qc}(X\times Y)$,
$G\in D_{qc}(Y\times Z)$ and $H\in D_{qc}(Z\times T)$ one has
$$(F\circ_Y G)\circ_Z H\simeq F\circ_Y (G\circ_Z H).$$

\noindent (ii) Assume that $X$ is flat over $S$.
Then for $F\in D_{qc}(X\times Y)$ and $G\in D_{qc}(Y)$ one has
$$F\circ_Y\De_*G\simeq F\otimes p_2^*G.$$
\end{lem}

The proof is a simple application of the projection and base change formulae (cf.
\cite{Mukai}, Prop. 1.3).

\begin{lem}\label{circle-push-pull-lem}
(i) Assume that we are given a morphism $g:Y\to Y'$ and let $X$ and $Z$ be schemes.
If $g$ is flat or $Z$ is flat over $S$ then for $F\in D_{qc}(X\times Y)$ and 
$G'\in D_{qc}(Y'\times Z)$ we have
a natural isomorphism
$$\left((\id_X\times g)_*F\right)\circ_{Y'}G\simeq F\circ_Y (g\times\id_Z)^*G.$$
Similarly, if $g$ is flat or $X$ is flat over $S$ 
then for $F'\in D_{qc}(X\times Y')$ and $G\in D_{qc}(Y\times Z)$ one
has
$$F'\circ_{Y'}(g\times\id_Z)_*G\simeq \left((\id_X\times g)^*F'\right)\circ_Y G.$$

\noindent
(ii) Assume that we are given morphisms $f:X\to X'$ and $h:Z\to Z'$.
If $f$ is flat or $Z$ is flat over $S$ then for $F\in D_{qc}(X\times Y)$ and $G\in D_{qc}(Y\times Z)$ we have a natural isomorphism
$$(f\times \id_Z)_*(F\circ_Y G)\simeq \left((f\times\id_Y)_*F\right)\circ_Y G.$$
Similarly, if $h$ is flat or $X$ is flat over $S$ then one has
$$(\id_X\times h)_*(F\circ_Y G)\simeq F\circ_Y (\id_Y\times h)_*G.$$

\noindent
(iii) Assume that we are given morphisms $f:X\to X'$ and $h:Z\to Z'$, and let $Y$ be a scheme.
Assume that either $f$ and $h$ are flat or $Y$ is flat over $S$. Then 
for $F\in D_{qc}(X'\times Y)$ and $G\in D_{qc}(Y\times Z')$ one has a natural isomorphism
\begin{equation}\label{pull-circle-eq}
(f\times h)^*(F\circ_Y G)\simeq \left((f\times\id_Y)^*F\right)\circ_Y \left((\id_Y\times h)^*G'\right).
\end{equation}

\noindent
(iv) Assume that we are given morphisms $f:X\to X'$, $g:Y\to Y'$, and
let $Z$ be a scheme. Assume that either $f$ and $g$ are flat or $Z$ is flat over $S$. Then for
$F\in D_{qc}(X\times Y)$ and $G\in D_{qc}(Y'\times Z)$ one has
a natural isomorphism
\begin{equation}\label{push-circle-eq}
\left((f\times g)_*F\right)\circ_{Y'} G\simeq (f\times \id_Z)_*\left(F\circ_Y (g\times \id_Z)^*G\right).
\end{equation}
\end{lem}

\Pf .
The proofs of (i), (ii) and (iii) are easy applications of the projection and base change formulae.
(iv) follows from (i) and (ii).
\ed

Assume we are given a triple of morphisms $f:X\to X'$, $g:Y\to Y'$, $h:Z\to Z'$, where either all these morphisms are flat, or $X$ and $Z'$ are flat over $S$, or $X'$ and $Z$ are flat over $S$. Then
we can define for $F\in D_{qc}(X\times Y)$ and $G\in D_{qc}(Y\times Z)$
a natural morphism
\begin{equation}\label{push-circle-mor}
(f\times g)_*F\circ_{Y'}(g\times h)_*G\to (f\times h)_*(F\circ_Y G).
\end{equation}
Indeed, say, assume that $X$ and $Z'$ are flat over $S$. 
Then using \eqref{push-circle-eq}
we get an isomorphism
$$(f\times g)_*F\circ_{Y'}(g\times h)_*G\simeq (f\times\id_Y)_*
\left(F\circ_Y (g\times\id_Z)^*(g\times h)_*G\right).$$
Now using the adjunction morphism 
$$(g\times\id_Z)^*(g\times h)_*G\simeq (g\times\id_Z)^*(g\times\id_Z)_*(\id_Y\times h)_*G\to
(\id_Y\times h)_*G$$
we get a morphism from $(f\times g)_*F\circ_{Y'}(g\times h)_*G$ to
$$(f\times\id_Y)_*\left(F\circ_Y (\id_Y\times h)_*G\right)\simeq (f\times h)_*(F\circ_Y G),$$
where the last isomorphism follows from Lemma \ref{circle-push-pull-lem}(ii).
We leave for the reader to formulate and
prove the associativity of this construction, say, under the assumption
that all schemes involved are flat over $S$.

\begin{defi} An object $K\in D_{qc}(X\times Y)$ is {\it $\circ$-flat over }$Y$  (or {\it $\circ_Y$-flat}) if
for every quasicoherent sheaf $F$ on $Y$ 
the object $K\circ_Y F\in D_{qc}(X)$ is a sheaf. In the case $X=S$ we have $K\in D_{qc}(Y)$
and the condition of $\circ_Y$-flatness means that the functor $F\mapsto \pi_*(F\ot K)$ on 
$D_{qc}$ is $t$-exact, where $\pi:Y\to S$ is the projection.
\end{defi}

Note that the condition of $\circ$-flatness over $Y$ is local over $X$. 

\begin{lem}\label{flat-lem}
(i) Let $f:X\to X'$ and $\tilde{g}:\tilde{Y}\to Y$ be affine morphisms and 
let $\tilde{f}:\tilde{X}\to X$ and $g:Y\to Y'$ be flat morphisms.
Then $(f\times g)_*$ sends $\circ_Y$-flat objects to $\circ_{Y'}$-flat objects. Assume in addition
that either $\tilde{g}$ is flat or $\tilde{X}$ is flat over $S$. Then
$(\tilde{f}\times\tilde{g})^*$ sends $\circ_Y$-flat objects to $\circ_{\tilde{Y}}$-flat objects.

\noindent
(ii) If $K\in D_{qc}(X\times Y)$ is $\circ_Y$-flat then
$K$ is a quasicoherent sheaf on $X\times Y$, flat over $Y$. 
If $Y$ is affine over $S$ then the converse is also true.

\noindent
(iii) Let $K$ be a $\circ_Y$-flat quasicoherent sheaf on $X\times Y$, and let $Z$ be a scheme.
Assume that either $X$ or $Z$ is flat over $S$. Then for any
quasicoherent sheaf $L$ on $Y\times Z$
the object $K\circ_Y L\in D_{qc}(X\times Z)$ is a sheaf.

\noindent
(iv) Let $K\in \Qcoh(X\times Y)$ be $\circ_Y$-flat and $L\in \Qcoh(Y\times Z)$ be $\circ_Z$-flat.
Assume that either $X$ or $Z$ is flat over $S$.
Then $K\circ_Y L\in D_{qc}(X\times Z)$ is also $\circ_Z$-flat.

\noindent
(v) Assume again that either $X$ or $Z$ is flat over $S$.
If $K\in \Qcoh(X\times Y)$ is $\circ_Y$-flat and
$L$ is a quasicoherent sheaf on $Y\times Z$, flat over $Z$, then 
$K\circ_Y L\in \Qcoh(X\times Z)$ is also flat over $Z$.

\noindent
(vi) Let $q:K^{\bullet}_1\to K^{\bullet}_2$ is a quasi-isomorphism of bounded complexes of quasicoherent sheaves on $X\times Y$ with $\circ_Y$-flat terms.
Then for every $L\in \Qcoh(Y\times Z)$ the induced morphism of complexes
$K^{\bullet}_1\circ_Y L\to K^{\bullet}_2\circ_Y L$ is a quasi-isomorphism.
\end{lem}

\Pf . (i) This follows from the isomorphisms of Lemma \ref{circle-push-pull-lem}.

\noindent
(ii) Without loss of generality we can assume $S$ to be affine.
Assume first that $Y$ is also affine. Then the functor $p_{1*}:D_{qc}(X\times Y)\to D_{qc}(X)$
is $t$-exact and conservative. Hence, the condition of $\circ_Y$-flatness is equivalent
to the condition that the functor $F\mapsto K\otimes p_2^*F:D_{qc}(Y)\to D_{qc}(X\times Y)$
is $t$-exact. Taking $F=\OO_Y$ we immediately see that $K$ itself is a sheaf. Furthermore,
this exactness condition is easily seen to be equivalent to flatness of $K$ over $Y$.

Now for arbitrary $Y$, using part (i), we see that $K|_{X\times U}$ is $\circ_U$-flat for any open affine
subset $U\sub Y$. As we have just shown this implies that $K|_{X\times U}$ is a sheaf, flat over $U$.
Hence, $K$ itself is a sheaf, flat over $Y$.

\noindent
(iii) The problem is local in $Z$, so we can assume that $Z$ is affine.
Then it is enough to check that $p^{XZ}_{X*}(K\circ_Y L)\in D_{qc}(X)$ is a sheaf, 
where $p^{XZ}_X:X\times Z\to X$ is the projection.
Since $X$ is flat over $S$, by Lemma \ref{circle-push-pull-lem}, one has
$$p^{XZ}_{X*}(K\circ_Y L)\simeq K\circ_Y(p^{YZ}_{Y*}L),$$
where $p^{YZ}:Y\times Z\to Y$ is the projection.
Now we observe that $p^{YZ}_{Y*}L$ is a sheaf on $Y$ and use the fact that $K$ is $\circ_{Y}$-flat.

\noindent
(iv) This follows from the associativity of the convolution in this case 
(see Lemma \ref{circle-prop-lem}(i)).

\noindent
(v) Without loss of generality we can assume $S$ and $Z$ to be affine. Now the statement follows from 
(ii) and (iv).

\noindent
(vi) It is enough to check that if a bounded complex $K^{\bullet}$ of quasicoherent
sheaves on $X\times Y$ with $\circ_Y$-flat terms represents an object 
$F\in D_{qc}(X\times Y)$ then the complex $K^{\bullet}\circ_Y L$ represents 
$F\circ_Y L\in D_{qc}(X\times Z)$. This can be easily proved by induction in the length of $K^{\bullet}$.
\ed

\begin{lem}\label{aff-flat-lem} 
(i) Let $f:Z\to X$ be an affine morphism and $g:Z\to Y$ a morphism of schemes. 
Then for any quasicoherent sheaf $F$ on $Z$, flat over $Y$, the 
quasicoherent sheaf $(f,g)_*F$ on $X\times Y$ is
$\circ$-flat over $Y$.

\noindent
(ii) The same assertion is true for the sheaf $(f,g)_{\sigma*}F$ 
if $X$ and $Y$ are schemes, $Z$ is a formal scheme, 
$f$ is {\it ldu}-affine, and $F$ is a sheaf in $A_{qct}(Z)$, flat over $Y$. 
\end{lem}

\Pf . (i) The fact that $F$ is flat over $Y$ implies that the functor
$$D_{qc}(Y)\to D_{qc}(Z): G\mapsto F\otimes g^*G$$
is $t$-exact. Now the assertion follows from the isomorphism 
$$(f,g)_*F\circ_Y G\simeq f_*(F\otimes g^*G)$$
for $G\in\Qcoh(Y)$ (that follows from the projection formula for the map $(f,g)$).

\noindent
(ii) One can repeat the above argument using the projection formula for sheaves on
formal schemes (see Theorem \ref{formal-proj-formula}) and the fact that
$f_{\sigma *}$ is exact for {\it ldu}-affine $f$.
\ed

\begin{lem}\label{lim-circ-lem} Let $X$, $Y$ and $Z$ be schemes, where either $X$ or $Z$
is flat over $S$.

\noindent
(i) Let $K$ be a $\circ_Y$-flat
quasicoherent sheaf on $X\times Y$.
Then the functor 
$$\Qcoh(Y\times Z)\to \Qcoh(X\times Z): F\mapsto K\circ_Y F$$
commutes with inductive limits.

\noindent (ii) Let $(K_i)$ be an inductive system of $\circ_Y$-flat quasicoherent sheaves on $X\times Y$. Then $\lim_i K_i$ is still $\circ_Y$-flat, and for every quasicoherent
sheaf $F$ on $Y\times Z$ the natural map 
$$\lim (K_i\circ_Y F)\to (\lim K_i)\circ_Y F$$
is an isomorphism.
\end{lem}

\Pf . (i) By definition, $K\circ_Y F=p_{XZ*}(p_{XY}^*K\otimes p_{YZ}^*F)$,
where we use projections from $X\times Y\times Z$ to double products of factors. Recall that
$K\circ_Y F$ is a sheaf by Lemma \ref{flat-lem}(iii), and $K$ is 
flat over $Y$ by Lemma \ref{flat-lem}(ii). It follows that the sheaf $p_{XY}^*K$ is flat over $Y\times Z$,  so we can use the underived functors in the above
formula for $K\circ_Y F$. All of them commute with inductive limits (for the push-forward
see Lemma \ref{lim-push-lem}).

\noindent
(ii) Note that $\lim K_i$ is flat over $Y$, so for $F\in \Qcoh(Y\times Z)$ we
have 
$$p_{XY}^*(\lim K_i)\otimes F\simeq (\lim p_{XY}^*K_i)\otimes F\simeq \lim (p_{XY}^*K_i\otimes F),$$
where we can use the underived functors. Now using Lemma \ref{lim-push-lem} we find that
the derived push-forward of this sheaf by $p_{XZ}$ is concentrated in degree $0$ and is isomorphic
to $\lim K_i\circ_Y F$.
\ed

\subsection{Kernel algebras}
\label{ker-alg-sec}

Let $X$ be a flat scheme over $S$. We denote by $\De:X\to X\times X$ the (relative over $S$) diagonal
embedding. Note that this is an affine morphism, since $X$ and $S$ are semi-separated.

\begin{defi} (i) A {\it kernel algebra} 
over $X$ is an object $\AA\in D_{qc}(X\times X)$ equipped with a morphism 
$\mu:\AA\circ_X \AA\to \AA$ (product) and a morphism $u:\De_*\OO_X\to \AA$ (unit),
subject to the usual associativity and unit axioms. Note that in these axioms we use
the isomorphisms 
$$(\AA\circ_X\AA)\circ_X\AA\simeq \AA\circ_X(\AA\circ_X\AA) \text{ and }
\AA\circ_X\De_*\OO_X\simeq\De_*\OO_X\circ_X\AA\simeq\AA$$ 
(see Lemma \ref{circle-prop-lem}).
Homomorphisms between kernel
algebras over the same scheme $X$ are defined in an obvious way.
A kernel algebra $\AA$ over $X$ is called {\it pure} if $\AA$ is 
a quasicoherent sheaf on $X\times X$, $\circ$-flat over $X$ with respect to both projections.

\noindent
(ii) A left (resp., right) module over a pure kernel algebra $\AA$ is a quasicoherent sheaf $F$ on $X$
equipped with a morphism $\AA\circ_X F\to F$ (resp., $F\circ_X \AA\to F$) satisfying the usual
associativity and unitality axioms. We denote by $\AA-\mod$ the category of left $\AA$-modules.

\noindent
(iii) We say that a pure kernel algebra $\AA$ over $X$ is {\it finite} if 
$\AA$ is a coherent sheaf on $X\times X$.
In this case a (left or right) module $F$ over $\AA$ is called {\it coherent} if $F$ is a coherent
sheaf on $X$. We denote by $\AA-\mod^c$ the category of coherent $\AA$-modules.
\end{defi}

Henceforward, whenever we consider kernel algebras over $X$ we implicitly assume $X$
to be flat over $S$.

Note that if $\AA$ is a pure kernel algebra then
by Lemma \ref{flat-lem}(ii), the object $\AA\circ_X \AA\in D_{qc}(X\times X)$ (resp.,
$\AA\circ_X F$, $F\circ_X \AA$) appearing in the above definition is a sheaf. Using exactness of
the functor $F\mapsto \AA\circ_X F$ 
one can immediately check that $\AA-\mod$ (resp, $\AA-\mod^c$ if $\AA$
is finite) is an abelian category.

\begin{lem}\label{lim-mod-lem}
Let $\AA$ be a pure kernel algebra over $X$. Then for any inductive system of
left $\AA$-modules $(F_i)$ there is a natural $\AA$-module structure on the quasicoherent
sheaf $F=\lim_i F_i$. The obtained $\AA$-module $F$ is the limit of $(F_i)$ in the category
of $\AA$-modules. 
\end{lem}

\Pf . For each map $F_i\to F_j$ in this inductive system there is an induced map of
quasicoherent sheaves $\AA\circ_X F_i\to \AA\circ_X F_j$, so we get a map of inductive systems
$(\AA\circ_X F_i)\to (F_i)$. Passing to the limit and using Lemma \ref{lim-circ-lem}(i) we get a map 
$$\AA\circ_X\lim F_i\simeq\lim \AA\circ_X F_i\to \lim F_i,$$
so we get an action of $\AA$ on $F=\lim F_i$. It is easy to check that this is a structure of an $\AA$-module on $F$, and that it is the limit of $(F_i)$ in the category of $\AA$-modules.
\ed

\begin{lem}\label{ex-com-lem} 
Let $\AA$ be a pure kernel algebra over $X$, $F$ a left $\AA$-module.
Then the complex of left $\AA$-modules
$$\ldots\stackrel{d_2}{\ra} \AA\circ_X \AA\circ_X F\stackrel{d_1}{\ra} 
\AA\circ_X F\stackrel{d_0}{\ra} F\to 0$$
is exact, where the differentials $d_i$ are alternating sums of the appropriate operations in
two consecutive factors.
\end{lem}

\Pf . Similarly to the case of associative algebras, we can 
consider the map $s_n:\AA^{\circ_X n}\circ_X F\to \AA^{\circ_X (n+1)}\circ_X F$
induced by the unit morphism $u:\De_*\OO_X\to\AA$ in the first factor.
One immediately checks that $(s_n)$ is a contracting homotopy for our complex, hence
it is acyclic.
\ed

\begin{lem}\label{coh-ker-lem}
Let $\AA$ be a finite pure kernel algebra over $X$. Then the natural functor
$D^-(\AA-\mod^c)\to D^-(\AA-\mod)$ is fully faithful, and its essential image
consists of all complexes of $\AA$-modules with bounded above coherent cohomology.
\end{lem}

\Pf . This is proved in the same way as for usual sheaves of $\OO_X$-modules 
(see e.g., \cite{Bezr}, Cor.~1).
One only has to check that for every surjection of $\AA$-modules $F\to G$,
where $G$ is coherent, there exists a coherent $\AA$-submodule $F_0\sub F$ that still
surjects onto $G$. We can start with a coherent subsheaf $F'\sub F$ surjecting onto $G$
and then replace it with the image of the corresponding morphism of $\AA$-modules
$\AA\circ_X F'\to F$.
\ed

\begin{exs}
1. If $\AA$ is a D-algebra on $X$ (see \cite{BB}), 
flat as a left and as a right $\OO_X$-module, then the
associated sheaf $\delta\AA$ on $X\times X$ has a structure of pure kernel algebra, and
the modules over $\delta\AA$ are exactly modules over $\AA$ in the usual sense.
Note that $\circ$-flatness of $\delta\AA$ over $X$ (with respect to either projection)
reduces to the usual flatness, since $\delta\AA$ is supported on the diagonal.

\noindent
2. Let $G$ be a discrete group acting on a scheme $X$. Then we have a natural pure kernel algebra
structure on
$$\AA^G_X:=\oplus_{g\in G}\OO_{\Ga_g},$$
where $\Ga_g\sub X\times X$ is the graph of the action of $g^{-1}\in G$ on $X$, i.e.,
$\Ga_g=\{(gx,x)\ |\ x\in X\}$. The product
is induced by the natural isomorphisms $\OO_{\Ga_g}\circ_X\OO_{\Ga_{g'}}\simeq\OO_{\Ga_{gg'}}$.
For every $F\in \Qcoh(X)$ we have natural isomorphisms $\OO_{\Ga_g}\circ_X F\simeq (g^{-1})^*F$.
It is easy to see that $\AA^G_X$-modules are exactly $G$-equivariant quasicoherent sheaves
on $X$. Later we will generalize this example to actions of finite group schemes and of formal
groups (see Corollary \ref{action-ker-cor}).

\noindent
3. Assume that $X$ and $S$ are affine, $X=\Spec A$, $S=\Spec R$. Then a pure kernel
algebra $\AA$ over $X$ corresponds to an associative ring $\Ga(\AA)$ equipped with a homomorphism
$A\to \Ga(\AA)$ such that the image of $R$ is in the centre of $\Ga(\AA)$, and such that
$\Ga(\AA)$ is flat as a left (resp., right) $A$-module. 
Modules over $\AA$ are usual modules over $\Ga(\AA)$.
\end{exs}

We are going to describe several basic operations with kernel algebras.
Let us denote by $\si:Y\times X\to X\times Y$ the permutation isomorphism.
Then for $K\in D_{qc}(X\times Y)$ and $L\in D_{qc}(Y\times Z)$ we have a natural isomorphism
\begin{equation}\label{perm-ker-eq}
\si^*(K\circ_Y L)\simeq \si^*L\circ_Y \si^*K
\end{equation}
on $Z\times X$.

\begin{defi}
For a kernel algebra 
$\AA$ over $X$ the {\it opposite} kernel algebra
$\AA^{opp}$ is $\si^*\AA$, where $\si$ is the permutation of factors in $X\times X$.
The product structure is induced by the product structure on $\AA$ using the isomorphism
\eqref{perm-ker-eq}. Note that the opposite of a pure kernel algebra is pure.
The isomorphism \eqref{perm-ker-eq} also shows that in the pure case 
there is an equivalence of the category of left $\AA^{opp}$-modules
with the category of right $\AA$-modules (and vice versa).
\end{defi}

For $K\in D_{qc}(X\times Y)$ and $K'\in D_{qc}(X'\times Y')$ we denote the 
{\it external product} 
$$K\Box K'=\si_{23}^*(p_{XY}^*K\otimes p_{X'Y'}^*K')\in D_{qc}(X\times X'\times Y\times Y'),$$
where 
$\si_{23}:X\times X'\times Y\times Y'\to X\times Y\times X'\times Y'$ is
the transposition. 

\begin{lem}\label{ext-prop-lem} 
(i) For $K\in D_{qc}(X\times Y)$, $L\in D_{qc}(Y\times Z)$,
$K'\in D_{qc}(X'\times Y')$ and $L'\in D_{qc}(Y'\times Z')$ one has
$$(K\Box K')\circ_{Y\times Y'}(L\Box L')\simeq (K\circ_Y L)\Box (K'\circ_{Y'} L')$$
provided $Y$ and $Y'$ are flat over $S$.

\noindent (ii) For $K$ and $K'$ as above and for $F\in D_{qc}(Y\times Y')$ one has a natural
isomorphism on $X\times X'$:
$$(K\Box K')\circ_{Y\times Y'}F\simeq (K\circ_Y F)\circ_{Y'}(\si^*K')$$
provided either $X'$ or $Y$ is flat over $S$.

\noindent (iii) Assume that $X$ and $X'$ are flat over $S$.
If $K$ (resp., $K'$) is a quasicoherent sheaf on $X\times Y$ (resp., $X'\times Y'$), 
$\circ$-flat over $Y$ (resp., $Y'$) then $K\Box K'$ is $\circ$-flat over $Y\times Y'$.
\end{lem}

\Pf . Parts (i) and (ii) are easy applications of the projection and base change formulae.
Part (iii) follows from (ii) and from Lemma \ref{flat-lem}(iii).
\ed

\begin{defi}
If $\AA$ is a kernel algebra over $X$ and $\BB$ is a kernel algebra over $Y$ then 
their {\it external product} $\AA\Box \BB$ has a natural structure of the kernel algebra over $X\times Y$ 
induced by an isomorphism of Lemma \ref{ext-prop-lem}(i). 
If $\AA$ and $\BB$ are pure then so is $\AA\Box \BB$ (the required $\circ$-flatness 
follows from Lemma \ref{ext-prop-lem}(iii)).
\end{defi}


\begin{lem}\label{prod-mod-lem} 
Let $\AA$ (resp., $\BB$) be a pure kernel algebra over $X$ (resp., $Y$). 
Let $F$ be a quasicoherent sheaf on $X\times Y$. 
Then
an $\AA\Box \BB^{opp}$-module structure on $F$ is determined by a pair of
morphisms $\a:\AA\circ_X F\to F$, $\b:F\circ_Y \BB\to F$ that satisfy the left and right module axioms
and that commute with each other, i.e., the following diagram is commutative
\begin{diagram}
\AA\circ_X F\circ_Y \BB &\rTo{}& \AA\circ_X F\\
\dTo{}&&\dTo{}\\
F\circ_Y\BB &\rTo{}& F
\end{diagram}
\end{lem}

\Pf . This follows easily from the natural isomorphism
$$(\AA\Box\BB^{opp})\circ_{X\times Y} F\simeq A\circ_X F\circ_Y\BB$$
provided by Lemma \ref{ext-prop-lem}(ii).
\ed

For example, taking $\BB=\De_*\OO_Y$ we see that an $\AA\Box\De_*\OO_Y$-module
structure on a quasicoherent sheaf $F$ on $X\times Y$ is given by a morphism
$\AA\circ_X F\to F$ on $X\times Y$ satisfying the module axioms.

\begin{defi} Let $f:X\to Y$ be a morphism of flat $S$-schemes.
If $\AA$ is a kernel algebra over $X$ then
$(f\times f)_*\AA$ has a natural structure of a kernel algebra over $Y$
induced by the morphism 
$$(f\times f)_*\AA\circ_Y (f\times f)_*\AA\to (f\times f)_*(\AA\circ_X \AA)$$
(see \eqref{push-circle-mor}). If in addition $f$ is
flat and affine, $Y$ is flat over $S$, and $\AA$ is pure, 
then $(f\times f)_*\AA$ is also pure. Indeed, the fact
that $(f\times f)_*\AA$ is $\circ$-flat over $Y$ follows from Lemma \ref{flat-lem}(i).
We call the kernel algebra $(f\times f)_*\AA$
{\it the push-forward of $\AA$ under $f$}.
\end{defi}

\begin{lem}\label{push-alg-lem} 
Let $f:X\to Y$ be a flat affine morphism, where $Y$ is flat over $S$.
Then for a pure kernel algebra $\AA$ over $X$ the functor $F\mapsto f_*F$ induces an equivalence
of the category of $\AA$-modules with the category of $(f\times f)_*\AA$-modules.
If in addition $f$ is finite and $\AA$ is finite then so is
$(f\times f)_*\AA$ and we have an equivalence
$\AA-\mod^c\simeq (f\times f)_*\AA-\mod^c$.
\end{lem}

\Pf . Note that the assertion is clear in the case when $Y$ and $S$ are affine.
In the general case, for an $\AA$-module $F$ the natural 
$(f\times f)_*\AA$-module structure on  $f_*F$ is given by the following composition
$$(f\times f)_*\AA\circ_Y f_*F\wt{\to} f_*(\AA\circ_X f^*f_*F)\to f_*(\AA\circ_X F)\to f_*F,$$
where the first isomorphism is derived from \eqref{push-circle-eq}, and the second arrow
is induced by the adjunction map $f^*f_*F\to F$. Conversely, assume we are given
an $(f\times f)_*\AA$-module structure on a quasicoherent sheaf $G$ over $Y$.
Then it induces an $f_*\OO_X$-module structure on $G$, so that $G\simeq f_*F$.
Furthermore, it is easy to see that the structure morphism
from $(f\times f)_*\AA\circ_Y f_*F\simeq f_*(\AA\circ_X f^*f_*F)$ to $f_*F$ commutes
with the $f_*\OO_X$-module structures on both sides, so it is induced by a morphism
$\AA\circ_X f^*f_*F\to F$. We claim that this morphism factors through the canonical
morphism $\a:\AA\circ_X f^*f_*F\to\AA\circ_X F$. Indeed, since $\a$ is surjective it is enough
to check this locally. But in the affine case we know this to be true.
The associativity of the obtained action $\AA\circ_X F\to F$ can also be checked locally.

In the case when $f$ is finite the functor $f_*$ preserves coherence, which implies 
the second assertion.
\ed

\begin{rem} Let $\AA$ be a pure kernel algebra on $X$, and let $g:Y\to Y'$ be an affine morphism
of flat $S$-schemes.
Then for an $\AA\Box\De_*\OO_Y$-module $F$ on $X\times Y$ we can introduce a natural
$\AA\Box\De_*\OO_{Y'}$-module structure on $(\id_X\times g)_*F$. Namely,
it is given by the map
$$\AA\circ_X(\id_X\times g)_*F\simeq (\id_X\times g)_*(\AA\circ_X F)\to (\id_X\times g)_*F,$$
where the first isomorphism follows from Lemma \ref{circle-push-pull-lem}(ii).
If in addition $g$ is flat then by Lemma \ref{push-alg-lem}, we can view $(\id_X\times g)_*F$
as a module over the pure kernel algebra $\AA\Box\De_*g_*\OO_Y$, and the above
$\AA\Box\De_*\OO_{Y'}$-module structure is induced by the homomoprhism of kernel algebras
$\AA\Box\De_*\OO_{Y'}\to\AA\Box\De_*g_*\OO_Y$.
\end{rem}

\subsection{Cech resolutions}\label{Cech-res-sec}

We are going to introduce a notion of compatibility of a kernel algebra $\AA$ with an open covering, so
that the Cech resolution of an $\AA$-module with respect to this covering would inherit
the $\AA$-module structure. Here is a general framework for this construction.

\begin{defi}
Let $(\CC,*)$ be a monoidal category with a unit object $I$, and let
$(A,\mu:A*A\to A, u:I\to A)$ be an algebra
in $\CC$, $O$ an object of $\CC$. We say that a morphism
$$\si:A*O\to O*A$$
{\it is compatible with the algebra structure on} $A$ if the following diagram is commutative
\begin{diagram}
A*A*O&\rTo{}& A*O*A &\rTo{}& O*A*A\\
\dTo{} &&&&\dTo{}\\
A*O&&\rTo{}&&O*A
\end{diagram}
where the horizontal (resp., vertical) arrows are induced by $\si$ (resp., $\mu$),
and in addition $\si(u*\id_O)=\id_O*u$ (equality of morphisms from $O$ to $O*A$). 
\end{defi}

The next lemma gives a pattern for constructing complexes similar to Cech resolutions. 

\begin{lem}\label{formal-lem}
Let $(\CC,*)$, $A$, $(O,\si)$ be as in the above definition.
Assume in addition that the category $\CC$ is additive (where the monoidal structure is given by an additive functor), 
and $O$ is equipped with a morphism $\iota:I\to O$ such that
$\si(\id_A*\iota)=\iota*\id_A$ (equality of morphisms from $A$ to $O*A$).
Then there is a natural structure of a dg-algebra (with a unit) over $(\CC,*)$ on 
$$T^{\bullet}(O)*A:=A\oplus O*A\oplus O*O*A\oplus\ldots,$$
where multiplication is induced by the product on $A$ and $\si$, the differential $d$
is characterized by 
$$d(A)=0, \ \ d|_O=(\id_O*\iota-\iota*\id_O)*\id_A:O*A\to O*O*A.$$
The morphism $\iota*\id_A:A\to O*A$ defines a map of complexes of $A$-bimodules
$$A\to T^{>0}(O)*A[1].$$
Similarly, for every left $A$-module $M$ there is a natural structure of a left (unital) dg-module over
$T^{\bullet}(O)*A$ on $T^{\bullet}(O)*M$ and a map of complexes of $A$-modules
$M\to T^{>0}(O)*M[1]$.
\end{lem}

The proof is left to the reader.

We will apply this Lemma to the monoidal structure given by the convolution. The algebra 
$A$ will be a kernel algebra and the object 
$O$ will be the sheaf $\De(\UU)$
concentrated on the diagonal associated with an open
covering $\UU$ (see \eqref{diag-U-eq} below). 
The obtained complex of $A$-modules
$T^{>0}(O)*M[1]$ will be the Cech resolution.

By an open covering $\UU=(U_i)$ of $X$ we always mean an open covering with respect
to the flat ({\it fppf}) topology. 
For a quasicoherent sheaf $F$ on $X$ and an open covering $(U_i)_{i\in I}$ of $X$, 
such that the maps
$j_i:U_i\to X$ are affine, let us denote by $C^{\bullet}_{\UU}(F)$ the corresponding 
Cech resolution of $F$ (where indices are allowed to coincide). Namely,
$$C^p_{\UU}(F)=\oplus_{J\sub I^{p+1}}j_{J,*}j_J^*F,$$
where for $J=(i_0,\ldots,i_p)$ we denote by $j_J:U_{i_0}\times_X\ldots\times_X U_{i_p}\hra X$
the natural embedding. 
Note that in this case all the fibered products of $j_i$'s are still affine maps, so the corresponding
push-forward functors are exact.
For example,
we can take Cech resolutions corresponding to
affine open coverings (the corresponding morphisms $j_i:U_i\to X$ will be automatically affine
since $X$ is semi-separated).

These Cech resolution can be viewed in the context of Lemma \ref{formal-lem} as follows.
Consider the following sheaf on $X\times X$: 
\begin{equation}\label{diag-U-eq}
\De(\UU):=\De_*C^0_{\UU}(\OO_X)=\oplus_i\De_*j_{i*}\OO_{U_i}.
\end{equation}
It is equipped with the natural map $\iota:\De_*\OO_X\to\De(\UU)$, hence
we have the corresponding dg-algebra $T^{\bullet}(\De(\UU))$ with respect to the circle product.
Set $\OO^{\bullet}(X,\UU):=T^{>0}(\De(\UU))[1]$.
It is easy to see that 
$$\OO^{\bullet}(X,\UU)=\De_*\CC^{\bullet}_{\UU}(\OO_X).$$
Furthermore, for $F\in\Qcoh(X\times Y)$ (resp., $G\in\Qcoh(Y\times X)$)
we have natural isomorphisms of resolutions of $F$ (resp., $G$)
\begin{equation}\label{Cech-eq1}
\OO^{\bullet}(X,\UU)\circ_X F\simeq C^{\bullet}_{(U_i\times Y)}(F),
\end{equation}
\begin{equation}\label{Cech-eq2}
G\circ_X \OO^{\bullet}(X,\UU)\simeq C^{\bullet}_{(Y\times U_i)}(G).
\end{equation}
Using Lemma \ref{formal-lem} we arrive to the following definition which specifies
when Cech resolutions for modules over a kernel algebra still carry an action of this kernel algebra.

\begin{defi}
We say that an open covering $(U_i)$ of $X$ (with affine maps $j_i:U_i\to X$) is {\it compatible} with 
a kernel algebra $\AA$ over $X$ if there is an
isomorphism in $D(X\times X)$
\begin{equation}\label{adapt-eq}
\si:\AA\circ_X\De(\UU)\wt{\ra}\De(\UU)\circ_X \AA
\end{equation}
compatible with the kernel algebra structure on $\AA$ and with the canonical
homomorphism $\De_*\OO_X\to\De(\UU)$. 
We say that a pure kernel algebra $\AA$ over $X$ is of {\it affine type} if there exists an
open {\it affine} covering $(U_i)$ of $X$ with respect to the flat topology, compatible with $\AA$.
\end{defi}


\begin{exs} 1. Let $\pi:X\to S$ be the projection. Then for any open covering $(S_i)$ of $S$ such that
the maps $S_i\to S$ are affine,
the induced open covering 
$(U_i=S_i\times_S X)$ of $X$ is compatible with any kernel algebra over $X$.
Indeed, this follows from the natural isomorphisms of $X\times_S X$-schemes for each $i$:
$$U_i\times_S X\simeq S_i\times_S (X\times_S X)\simeq (X\times_S X)\times_S S_i\simeq X\times_S U_i.$$ 

\noindent 2. Let $\AA$ be a pure kernel algebra over $X$. If $(U_i)$ is an open affine covering of $X$ 
with respect to the Zariski topology and the support of 
$\AA|_{U_i\times X}$ (resp., $\AA|_{X\times U_i}$) is contained in 
$U_i\times U_i$ for every $i$, then $(U_i)$ is compatible with $\AA$. 
For example, if $\AA$ is supported on the diagonal in $X\times X$ (i.e., comes from a D-algebra)
then $\AA$ is of affine type.

\noindent 
3. Assume that $\AA$ is associated with an action of a discrete group $G$ on $X$ (see Example 2
in section \ref{ker-alg-sec}).
Take an open affine covering $(U_i)_{i\in I}$ of $X$ such that $G$ permutes the open subsets $U_i$.
In other words, we assume that there is an action of $G$ on the set of indices $I$ and isomorphisms
$a_g:U_i\wt{\to} U_{g(i)}$ defining the action of $G$ on $\sqcup_i U_i$ such that for each $g\in G$
the diagram 
\begin{diagram}
U_i &\rTo{a_g}&U_{g(i)}\\
\dTo{j_i} &&\dTo{j_{g(i)}}\\
X &\rTo{g}& X
\end{diagram}
is commutative.
Then it is easy to check that the covering $(U_i)$ is compatible with $\AA$.
For a generalization, see Proposition \ref{aff-ker-prop}.
\end{exs}

Assume that an open affine covering $\UU=(U_i)$ is compatible with a pure kernel algebra $\AA$ over $X$. Then by Lemma \ref{formal-lem},
we have a natural resolution of $\AA$ by $\AA\Box\AA^{opp}$-modules on $X\times X$
\begin{equation}\label{ker-alg-res-eq}
\AA^{\bullet}(\UU):=\OO^{\bullet}(X,\UU)\circ_X\AA\simeq\AA\circ_X \OO^{\bullet}(X,\UU)
\end{equation}
that can be identified with the
Cech resolution of $\AA$ with respect to either covering: $(U_i\times X)$ or
$(X\times U_i)$. 

We will also need truncated Cech resolutions. First, for every $N>0$ consider the complex
$\tau_{\le N}C^{\bullet}_{(U_i)}(\OO_X)$. Clearly it is still a resolution of $\OO_X$.
We are going to show that for sufficiently large $N$ this truncation will be a sufficiently good
replacement for the Cech complex (see Lemma \ref{flat-res-lem} below).

\begin{lem}\label{flat-Cech-lem} 
Let $\UU=(U_i)$ be an open covering of $X$ such that the morphisms
$j_i:U_i\to X$ are affine.

\noindent
(i) For every $N>0$ the kernel of the differential
$K^N=\ker(d^N:C^N_{\UU}(\OO_X)\to C^{N+1}_{\UU}(\OO_X))$ is locally 
a direct summand of $C^{N-1}_{\UU}(\OO_X)$.

\noindent
(ii) If $F$ is a quasicoherent sheaf on $X$ then we have
$$\tau_{\le N}C^{\bullet}_{\UU}(F)\simeq (\tau_{\le N}C^{\bullet}_{\UU}(\OO_X))\otimes F\simeq
\tau_{\le N}\OO^{\bullet}(X,\UU)\circ_X F.$$
\end{lem}

\Pf . (i) Localizing we can assume that one of the open subsets is $X$ itself.
In this case we have a contracting homotopy $h$ for the complex
\begin{equation}\label{Cech-com-O-eq}
\OO_X\to C^0_{\UU}(\OO_X)\to C^1_{\UU}(\OO_X)\to\ldots.
\end{equation}
Restricting $h$ to $K^N$ we get a splitting $K^N\to C^{N-1}_{\UU}(\OO_X)$ of the
surjection $C^{N-1}_{\UU}(\OO_X)\to K^N$ induced by $d^{N-1}$. Hence, $K^N$ is a direct summand
of $C^{N-1}_{\UU}(\OO_X)$. 

\noindent
(ii) As we have seen in (i), the short exact triples $0\to K^{N-1}\to C^{N-1}\to C^N\to 0$
are locally split. Therefore, tensoring them with any quasicoherent sheaf we still get
exact triples.
\ed

From our point of view the key feature of the Cech resolution is that it transforms flat sheaves
into complexes of $\circ$-flat sheaves. Let us show that the same is true for sufficiently big
truncated Cech resolutions.

\begin{lem}\label{flat-res-lem} 
Let $V$ be a flat quasicoherent sheaf on $X$, and let
$\UU=(U_i)$ be an open affine covering of $X$ (with respect to the flat topology).
Let $n_0$ be the integer such that $H^i(X,F)=0$ for $i>n_0$ for all quasicoherent sheaves $F$.
Then there exists $n_0$ such that for $N>n_0$ all terms of the truncated Cech resolution
$\tau_{\le N}C^{\bullet}_{\UU}(V)$ are $\circ$-flat over $X$.
\end{lem}

\Pf . First, we note that by Lemma \ref{flat-lem}(i),(ii),
all the terms $C^n_{\UU}(V)$
are $\circ$-flat over $X$ (because $V$ is flat over $X$). Therefore, we only have to check that
$K^N(V)=\ker(C^N_{\UU}(V)\to C^{N+1}_{\UU}(V))$ is $\circ$-flat over $X$. 
Recall that by Lemma \ref{flat-Cech-lem}, 
$K^N(V)$ is the tensor product of $V$ with a flat sheaf $K^N(\OO_X)$, so it is flat over $X$.
Therefore, for every quasicoherent sheaf $G$ on $X$ 
the sequence of sheaves on $X$
$$0\to V\otimes G\to C^0_{\UU}(V)\otimes G\to\ldots\to C^{N-1}_{\UU}(V)\otimes G\to 
K^N(V)\otimes G\to 0$$
is exact. 
Since in this sequence all the sheaves $C^i_{\UU}(V)\otimes G$ are $\pi$-acyclic,
where $\pi:X\to S$ is the projection,
our assumption $N>n_0$ implies that $K^N(V)\otimes G$ is also $\pi$-acyclic. Hence,
$$G\circ_X K^N(V)=\pi_*(K^N(V)\otimes G)$$ 
is concentrated in degree $0$, i.e., $K^N(V)$ is $\circ$-flat over $X$. 
\ed

\subsection{Convolution for modules over kernel algebras}

Let $(X,\AA)$, $(Y,\BB)$, $(Z,\CC)$ be three (flat) $S$-schemes equipped with pure kernel algebras.
We would like to find an analog of the convolution operation considered in section \ref{circle-sec}
for modules. Namely, 
for $F\in D(\AA\Box\BB^{opp}-\mod)$ and $G\in D(\BB\Box\CC^{opp}-\mod)$ we would like to
define $F\circ_{\BB} Y\in D(\AA\Box\CC^{opp}-\mod)$
that would globalize the operation of tensor product of bimodules. Below we will show how to do this under a technical assumption on $\BB$ (see Proposition \ref{ass-ker-alg-prop}).
We start by defining an underived
version and then will use appropriate resolutions.  

Let $F$ be a (left) $\AA\Box\BB^{opp}$-module, and let $G$ be a $\BB\Box\CC^{opp}$-module.
Assume that either $F$ or $G$ is $\circ$-flat over $Y$. Then we define
$$F\ov{\circ}_\BB G=\coker(F\circ_Y \BB\circ_Y G\stackrel{\a}{\ra} F\circ_Y G),$$
where $\a$ is the difference of two natural maps, one induced by the left action
$\BB\circ_Y G\to G$ and the other by the right action $F\circ_Y \BB\to F$.
Note that $F\circ_Y G$ and $F\circ_Y \BB\circ_Y G$ are sheaves because of our $\circ$-flatness
assumption.
The left action of $\AA$ on $F$ and the right action of $\CC$ on $G$ induce 
the structure of the $\AA\Box\CC^{opp}$-module on $F\ov{\circ}_\BB G$ (see Lemma 
\ref{prod-mod-lem}).

Note that if $\AA'\to\AA$ (resp., $\CC'\to\CC$) is a homomorphism of pure kernel algebras
on $X$ (resp., $Z$) then the operation $F\circ_Y G$ is compatible with the natural
restriction of scalars from $\AA$ to $\AA'$ (resp., from $\CC$ to $\CC'$).
Also, in the case $\BB=\De_*\OO_Y$ we get the usual circle operation over $Y$ (since in this
case $\a=0$).

\begin{lem}\label{ass-circ-alg-lem} 
Let $(X,\AA)$, $(Y,\BB)$, $(Z,\CC)$, $(T,\DD)$ be 
schemes equipped with pure kernel algebras, and let 
$F\in \AA\Box\BB^{opp}-\mod$, $G\in \BB\Box\CC^{opp}-\mod$ and 
$H\in \CC\Box\DD^{opp}-\mod$.

\noindent
(i) One has an isomorphism of $\AA\Box\BB^{opp}$-modules $F\ov{\circ}_{\BB}\BB\simeq F$.

\noindent
(ii) Assume that $F$ is $\circ$-flat over $Y$ and $H$ is $\circ$-flat over
$Z$. Then one has a functorial isomorphism of $\AA\Box\DD^{opp}$-modules
$$(F\ov{\circ}_\BB G)\ov{\circ}_\CC H\simeq F\ov{\circ}_\BB (G\ov{\circ}_\CC H).$$
\end{lem}

\Pf . (i) This follows easily from Lemma \ref{ex-com-lem}.

\noindent (ii)
We have
\begin{align*}
&(F\ov{\circ}_\BB G)\ov{\circ}_\CC H\simeq
\coker((F\ov{\circ}_\BB G)\circ_Z \CC\circ_Z H\to (F\ov{\circ}_\BB G)\circ_Z H)\simeq\\
&\coker(F\circ_Y G\circ_Z \CC\circ_Z H\oplus F\circ_Y \BB\circ_Y G\circ_Z H\stackrel{\de}{\ra}
F\circ_Y G\circ_Y H),
\end{align*}
where $\de$ is induced by the appropriate operations between the consecutive factors in the $\circ$-products. Similarly, we get an isomorphism $F\ov{\circ}_\BB (G\ov{\circ}_\CC H)\simeq\coker(\de)$.
\ed

\begin{defi} Let $F$ be a left $\AA\Box\BB^{opp}$-module on $X\times Y$, where
$\AA$ (resp., $\BB$) is a pure kernel algebra over $X$ (resp., $Y$). 

\noindent (i)
We say that $F$ is {\it $\circ$-flat over $\BB$} (or simply {\it $\circ_{\BB}$-flat}) if 
it is $\circ$-flat over $Y$ and the functor 
$$\BB-\mod\to \AA-\mod: G\mapsto F\ov{\circ}_{\BB} G$$
is exact.

\noindent (ii)
$F$ is {\it $\circ$-free over }$\BB$ if $F\simeq F_0\circ_Y\BB$, where $F_0$ is a 
$\AA\Box (\De_*\OO_Y)$-module, $\circ$-flat over $Y$.

\noindent (iii)
A complex $F^{\bullet}$ of left $\AA\Box\BB^{opp}$-modules on $X\times Y$
is called $q-\circ$-flat over $\BB$ (or $q-\circ_{\BB}$-flat) if it has $\circ_Y$-flat terms and
for any exact complex $G^{\bullet}$ of $\BB$-modules on $Y$ the complex 
$F^{\bullet}\ov{\circ}_{\BB} G^{\bullet}$ is exact (where the definition of $\ov{\circ}_{\BB}$ is
extended to complexes in the standard way).
\end{defi}

For example, if $\BB=\De_*\OO_Y$ then $F$ is automatically
$\circ$-free over $\BB$ (and $\circ_{\BB}$-flat) provided it is $\circ_Y$-flat.

\begin{lem}\label{flat-alg-lem} 
(i) Assume that $F\in \AA\Box\BB^{opp}-\mod$ is $\circ$-flat over $\BB$.
Then for every scheme $Z$ equipped with a pure kernel algebra $\CC$ the functor
$$\BB\Box\CC^{opp}-\mod\to \AA\Box\CC^{opp}-\mod: G\mapsto F\ov{\circ}_{\BB} G$$
is exact. Similarly, if $F$ is a $q-\circ_{\BB}$-flat complex of $\AA\Box\BB^{opp}$-modules
then the natural extension of the above functor to complexes sends exact complexes to exact complexes.

\noindent (ii) 
Assume that $F\in \AA\Box\BB^{opp}-\mod$ is $\circ$-free over $\BB$. Then 
$F$ is $\circ$-flat over $\BB$.
\end{lem}

\Pf . (i) Composing with the natural forgetful functor we can assume that $\CC=\De_*\OO_Z$.
Also, we can assume that $Z$ is affine. Then it is enough to check the exactness of the functor
$G\mapsto p_{1*}(F\ov{\circ}_{\BB}G)$, where $p_1:X\times Z\to X$ is the projection.
Now the assertion follows from the isomorphism 
$p_{1*}(F\ov{\circ}_{\BB}G)\simeq F\ov{\circ}_{\BB}(p_{1*}G)$ that one can easily derive from
Lemma \ref{circle-push-pull-lem}(ii) (recall that $p_{1*}G$ has a natural $\BB\otimes\De_*\OO_Z$-module structure, see Remark after Lemma \ref{push-alg-lem}). The same argument works for the second statement.

\noindent (ii)
Assume that
$F=F_0\circ_Y \BB$, where $F_0$ is a $\AA\Box(\De_*\OO_Y)$-module on $X\times Y$,
$\circ$-flat over $Y$.
Note that in this case $F$ is also $\circ_Y$-flat by Lemma \ref{flat-lem}(iv).
It is enough to construct a functorial isomorphism of $\AA$-modules
$$(F_0\circ_Y \BB)\ov{\circ}_\BB G\simeq F_0\circ_Y G$$
for $G\in\BB-\mod$. Using the definition and the $\circ_Y$-flatness of $F_0$ we get
$$(F_0\circ_Y \BB)\ov{\circ}_\BB G\simeq \coker(F_0\circ_Y\BB\circ_Y\BB\circ_Y G\to
F_0\circ_Y\BB\circ_Y G)\simeq
F_0\circ_Y \coker(\BB\circ_Y \BB\circ_Y G \ra \BB\circ_Y G).$$
By Lemma \ref{ex-com-lem}, this is isomorphic to $F_0\circ_Y G$. 
\ed

The following lemma is an analogue of the well known properties of $q$-flatness
in the case of the usual tensor product of sheaves (cf. \cite{Lipman}, (2.5.4)).

\begin{lem}\label{lim-flat-lem} 
Let $X$ and $Y$ be schemes equipped with pure kernel algebras $\AA$ and $\BB$.

\noindent
(i) Let $(F^{\bullet}_i)$ be an inductive system of complexes of left $\AA\Box\BB^{opp}$-modules on $X\times Y$ (connected by chain maps). Assume that 
$F^{\bullet}_i$ is $q-\circ_{\BB}$-flat for every $i$. 
Then the same is true for the complex $\lim_i F^{\bullet}_i$.

\noindent
(ii) If $F^{\bullet}$ is a bounded above complex of $\AA\Box\BB^{opp}$-modules with
$\circ_{\BB}$-flat terms then $F^{\bullet}$ is $q-\circ_{\BB}$-flat.
\end{lem}

\Pf . (i) First, we note that the terms of $\lim F^{\bullet}_i$ are still $\circ_Y$-flat by Lemma 
\ref{lim-circ-lem}(ii). Furthermore, using the same Lemma one can easily see that for 
every complex $G^{\bullet}$ of $\BB$-modules the natural
map  
\begin{equation}\label{lim-ker-alg-eq}
\lim (F^{\bullet}_i\ov{\circ}_{\BB} G^{\bullet})\to (\lim F^{\bullet}_i)\ov{\circ}_{\BB} G^{\bullet}
\end{equation}
is an isomorphism. This immediately implies the result.

\noindent (ii) Such a complex is the inductive limit of bounded complexes with $\circ_{\BB}$-flat
terms for which the statement is clear (cf. [Lipman](2.5.4)).
\ed

Since we will need resolutions with terms that are $\circ$-flat over a given pure kernel algebra,
we give the following technical definition.

\begin{defi} (i) Let $\AA$ be a pure kernel algebra over $X$, 
and let $F$ be a $\AA\Box\De_*\OO_Y$-module on $X\times Y$. We say that
$F$ is {\it $\circ$-flattening over} $\AA$ if it is $\circ_Y$-flat, and
for every flat quasicoherent sheaf $V$ over $Y$ the $\AA$-module $F\circ_Y V$ is
$\circ$-flat over $\AA$.

\noindent
(ii) A kernel algebra $\AA$ over $X$ is called {\it left (resp., right) admissible} if it is pure and there exist
a quasi-isomorphism (given by a chain map)
$\AA\to\MM^{\bullet}$, where $\MM^{\bullet}$ is a bounded complex
of $\AA\Box\AA^{opp}$-modules on $X\times X$, such that each term $\MM^n$,
viewed as an $\AA\Box\De_*\OO_X$-module (resp., $\De_*\OO_X\Box\AA^{opp}$) is 
$\circ$-flattening over $\AA$ (resp., $\AA^{opp}$).
We say that a kernel algebra is {\it admissible} if it is left and right admissible.
\end{defi}


The main point of admissibility assumption is that it allows to construct
resolutions by modules that are $\circ$-flat over our kernel algebra (see Proposition
\ref{quasiisom-prop} below).
We will show also that every pure kernel algebra of affine type 
is admissible,
and that the push-forward of an admissible algebra under a flat affine morphism is admisible.

\begin{lem}\label{flatten-lem} 
Let $\AA$ be a pure kernel algebra over a flat $S$-scheme $X$, and let $Y$ and $Z$ be 
flat $S$-schemes.
If $F$ is a $\AA\Box\De_*\OO_Y$-module, $\circ$-flattening over $\AA$, then for 
every quasicoherent sheaf $V$ over $Y\times Z$, flat over $Y$, the $\AA\Box\De_*\OO_Z$-module
$F\circ_Y V$ is $\circ$-flat over $\AA$.
\end{lem}

\Pf . We have to check exactness of the functor
$$\AA-\mod\to\Qcoh(Z): G\mapsto G\ov{\circ}_{\AA}(F\circ_Y V).$$
Without loss of generality we can assume that $Z$ is affine.
Then it is enough to check the exactness of the functor
$$G\mapsto \pi_*(G\ov{\circ}_{\AA}(F\circ_Y V))\simeq G\ov{\circ}_{\AA}p_{X*}(F\circ_Y V)\simeq
G\ov{\circ}_{\AA}(F\circ_Y p_{Y*}(V)),$$
where $\pi:X\to S$, $p_X:X\times Z\to X$ and $p_Y:Y\times Z\to Y$ are the projections. But this follows
from the condition that $F$ is $\circ$-flattening over $\AA$, since $p_{Y*}(V)$ is a flat sheaf on $Y$.
\ed

\begin{prop}\label{quasiisom-prop} 
Let $(X,\AA)$ and $(Y,\BB)$ be schemes (flat over $S$) equpped with pure kernel algebras.
Assume that $\AA$ is left admissible. 
Then every complex in $D(\AA\Box\BB^{opp})$
is quasi-isomorphic to a $q-\circ_{\AA}$-flat complex of $\AA\Box\BB^{opp}$-modules.
\end{prop}

\Pf . 
First, let $F$ be an arbitrary $\AA\Box\BB^{opp}$-module on $X\times Y$.
Then we can find a functorial
surjection of quasicoherent sheaves $P\to F$, where $P$ is flat over $X\times Y$
(by Lemma \ref{flat-cover-lem}; note that $X\times_S Y$ is semi-separated). 
We have the induced surjection of $\AA\Box\BB^{opp}$-modules $\AA\circ_X P\circ_Y \BB\to F$.
Note that $P\circ_Y \BB$ is still flat over $X$ by Lemma \ref{flat-lem}(v).
This implies that for every bounded above complex of $\AA\Box\BB^{opp}$-modules $F^{\bullet}$ on $X\times Y$ 
there exists a quasi-isomorphism of the form
\begin{equation}\label{quasiisom-eq1}
c_1(F^{\bullet})=\AA\circ_X V^{\bullet}\to F^{\bullet},
\end{equation}
where $V^{\bullet}$ is a bounded above complex of $(\De_*\OO_X)\Box\BB^{opp}$-modules, 
flat over $X$.
Next, let $\AA\to\MM^{\bullet}$ be a quasi-isomorphism of $\AA$ with a bounded complex of 
$\AA\Box\AA^{opp}$-modules, $\circ$-flattening
over $\AA$ on the left. Consider the natural map of complexes of $\AA\Box\BB^{opp}$-modules
$$\a:c_1(F^{\bullet})\to c_2(F^{\bullet}):=\tot[\MM^{\bullet}\ov{\circ}_{\AA} c_1(F^{\bullet})]$$
(where $\tot$ denotes passing to the total complex of a bicomplex), induced by the map
$\AA\to\MM^{\bullet}$. By Lemma \ref{ass-circ-alg-lem}, we have isomorphisms 
$$\MM^{\bullet}\ov{\circ}_{\AA}(\AA\circ_X V^n)\simeq \MM^{\bullet}\circ_X V^n$$
of complexes of $\AA\Box\BB^{opp}$-modules.
Since for each $n$ the natural map 
$\AA\circ_X V^n\to \MM^{\bullet}\circ_X V^n$ is a quasi-isomorphism,
it follows that $\a$ is a quasi-isomorphism. On the other hand,
applying Lemma \ref{flatten-lem} we see that the terms of $c_2(F^{\bullet})$ are $\circ$-flat over $\AA$.
Hence, by Lemma \ref{lim-flat-lem}(ii), $c_2(F^{\bullet})$ is $\circ_{\AA}$-flat.

Finally, we use the standard trick to extend the above construction to the case of
an unbounded complex of $\AA\Box\BB^{opp}$-modules $F^{\bullet}$. First, we consider
for each $n\ge 0$ the truncated complex $\tau_{\le n}F^{\bullet}$ and consider the
quasi-isomorphisms of bounded above complexes
$$\tau_{\le n}F^{\bullet}\leftarrow c_1(\tau_{\le n}F^{\bullet})\rightarrow c_2(\tau_{\le n}F^{\bullet}).$$
By functoriality we have chain maps $c_i(\tau_{\le n}F^{\bullet})\to c_i(\tau_{\le n+1}F^{\bullet})$
for $i=1,2$, commuting with the above quasi-isomorphisms.
Now setting $c_i(F^{\bullet})=\lim_n c_i(\tau_{\le n}F^{\bullet})$ for $i=1,2$ gives us quasi-isomorphisms
$$F^{\bullet}\leftarrow c_1(F^{\bullet})\rightarrow c_2(F^{\bullet}).$$
Furthermore, by Lemma \ref{lim-flat-lem}(i), the complex $c_2(F^{\bullet})$ is $q-\circ_{\AA}$-flat. 
\ed

Returning to the situation in the beginning of this subsection, 
let $F^{\bullet}$ (resp. $G^{\bullet}$)
be a complex of $\AA\Box\BB^{opp}$-modules (resp., $\BB\Box\CC^{opp}$-modules).
Then assuming that $\BB$ is right admissible, we define $F^{\bullet}\circ_\BB G^{\bullet}$ 
as the total complex associated with the bicomplex
$\wt{F}^{\bullet}\ov{\circ}_\BB G^{\bullet}$, where $\wt{F}^{\bullet}$ is a $q-\circ_{\BB}$-flat
complex quasi-isomorphic to $F^{\bullet}$. Alternatively, if $\BB$ is left admissible then
we can use a $q-\circ_{\BB}$-flat
resolution of $G^{\bullet}$.

\begin{prop}\label{ass-ker-alg-prop}
Let $(X,\AA)$, $(Y,\BB)$, $(Z,\CC)$ be 
flat $S$-schemes equipped with pure kernel algebras, where 
$\BB$ is left (resp., right) admissible. Then there is a biexact
functor
$$D(\AA\Box\BB^{opp}-\mod)\times D(\BB\Box\CC^{opp}-\mod)\to
D(\AA\Box\CC^{opp}-\mod): (F^{\bullet},G^{\bullet})\mapsto F^{\bullet}\circ_\BB G^{\bullet}$$
defined as above using $q-\circ_{\BB}$-flat resolution of $G^{\bullet}$ (resp., $F^{\bullet}$).
Assume now that $\BB$ is right admissible and $\CC$ is left admissible, and let
$(T,\DD)$ be another flat $S$-scheme equipped with a pure kernel algebra. 
Then for $F\in D(\AA\Box\BB^{opp}-\mod)$, $G\in D(\BB\Box\CC^{opp}-\mod)$ and 
$H\in D(\CC\Box\DD^{opp}-\mod)$
we have a functorial isomorphism in $D(\AA\Box\DD^{opp}-\mod)$:
$$(F\circ_\BB G)\circ_\CC H\simeq F\circ_\BB (G\circ_\CC H).$$
\end{prop}

\Pf . To prove associativity we choose 
complexes quasiisomorphic to $F$ and $H$ 
that are $q-\circ$-flat over $\BB$ and $\CC$,
respectively.
\ed

The starting point for constructing examples of admissible kernel algebras is the following result.

\begin{lem}\label{adm-aff-lem}
Let $X$ be a flat $S$-scheme.
Then a pure kernel algebra of affine type over $X$ is admissible.
\end{lem}

\Pf . Let $\UU=(U_i)$ be an open affine covering compatible with a kernel algebra $\AA$.
We claim that for sufficiently large $N$ the truncated Cech
resolution $\tau_{\le N}\AA^{\bullet}(\UU)=\AA\circ_X\tau_{\le N}\OO^{\bullet}_X(\UU)$ 
(see \eqref{ker-alg-res-eq})
has terms that are $\circ$-flattening over $\AA$ on the left. 
Since $\tau_{\le N}\OO^{\bullet}_X(\UU)$ is $\circ$-flat over $X$ (recall that it is supported
on the diagonal), it follows that $\tau_{\le N}\AA^{\bullet}(\UU)$ is also $\circ$-flat over $X$
(on either side). Next, for a flat quasicoherent sheaf $V$ on $X$ we have
$$\tau_{\le N}\AA^{\bullet}(\UU)\circ_X V\simeq \AA\circ_X\tau_{\le N}\OO^{\bullet}_X(\UU)\circ_X V
\simeq\AA\circ_X\tau_{\le N}C^{\bullet}_{\UU}(V).$$
By Lemma \ref{flat-res-lem}, the terms of the complex $\tau_{\le N}C^{\bullet}_{\UU}(V)$ are $\circ$-flat
over $X$. Hence, the terms of $\tau_{\le N}\AA^{\bullet}(\UU)\circ_X V$ are $\circ$-free over $\AA$, so
they are $\circ$-flat over $\AA$ by Lemma \ref{flat-alg-lem}(ii).
A similar argument shows
that the terms of $\tau_{\le N}\AA^{\bullet}(\UU)$ are $\circ$-flattening over $\AA$ on the right.
\ed

Our next goal is to show that admissibility is preserved under push-forwards with respect to
affine flat morphisms.

\begin{lem}\label{subalg-circ-lem}
Let $\BB_0\to\BB$ be a homomorphism of pure kernel algebras over $Y$, and let $F$ be an $\AA\Box\BB^{opp}$-module. Assume that $\BB$ 
and $F$ are $\circ$-flat over $\BB_0$ acting on the right.
Then for every $G\in\BB\Box\CC^{opp}-\mod$ there is a natural isomorphism of 
$\AA\Box\CC^{opp}$-modules
$$F\ov{\circ}_\BB G\simeq\coker(F\ov{\circ}_{\BB_0} \BB\ov{\circ}_{\BB_0} G\to F\ov{\circ}_{\BB_0} G).$$
\end{lem}

\Pf . This follows easily from the definitions.
\ed


\begin{lem}\label{simple-push-forward-lem} 
Let $f:Y\to Y'$ be a flat affine morphism, and let $\BB_f:=\De_*(f_*\OO_Y)$ be the corresponding 
pure kernel algebra over $Y'$. Then $\BB_f$ is admissible, and
for $F\in D_{qc}(X\times Y)$ and $G\in D_{qc}(Y\times Z)$ one has
a natural isomorphism
\begin{equation}\label{push-forward-circle-isom}
(\id_X\times f)_*F\circ_{\BB_f} (f\times\id_Z)_*G\wt{\ra}
F\circ_Y G .
\end{equation}
If $F$ is a $\circ_Y$-flat quasicoherent sheaf on 
$X\times Y$ then $(\id_X\times f)_*F$ is $\circ$-flat over
$\BB_f$.
\end{lem}

\Pf . The pure kernel algebra $\BB_f$ is of affine type (being supported on the diagonal), so it is admissible. Hence, the right-hand side of \eqref{push-forward-circle-isom} is well-defined.
We have a natural morphism 
\begin{equation}\label{fl-af-push-forward-mor}
(\id_X\times f)_*F\circ_{Y'}(f\times\id_Z)_*G\to F\circ_Y G
\end{equation}
(see \eqref{push-circle-mor}). 
If $F$ is a sheaf, $\circ$-flat over $Y$, then
$(\id_X\times f)_*F$ is $\circ$-flat over $Y'$ by Lemma \ref{flat-lem}(i).
It is easy to see that for such $F$ and for a quasicoherent sheaf $G$ on $Y\times Z$ the morphism
\eqref{fl-af-push-forward-mor} induces a morphism
\begin{equation}\label{push-forward-circle-isom2}
(\id_X\times f)_*F\ov{\circ}_{\BB_f} (f\times\id_Z)_*G\to F\circ_Y G.
\end{equation}
We claim that in fact this is an isomorphism. Indeed, note that the functor of $G$ in the right-hand side
(resp., the left-hand side) is exact (resp., right exact). 
Also, for every quasicoherent sheaf $G$ on $Y\times Z$
the natural map $(f\times\id_Z)^*(f\times\id_Z)_*G\to G$ is surjective. Thus, it is enough to check
that \eqref{push-forward-circle-isom2} is an isomorphism for $G=(f\times\id_Z)^*G'$, where
$G'$ is a sheaf on $Y\times Z'$. To this end we use the isomorphism
$$\BB_f\circ_{Y'}G'\simeq\left((f\times f)_*\De_*\OO_Y\right)\circ_{Y'}G'\simeq
(f\times\id_Z)_*\left(\De_*\OO_Y\circ_Y (f\times\id_Z)^*G'\right)\simeq 
(f\times\id_Z)_*(f\times\id_Z)^*G'
$$
that follows from \eqref{push-circle-eq}. Together with Lemma \ref{ass-circ-alg-lem}(i) this leads to
$$(\id_X\times f)_*F\ov{\circ}_{\BB_f} (f\times\id_Z)_*(f\times\id_Z)^*G'\simeq (\id_X\times f)_*F
\circ_{Y'} G,$$
which is isomorphic to $F\circ_Y(f\times\id_Z)^*G'$ by Lemma \ref{circle-push-pull-lem}(i).
Note that the isomorphism
\eqref{push-forward-circle-isom2} for $Z=S$ implies that $(\id_X\times f)_*F$ is $\circ$-flat over $\BB_f$
(since every $\BB_f$-module on $Y'$ is of the form $f_*G$, where $G$ is a quasicoherent sheaf on $Y$). From the considered case we can immediately extend the isomorphism
\eqref{push-forward-circle-isom} to the case when $F$ is a bounded above complex of $\circ_Y$-flat
sheaves on $X\times Y$ and $G$ is arbitrary complex in $D_{qc}(Y\times Z)$.
Finally, if $F$ is unbounded we can consider the $\circ_Y$-flat resolutions of the truncations
$\tau_{\le n}F$ and pass to the limit using the isomorphism \eqref{lim-ker-alg-eq}.
\ed

\begin{prop}\label{push-forward-prop} 
Let $f:Y\to Y'$ be a flat affine morphism,
$\BB$ a pure kernel algebra over $Y$, and let
$\BB':=(f\times f)_*\BB$ be the corresponding pure kernel algebra over $Y'$.
Let also $\AA$ (resp., $\CC$) be a pure kernel algebra over $X$ (resp., $Z$).

\noindent
(i) If $F$ is an $\AA\Box\BB^{opp}$-module, $\circ$-flat over $Y$, then $(\id_X\times f)_*F$ is
$\circ$-flat over $Y'$ and for $G\in\BB\Box\CC^{opp}-\mod$ we have a natural isomorphism
\begin{equation}\label{push-forward-ker-isom}
(\id_X\times f)_*F\ov{\circ}_{\BB'}(f\times\id_Z)_*G\wt{\to} F\ov{\circ}_{\BB} G
\end{equation}
In particular, if $F$ is
$\circ$-flat over $\BB$, then $(\id_X\times f)_*F$ is
$\circ$-flat over $\BB'$. 

\noindent
(ii) If $\BB$ is left (resp., right) admissible then so is $\BB'$. In this case for  
$F\in D(\AA\Box\BB^{opp}-\mod)$ and $G\in D(\BB\Box\CC^{opp}-\mod)$ we have
a natural isomorphism
\begin{equation}\label{der-push-forward-ker-isom}
(\id_X\times f)_*F\circ_{\BB'}(f\times\id_Z)_*G\simeq F\circ_{\BB} G.
\end{equation}
\end{prop}

\Pf . (i) We have already seen that $(\id_X\times f)_*F$ is $\circ$-flat over $\BB_f$ (and over $Y'$) in Lemma \ref{simple-push-forward-lem}.
Note also that by the same Lemma, $\BB'=(\id_X\times f)_*(f\times\id_Y)_*\BB$ 
is $\circ$-flat over $\BB_f$ (since $(f\times\id_Y)_*\BB$ is $\circ$-flat over $Y$ by Lemma 
\ref{flat-lem}(i)).
Therefore, using Lemma \ref{subalg-circ-lem} 
we can represent the left-hand side of \eqref{push-forward-ker-isom} 
as the cokernel of
$$(\id_X\times f)_*F\circ_{\BB_f}\BB'\circ_{\BB_f} (f\times\id_Z)_*G\to
(\id_X\times f)_*F\circ_{\BB_f}\circ (f\times\id_Z)_*G.$$
Here the target is isomorphic to $F\circ_{Y} G$ by Lemma \ref{simple-push-forward-lem}.
Also, using the same Lemma, we can rewrite
$$(\id_X\times f)_*F\circ_{\BB_f}(f\times f)_*\BB\circ_{\BB_f} (f\times\id_Z)_*G\simeq
(\id_X\times f)_*F\circ_{\BB_f}(f\times\id_Y)_*\BB\circ_Y G\simeq F\circ_Y\BB\circ_Y G,$$
so that the abovve map is identified with the map
$F\circ_Y\BB\circ_Y G\to F\circ_Y G$ used to define $F\ov{\circ}_\BB G$.
This gives the isomorphism \eqref{push-forward-ker-isom}.

Since the functor $(f\times\id_Z)_*:\BB\Box\CC^{opp}-\mod\to\BB'\Box\CC^{opp}-\mod$ is an equivalence (see Lemma \ref{push-alg-lem}), the second assertion follows.

\noindent
(ii) Let $\BB\to\MM^{\bullet}$ be a quasi-isomorphism with a bounded complex
of $\BB\Box\BB^{opp}$-modules on $Y\times Y$, such that each term $\MM^n$ is 
$\circ$-flattening over $\BB$ on the left. We claim that $\BB'\to (f\times f)_*\MM^{\bullet}$
is a $\circ$-flattening resolution for $\BB'$. Indeed, let us check that $(f\times f)_*\MM^n$ is $\circ$-flattening
over $\BB'$ on the left. First, we note that $(f\times f)_*\MM^n$ is $\circ$-flat over $Y'$ (on the right)
by Lemma \ref{flat-lem}(i).
Next, given a flat quasicoherent sheaf $V$ on $Y'$ we have to check
that $(f\times f)_*\MM^n\circ_{Y'}V$ is $\circ$-flat over $\BB'$. Using \eqref{push-circle-eq} we obtain
$$(f\times f)_*\MM^n\circ_{Y'}V\simeq (f\times\id_Y)_*(\MM^n\circ_Y f^*V)$$
which is $\circ$-flat over $\BB'$ by part (i).
Finally, the isomorphism \eqref{der-push-forward-ker-isom} follows
from \eqref{push-forward-ker-isom} by taking a $q-\circ_{\BB}$-flat resolution of $G$ (resp., $F$). 
\ed

\subsection{Extending an equivalence of derived categories to modules over kernel algebras}
\label{equiv-ker-sec}


Recall that a pair of adjoint functors $(\Phi,\Psi)$, where $\Phi:\CC\to\CC'$, $\Psi:\CC'\to \CC$,
can be characterized by two morphisms of functors $\a:\Phi\Psi\to\Id_{\CC'}$ and 
$\b:\Id_{\CC}\to\Psi\Phi$, such that the compositions $\Phi\to\Phi\Psi\Phi\to\Phi$ and 
$\Psi\to\Psi\Phi\Psi\to\Psi$ are the identity morphisms. Furthermore, $\Phi$ is fully faithful (resp.,
equivalence) if and only if $\b$ is an isomorphism (resp., $\a$ and $\b$ are isomorphisms).
For functors given by kernels we
can consider the following version of this picture.

\begin{defi}
Let $X$ and $Y$ be flat $S$-schemes. Assume that we have 
objects $\PP\in D_{qc}(X\times Y)$ and $\QQ\in D_{qc}(Y\times X)$ 
equipped with a pair of morphisms
\begin{equation}\label{ker-eq1}
\a:\PP\circ_Y\QQ\to\De_*\OO_X,
\end{equation}
\begin{equation}\label{ker-eq2}
\b:\De_*\OO_Y\to\QQ\circ_X\PP,
\end{equation}
such that the compositions
$$\PP\stackrel{\id_{\PP}\circ\b}{\to} \PP\circ_Y\QQ\circ_X\PP\stackrel{\a\circ\id_{\PP}}{\to}\PP$$
$$\QQ\stackrel{\b\circ\id_{\QQ}}{\to} \QQ\circ_X\PP\circ_Y\QQ\stackrel{\id_{\QQ}\circ\a}{\to}\QQ$$
are the identity morphisms. In this situation we say that $(\PP,\QQ;\a,\b)$ is an {\it adjoint kernel
data}.
Consider the corresponding functors
$$\Phi_{\PP}: D_{qc}(Y)\to D_{qc}(X): G\mapsto\PP\circ_Y G,$$
$$\Phi_{\QQ}:D_{qc}(X)\to D_{qc}(Y): F\mapsto\QQ\circ_X F.$$
Then Lemma \ref{circle-prop-lem} implies that $(\Phi_{\PP},\Phi_{\QQ})$ is an adjoint pair of functors.
\end{defi}

\begin{ex} Assume $Y$ is smooth and proper over $S$ of relative dimension $n$, 
and let $\PP$ be a perfect complex on 
$X\times Y$. Then by the duality theory
we have a natural adjoint kernel data $(\PP,\QQ)$ with
$$\QQ=\si_*(\PP^{\vee}\ot p_Y^*\om_Y)[n],$$
where $\si:Y\times X\to X\times Y$ is the permutation of factors, $p_Y:X\times Y\to Y$ is the projection.
\end{ex}

Now let $\AA$ be a kernel algebra over $X$. Set
$\BB=\QQ\circ_X \AA \circ_X\PP$. Then we have a natural morphism
$$\BB\circ_Y \BB\simeq\QQ\circ_X \AA\circ_X\PP\circ_Y\QQ\circ_X \AA\circ_X\PP\to
\QQ\circ_X (\AA\circ_X \AA)\circ_X\PP,$$
induced by \eqref{ker-eq1}.
Therefore, the product $\mu:\AA\circ_X \AA\to \AA$ induces an associative operation
$\BB\circ_Y \BB\to \BB$. Similarly, the unit $u:\De_*\OO_X\to \AA$ together with the
map \eqref{ker-eq2} induce a unit morphism
$\De_*\OO_Y\to \BB$. One can check that the unit axioms for this morphism follow from the
compatibilities between $\a$ and $\b$ that we assumed.
Thus, we get a structure of a kernel algebra on $\BB$.

Note also that \eqref{ker-eq1} induces morphisms
\begin{equation}\label{til-P-Q-eq}
\begin{array}{l}
\PP\circ_Y \BB\simeq \PP\circ_Y\QQ\circ_X \AA\circ_X\PP\to \AA\circ_X\PP,\\
\BB\circ_Y\QQ\simeq\QQ\circ_X\AA\circ_X\PP\circ_Y\QQ\to\QQ\circ_X \AA,
\end{array}
\end{equation}
that are compatible with the products (resp., units) on $\BB$ and $\AA$.

Assume in addition that $\AA$ and $\BB$ are pure kernel algebras,
and that each of the complexes $\PP$ and $\QQ$ is concentrated in single cohomological degree. 
Then by Lemma \ref{flat-lem}(ii), $\wt{\PP}:=\AA\circ_X\PP$ and $\wt{\QQ}:=\QQ\circ_X\AA$ are also concentrated
in single cohomological degree. Using morphisms \eqref{til-P-Q-eq} we get the 
commuting left $\AA$-action and right $\BB$-action on $\wt{\PP}$ (resp., left $\BB$-action and right $\AA$-action on $\wt{\QQ}$).
Therefore, using Lemma \ref{prod-mod-lem} we can view
$\wt{\PP}$ (resp., $\wt{\QQ}$) as an object of
$D(\AA\Box\BB^{opp}-\mod)$ (resp., $D(\BB\Box\AA^{opp}-\mod)$).
Note that the unit in $\AA$ induces morphisms
$$\iota_{\PP}:\PP\to\wt{\PP},\ \ \ \iota_{\QQ}:\QQ\to\wt{\QQ}$$
in $D_{qc}(X\times Y)$ and $D_{qc}(Y\times X)$.

\begin{thm}\label{PR-thm} Assume $S$ is semi-separated.
Let $X$ and $Y$ be semi-separated schemes, flat over $S$, and let
$\PP\in D_{qc}(X\times Y)$, $\QQ\in D_{qc}(Y\times X)$ be objects, each concentrated in
single cohomological degree, equipped with a pair of morphisms
\eqref{ker-eq1} and \eqref{ker-eq2}, such that $(\PP,\QQ,\a,\b)$ is an adjoint kernel data.

\noindent
(i) Let $\AA$ be an admissible kernel algebra over $X$.
Assume that the corresponding kernel algebra 
$\BB=\QQ\circ_X \AA\circ_X\PP$ over $Y$
is admissible (in particular, pure). Define
$\wt{\PP}\in D(\AA\Box\BB^{opp}-\mod)$ and 
$\wt{\QQ}\in D(\BB\Box\AA^{opp}-\mod)$ as above, and consider the functors
$$\Phi_{\wt{\PP}}: D(\BB-\mod)\to D(\AA-\mod): G\mapsto \wt{\PP}\circ_{\BB} G,$$
$$\Phi_{\wt{\QQ}}: D(\AA-\mod)\to D(\BB-\mod): F\mapsto \wt{\QQ}\circ_{\AA} F.$$ 
Then the natural map in $D_{qc}(Y\times Y)$ given as the composition
\begin{equation}\label{QP-iota-eq}
\BB=\QQ\circ_X\AA\circ_X\PP\stackrel{\iota_{\QQ}\iota_{\PP}}{\to}
\wt{\QQ}\circ_X\AA\circ_X\wt{\PP}\to\wt{\QQ}\circ_{\AA}\wt{\PP}
\end{equation}
lifts to an isomorphism in $D(\BB\Box\BB^{opp}-\mod)$. 
Hence, we have the induced isomorphism of functors
\begin{equation}\label{QP-fun-eq}
\Phi_{\wt{\QQ}}\circ \Phi_{\wt{\PP}}\simeq\Id.
\end{equation}
If in addition $\a$ and $\b$ are isomorphisms then $\Phi_{\wt{\PP}}$ and
$\Phi_{\wt{\QQ}}$ are mutually quasi-inverse equivalences.

\noindent
(ii) Let $\AA_1\to \AA$ be a homomorphism of kernel algebras over $X$, and consider the induced homomorphism of kernel algebras
$\BB_1=\QQ\circ_X \AA_1\circ\PP\to\BB=\QQ\circ_X\AA\circ_X\PP$.
Suppose $\AA,\AA_1,\BB,\BB_1$ are all admissible. Then the functors 
$\Phi_{\wt{\QQ}}$ (resp., $\Phi_{\wt{\PP}}$) constructed in (i) for
the pairs $(\AA,\BB)$ and $(\AA_1,\BB_1)$ commute with the  
corresponding restriction (resp., induction) functors between the derived categories of modules.

\noindent
(iii) Assume in addition that $X$ and $Y$ are proper over $S$,
$\PP$ and $\QQ$ are coherent sheaves of finite $\Tor$-dimension (up to a shift),
and $\AA$ is a finite admissible kernel algebra over $X$ such that $\BB=\QQ\circ_X\AA\circ_X\PP$
is admissible. Then $\BB$ is also finite and the functor $\Phi_{\wt{\QQ}}$
restricts to a functor from $D^b(\AA-\mod^c)$ to $D^b(\BB-\mod^c)$. 
In particular, when $\a$ and $\b$ are isomorphisms
these two categories are equivalent.
\end{thm}

\Pf . (i) Using the properties of the convolution operation for modules we obtain natural isomorphisms
\begin{equation}\label{QPB-eq1}
\wt{\QQ}\circ_{\AA}\wt{\PP}\simeq \wt{\QQ}\circ_{\AA}(\AA\circ_X\PP)\simeq
\wt{\QQ}\circ_X\PP\simeq\QQ\circ_X\AA\circ_X\PP=\BB
\end{equation}
in $D(\BB\Box\De_*\OO_Y-\mod)$ and
\begin{equation}\label{QPB-eq2}
\wt{\QQ}\circ_{\AA}\wt{\PP}\simeq (\QQ\circ_X\AA)\circ_{\AA}\wt{\PP}\simeq
\QQ\circ_X\wt{\PP}\simeq\QQ\circ_X\AA\circ_X\PP=\BB
\end{equation}
in $D(\De_*\OO_Y\Box\BB^{opp}-\mod)$.
In particular, $\wt{\QQ}\circ_{\AA}\wt{\PP}$ is concentrated in degree $0$. Also, it is easy to check that
both isomorphisms \eqref{QPB-eq1} and \eqref{QPB-eq2} are compatible with the map
\eqref{QP-iota-eq}. Hence, we get the required isomorphism in $\BB\Box\BB^{opp}-\mod$.
The isomorphism of functors \eqref{QP-fun-eq} follows from this by the associativity of the
convolution (see Proposition \ref{ass-ker-alg-prop}).

If both $\a$ and $\b$ are isomorphisms then we can reverse the roles of $(X,\AA)$ and $(Y,\BB)$
(resp., of $\PP$ and $\QQ$). This will lead to an isomorphism of functors
$\Phi_{\wt{\PP}}\circ\Phi_{\wt{\QQ}}\simeq \Id$ that together with \eqref{QP-fun-eq} proves that
$\Phi_{\wt{\PP}}$ and $\Phi_{\wt{\QQ}}$ are mutually quasi-inverse.

\noindent
(ii) Let $\wt{\PP}$ and $\wt{\QQ}$ have the same meaning as in (i), and let
$\wt{\PP}_1:=\AA_1\circ_X\PP$ and $\wt{\QQ}_1:=\QQ\circ_X\AA_1$ be similar objects
for $\AA_1$. Then we have an isomorphism $\wt{\QQ}\simeq\wt{\QQ}_1\circ_{\AA_1}\AA$
in $D(\BB\Box\AA^{opp}-\mod)$. Therefore, for $F\in D(\AA-\mod)$ we get
$$\Phi_{\wt{\QQ}}(F)=\wt{\QQ}\circ_{\AA} F\simeq(\wt{\QQ}_1\circ_{\AA_1}\AA)\circ_{\AA} F\simeq
\wt{\QQ}_1\circ_{\AA_1} F=\Phi_{\wt{\QQ}_1}(F),$$
which gives the compatibility of $(\Phi_{\wt{\QQ}},\Phi_{\wt{\QQ}_1})$ with the restriction functors
$D(\AA-\mod)\to D(\AA_1-\mod)$ and $D(\BB-\mod)\to D(\BB_1-\mod)$. On the other hand,
using the isomorphism $\wt{\PP}\simeq \AA\circ_{\AA_1}\wt{\PP}_1$,
for $G\in D(\BB_1-\mod)$ we obtain
$$\Phi_{\wt{\PP}}(\BB\circ_{\BB_1}G)=\wt{\PP}\circ_{\BB}(\BB\circ_{\BB_1}G)\simeq
\wt{\PP}\circ_{\BB_1}G\simeq \AA\circ_{\AA_1}\wt{\PP_1}\circ_{\BB_1}G=\AA\circ_{\AA_1}
\Phi_{\wt{\PP}_1}(G),$$
which gives the required compatibility with the induction functors.

\noindent
(iii) Our assumptions imply that $\BB$ is obtained by the proper push-forward from a bounded
complex with coherent cohomologies, hence, it is a coherent sheaf itself.
By part (ii), the functor $\Phi_{\wt{\QQ}}:D(\AA-\mod)\to D(\BB-\mod)$ commutes
with the forgetful functors to quasicoherent sheaves, so it preserves
the condition for a complex to have coherent cohomology sheaves. 
By Lemma \ref{coh-ker-lem}, we obtain a functor $D^-(\AA-\mod^c)\to D^-(\BB-\mod^c)$
commuting with forgetful functors. Since $\QQ$ has finite $\Tor$-dimension, it preserves
boundedness of cohomology. The last assertion follows from (i).
\ed


\section{Kernel representations and derived equivalences}
\label{kerrep-sec}

In this section we develop a framework for
constructing examples of the Fourier-Mukai
dual kernel algebras. 
The main idea is that some categories of twisted sheaves on global quotient stacks $X/G$ can be identified with modules over kernel algebras that are constructed from 
the corresponding action of $G$ on the category $D_{qc}(X)$. 
We introduce the notion of a kernel representation of a group $G$ over $X$
 that should be viewed as a refined version of an action of $G$ on $D_{qc}(X)$.
The point is that given such a kernel representation and a derived equivalence from $X$ to $Y$,
one gets a kernel representation of $G$ over $Y$. This helps to compute the Fourier-Mukai
transforms of the associated kernel algebras. An example of such computation in the context of
abelian varieties is Theorem \ref{twist-eq-thm}.

We denote by $S$ a fixed base scheme, and as before, all products are fibered over $S$.
{\it We also assume all schemes and formal schemes to be flat over $S$}. 
By a formal group scheme $G$ over $S$ we mean a group object in the category of formal
schemes over $S$ (flat over $S$). Sometimes, we will assume that a formal group scheme $G$ is
of {\it ldu}-pseudo-finite type (resp., {\it ldu}-pseudo-finite)
over $S$---we refer to this as {\it ldu-pft} (resp., {\it ldu-pf}) over $S$.

For a field $k$ of characteristic zero we will call simply {\it formal $k$-groups}  
({\it $k$-groupes formels} in terminology of \cite{Lau}),
commutative formal groups $G$ over $k$ such that
(i) the natural action of $\Gal(\ov{k}/k)$ on $G(\ov{k})$ factors through some quotient $\Gal(k'/k)$
with $k'\sub\ov{k}$ a finite Galois extension of $k$, (ii) the group $G(\ov{k})=G(k')$ is finitely generated,
and (iii) the topological $k$-algebra $\OO_{G,0}$ is a quotient of $k[[t_1,\ldots,t_n]]$ for some $n$.
The Cartier duality (see \cite{SGA3}, exp. ${\rm VII_B}$, (2.2.2), or \cite{Dem}, ch.II) 
gives an equivalence between 
the category of formal $k$-groups and the category of affine commutative algebraic
groups over $k$.
We will also use the Cartier duality for finite flat group schemes over an arbitrary base $S$.

Starting from section \ref{Cartier-sec},
whenever formal $k$-groups are mentioned it is assumed that $S=\Spec(k)$.



\subsection{Convolution of sheaves parametrized by formal schemes}
\label{formal-conv-sec}

We need a version of the convolution operation considered in section \ref{circle-sec} involving
formal schemes. More precisely, formal schemes will appear only as parameter spaces, so that the actual ``integration" is performed only along the direction of a usual scheme.
We start with the data consisting of 
three schemes $X,Y,Z$ and two formal schemes 
$P$ and $Q$. 
Let us denote by $D_{P,\clim}(X\times P\times Y)\sub D_{qc}(X\times P\times Y)$
the full subcategory of objects $F$ such that for some open covering $(U_i)$ of $P$
the restriction of $F$ to $X\times U_i\times Y$ is in $D_{\clim}(X\times U_i\times Y)$.
Assuming that $Y$ is proper, we can define the convolution operation
$$D_{P,\clim}(X\times P\times Y)\times D_{Q,\clim}(Y\times Q\times Z)\to 
D_{P\times Q,\clim}(X\times P\times Q\times Z)$$
sending $(K,L)$ to 
$$K\circ_Y L:=p_{1245*}(p_{123}^*K\otimes p_{345}^*L),$$
where $p_I$ denote projections from the product $X\times P\times Y\times Q\times Z$.
Note that 
$K\circ_Y L$ is in $D_{P\times Q,\clim}$ by Proposition \ref{formal-pull-push-prop}(vi).
More generally, the above operation is defined when $Y$ is not necessarily proper but
$K$ is the push-forward from a closed formal subscheme proper over $X\times P$
---we will say in this case that the support of $K$ is proper over $X\times P$ (alternatively,
it is enough to require the support of $L$ to be proper over $Q\times Z$).
Using appropriate versions of the  projection and base change formulae
for formal schemes (see Theorems \ref{formal-proj-formula} and \ref{formal-base-change}(i)) we can easily check the following two properties of this operation:

{\bf Associativity}: if $T$ (resp., $R$) is another scheme (resp., 
formal scheme) then for $M\in D_{R,\clim}(Z\times R\times T)$ one has a natural
isomorphism
$$(K\circ_Y L)\circ_Z M\simeq K\circ_Y (L\circ_Z M).$$
provided the support of $L$ 
is proper over $Y\times Q$ and over $Q\times Z$.

{\bf Change of parameter spaces}: if $\phi:P'\to P$ and $\psi:Q'\to Q$ are maps of 
formal schemes then we have
$$\left((\id_X\times \phi\times\id_Y)^*K\right)\circ_Y \left((\id_Y\times \psi\times\id_Z)^*L\right)\simeq 
(\id_X\times \phi\times\psi\times\id_{Z})^*(K\circ_Y L)
$$
provided the support of $L$ is proper over $Q\times Z$ (or the support of $K$ is proper over 
$X\times P$).


\subsection{Kernel representations and convolution algebras}
\label{ker-rep-alg-sec}

The following definition provides a natural framework for working with "nice" homomorphisms
from a group scheme (resp., formal group scheme)
to the group of autoequivalences of the 
derived category of quasicoherent sheaves on a scheme.

\begin{defi} 
Let $G$ be a monoid in the category of $S$-schemes, $X$ an $S$-scheme (recall
that flatness over $S$ is always assumed). 

\noindent (i)
A {\it kernel quasi-representation of $G$ over $X$}
is an object $V\in D_{qc}(X\times G\times X)$
equipped with a morphism
$$\mu:V\circ_X V\to (\id_X\times m \times\id_X)^*V$$
over $X\times G\times G\times X$, where $m:G\times G\to G$ is the group law,
satisfying the natural associativity condition over $X\times G\times G\times G\times X$,
and a morphism
$$u:\De_*\OO_X\to (\id_X\times e\times \id_X)^*V,$$
where $e:S\to G$ is the unit for $G$, compatible with $\mu$ in a natural way.
The associativity condition amounts to the commutativity of the following diagram:
\begin{diagram}
V\circ_X V\circ_X V &\rTo{\mu\circ_X\id} & (\id_X\times m\times\id_X)^*V\circ_X V\simeq&
(\id_X\times m\times\id_{G\times X})^*(V\circ_X V)\\
\dTo{\id\circ_X\mu}&&&\dTo{\mu'}\\
V\circ_X(\id_X\times m\times\id_X)^*V&\simeq
(\id_{X\times G}\times m\times\id_X)^*(V\circ_X V)&\rTo{\mu''}&(\id_X\times mm\times \id_X)^*V
\end{diagram}
where $mm:G^3\to G$ sends $(g_1,g_2,g_3)$ to $g_1g_2g_3$, the maps $\mu'$ and $\mu''$
are induced by $\mu$, the isomorphisms are particular cases of \eqref{pull-circle-eq}.
The compatibility of the unit morphism $u$ with $\mu$ means that the composition
$$V\simeq \De_*\OO_X\circ_X V\stackrel{u\circ_X\id}{\to} (\id_X\times e\times\id_X)^*V\circ_X V\to V,$$
where the second arrow is induced by $(\id_X\times e\times\id_G\times\id_X)^*\mu$
(resp., similar composition with $\id\circ_X u$) equals the identity map.

\noindent (ii)  
A {\it kernel representation} is a kernel quasi-representation such that
$\mu$ and $u$ are isomorphisms. 

\noindent
(iii) Let $G$ be a monoid in the category of formal schemes over $S$, and let
$X$ be an $S$-scheme. Then the definitions (i) and (ii) still make sense
with $V\in D_{G,\clim}(X\times G\times X)$ provided the support of $V$ is proper over $X\times G$
and over $G\times X$ (see section \ref{formal-conv-sec}).

\noindent
(iv) If $\rho:G'\to G$ is a homomorphism of monoids and $V$ is a kernel quasi-representation of
$G$ over $X$ then we define the {\it restriction of $V$ to $G'$} by 
$$\rho^*V:=(\id_X\times\rho\times\id_X)^*V$$
---it is a kernel quasi-representation of $G'$ over $X$. The same construction works when $G$ and $G'$
are formal group schemes 
provided $V\in D_{\clim}(X\times G\times X)$ has the support proper over $X\times G$
and over $G\times X$.
\end{defi}

If $G$ is a discrete group then to give a kernel representation of $G$ amounts to the
data of  kernels $V_g\in D_{qc}(X\times X)$ parametrized by $g\in G$ along with the isomorphisms
$V_{g_1}\circ V_{g_2}\simeq V_{g_1g_2}$ satisfying the natural axioms (in particular,
the functors corresponding to $V_g$ and $V_{g^{-1}}$ are quasi-inverse of each other). Thus, we
have an action of $G$ on $D_{qc}(X)$ by the corresponding autoequivalences (in the strict sense
as defined in \cite{Verdier}). A kernel quasi-representation in the case of a discrete group
$G$ can be viewed as a $G$-graded kernel algebra $\oplus_{g\in G} V_g$. 
One can try to imitate this construction for general $G$
by associating with a kernel quasi-representation
a kernel algebra via some kind of ``integration" over $G$. Below we will show that this
is possible provided one fixes an appropriate analogue of ``measure" on $G$. 

The passage from kernel representations to kernel algebras goes through an intermediate notion
that we call a {\it convolution algebra}. We will first construct kernel algebras associated with
convolution algebras and then explain the connection with kernel representations
in Proposition \ref{ker-G-rep-prop}.

\begin{defi} (i) Let $X$ be a scheme and let $M\stackrel{(s,t)}{\to} X\times X$ be a 
monoid in the category of $X\times X$-schemes, such that the maps $s$ and $t$ from $M$ to $X$ are flat.
Let $m: M\times_X M\to M$ and $e:X\to M$ denote the product and the unit maps.
A {\it convolution algebra on} $M$ is an object $\VV\in D_{qc}(M)$
equipped with a product morphism
$$\mu:m_*(p_1^*\VV\otimes p_2^*\VV)\to\VV$$
and a unit morphism 
$$e_*\OO_X\to\VV$$
such that 

\noindent
(a) the diagram
\begin{diagram}
m_*(\id_M\times m)_*(p_1^*\VV\otimes p_2^*\VV\otimes p_3^*\VV)&\simeq
m_*(m\times\id_M)_*(p_1^*\VV\otimes p_2^*\VV\otimes p_3^*\VV)\rTo{\mu\ot\id}&
m_*(p_1^*\VV\otimes p_2^*\VV)\\
\dTo{\id\ot\mu} & &\dTo{\mu}\\
m_*(p_1^*\VV\otimes p_2^*\VV)&\rTo{\mu}&\VV
\end{diagram}
is commutative, where we use the base change isomorphism for the cartesian diagram
\begin{diagram}
M\times_X M\times_X M &\rTo{p_{12}}& M\times_X M\\
\dTo{m\times\id_M}&&\dTo{m}\\
M\times_X M&\rTo{p_1}& M
\end{diagram} 
along with the projection formula for $m\times\id_X$ 
(plus similar isomorphisms for $\id_M\times m$); 

\noindent
(b) the composition
$$\VV\simeq m_*(e\times\id_M)_*\VV\simeq 
m_*((e\times\id_M)_*\OO_M\otimes p_2^*\VV)\simeq m_*(p_1^*e_*\OO_X)\otimes p_2^*\VV)\to
m_*(p_1^*\VV\otimes p_2^*\VV)\stackrel{\mu}{\to}\VV$$ 
(along with a similar composition going through
$m_*(p_1^*\VV\otimes p_2^*(e_*\OO_X))$) is the identity. Note that here we use the projection formula
for the map $e\times\id_M: M\to M\times_X M$ and the base change formula for the cartesian diagram
\begin{diagram}
M &\rTo{s}& X\\
\dTo{e\times\id_M} && \dTo{e}\\
M\times_X M &\rTo{p_1}& X
\end{diagram}

\noindent
(ii) Now let $X$ be a formal scheme, $M\stackrel{(s,t)}{\to} X\times X$ a monoid in the
category of formal schemes over $X\times X$ ({\it formal monoid over} $X\times X$). Assume in addition 
that the maps $s$ and $t$ are flat and of {\it ldu}-pseudo-finite type, and $m$ is {\it ldu}-quasicompact. A {\it convolution algebra on} $M$ is
an object $\VV\in D_{qct}(M)$ equipped with the same structure as before with $m_*$ replaced
by $m_{\sigma *}$ (where we use 
Theorems \ref{formal-proj-formula} and \ref{formal-base-change}(ii),(iii) to formulate the above conditions). Note that in the case when $M$ is a groupoid
the condition for $m$ to be {\it ldu}-quasicompact follows from the assumption that $s$ is
of {\it ldu}-pseudo-finite type.
\end{defi}

It is clear that for $M=X\times X$, where $X$ is a scheme, we get exactly the notion of a kernel algebra.
On the other hand, for $M=X$ we get the notion of an algebra in $D_{qc}(X)$ with respect
to the tensor product.
In the case when $X=S$ (the base scheme) and 
$M=G$ is a group scheme over $S$, 
the {\it convolution algebra on $G$}, is an algebra object in $D_{qc}(G)$ with respect to the Pontryagin product. Similarly, if $G$ is an {\it ldu-pft} formal group scheme over $S$ (i.e.,
of {\it ldu}-pseudo-finite type over $S$) then 
we have a notion of a convolution algebra structure on an object of $D_{qct}(G)$. 

\begin{lem}\label{push-convolution-lem} 
Let $f:M\to M'$ is the 
morphism of monoids over $X\times X$ satisfying the conditions of the above definition,
and let $\VV$ be a convolution
algebra on $M$. Then $f_*\VV$ has a natural structure of a convolution algebra
on $M'$. If $f':M'\to M''$ is another such morphism then the natural isomorphism
$(f'f)_*(M)\simeq f'_*f_*(M)$ is compatible with the convolution algebra structures.
Similar assertions hold in the case of formal monoids, where $f$ is assumed to be
{\it ldu}-quasicompact and $f_*$ is replaced with $f_{\sigma *}$
\end{lem}

\Pf . Set $\VV'=f_*\VV$, 
and let $m':M'\times_X M'\to M'$ denote the product map for $M'$. 
First, let us observe that using the adjoint pair $((f\times f)^*, (f\times f)_*)$ and the canonical
maps $f^*\VV'\to\VV$ we get a morphism on $M'\times M'$
\begin{equation}\label{ffVV-eq}
p_1^*\VV'\otimes p_2^*\VV'\to (f\times f)_*(p_1^*\VV\otimes p_2^*\VV).
\end{equation}
We claim that it is an isomorphism. Indeed, applying the base change formula
to the cartesian square
\begin{diagram}
M\times M' &\rTo{f\times\id}& M'\times M'\\
\dTo{p_1}&&\dTo{p_1}\\
M &\rTo{f}& M'
\end{diagram}
we get an isomorphism $p_1^*\VV'\simeq(f\times\id_{M'})_*p_1^*\VV$ on $M'\times M'$.
Together with the projection formula this gives
$$p_1^*\VV'\otimes p_2^*\VV'\simeq (f\times\id_{M'})_*(p_1^*\VV\otimes p_2^*\VV').$$
Next, we apply the base change formula to the cartesian square
\begin{diagram}
M\times M &\rTo{\id\times f}& M\times M'\\
\dTo{p_2}&&\dTo{p_2}\\
M &\rTo{f}& M'
\end{diagram}
to get an isomorphism $p_2^*\VV'\simeq(\id_M\times f)_*p_2^*\VV$. Combined with the previous
isomorphism and using the projection formula again we obtain an isomorphism
$$p_1^*\VV'\otimes p_2^*\VV'\simeq (f\times\id_{M'})_*(p_1^*\VV\otimes (\id_M\times f)_*p_2^*\VV)
\simeq (f\times f)_*(p_1^*\VV\otimes p_2^*\VV),$$
that is given by the map \eqref{ffVV-eq}.
Therefore,
$$m'_*(p_1^*\VV'\otimes p_2^*\VV')\simeq m'_*(f\times f)_*(p_1^*\VV\otimes p_2^*\VV)\simeq
f_*m_*(p_1^*\VV\ot p_2^*\VV)$$
over $M'\times_X M'$ that allows to view $f_*\mu$ as a product structure on $\VV'$.
The associativity of this product can be easily deduced from the associativity of $\mu$ using the 
canonical isomorphism on $M'\times M'\times M'$
$$p_1^*\VV'\otimes p_2^*\VV'\otimes p_3^*\VV'\wt{\to} 
(f\times f\times f)_*(p_1^*\VV\otimes p_2^*\VV\otimes p_3^*\VV)$$
defined similarly to \eqref{ffVV-eq}. The case of formal monoids is very similar.
\ed

For us the most important case of the above construction is when $M'=X\times X$ and $f$ is the
structure morphism $(s,t):M\to X\times X$.


\begin{defi} Let $X$ be a scheme,
$M$ a monoid (resp., formal monoid) in the category of $X\times X$-schemes such that
the structure maps $s,t:M\to X$ are flat (resp, flat and of {\it ldu}-pseudo-finite type; 
in addition assume that $m$ is {\it ldu}-quasicompact). 
We associate with every convolution algebra $\VV$ over $M$ a kernel
algebra over $X$ by setting
\begin{equation}\label{conv-alg-ker-alg}
\AA_M(\VV)=(s,t)_*\VV
\end{equation}
(resp., $\AA_M(\VV)=(s,t)_{\sigma *}\VV$).
\end{defi}
 
The following lemma gives a sufficient condition 
for checking when the kernel algebra of the form $\AA_M(\VV)$ is pure.

\begin{lem}\label{aff-monoid-lem} Assume $X$ and $S$ are semi-separated.
Let $M\stackrel{(s,t)}{\to} X\times X$ be a monoid (resp., formal monoid) over $X\times X$  
satisfying the assumptions of the above definition, and let $\VV$ be a convolution algebra over $X$. 
Assume in addition that both maps $s,t:M\to X$ are affine (resp., {\it ldu}-affine),
and that $\VV$ is a quasicoherent sheaf (resp., torsion quasicoherent sheaf),
flat over $X$ with respect to both $s$ and $t$.  Then the kernel algebra
$\AA_M(\VV)$ is pure.
\end{lem}

\Pf . We only need to check that $(s,t)_*\VV$ (resp., $(s,t)_{\sigma *}\VV$)
is $\circ$-flat over $X$ with respect to both projections.
But this follows immediately from Lemma \ref{aff-flat-lem}. 
\ed

If $G$ is a group scheme (resp., formal group scheme) over $S$
then we can view $X\times G\times X$ as a monoid (resp., formal monoid) over $X\times X$ 
using the product on $G$ and the unit in $G$. 
In the following proposition we denote by $p_i$ (resp., $p_{ij}$) the projection of the product
$X\times G\times X$ onto its factors (resp., double products of factors).

\begin{prop}\label{ker-G-rep-prop} 
Let $G$ be a group scheme (resp., {\it ldu-pft} formal group scheme) over $S$, $X$ an $S$-scheme, 
$V$ a kernel quasi-representation of $G$ over $X$ (recall that in the formal case we assume that 
$V\in D_{G,\clim}(X\times G\times X)$ and the support of $V$ is proper over $X\times G$ and
$G\times X$), $P$ a convolution algebra on $G$. 
Then $V\otimes p_2^*P$ 
has a natural structure of a convolution algebra on $X\times G\times X$. 
Hence, we also have the corresponding kernel algebra over $X$
$$\AA(V,P):=\AA_{X\times G\times X}(V\otimes p_2^*P)=
p_{13*}(V\otimes p_2^*P).$$
(in the formal case we should use $(p_{13})_{\sigma*}$)
\end{prop}

\Pf . We will only consider the case when $G$ is a group scheme (the case of a formal group scheme
is analogous).
Note that the product map for $X\times G\times X$ is the composition of the 
projection $X\times G\times X\times G\times X\to X\times G\times G\times X$ with
the map $(\id_X\times m\times\id_X):X\times G\times G\times X\to X\times G\times X$.
Applying the projection formula for the former map we see that the 
convolution product on $V\otimes p_2^*P$ should correspond to a map
$$(\id_X\times m\times\id_X)_*[(V\circ_X V)\otimes p_2^*P\otimes p_3^*P]\to V\otimes p_2^*P.$$
To construct such map we start with the morphism
$$(V\circ_X V)\otimes p_2^*P\otimes p_3^*P \to
(\id_X\times m\times\id_X)^*V\otimes p_2^*P\otimes p_3^*P$$
on $X\times G\times G\times X$ induced by the map $\mu:V\circ_X V\to (\id_X\times m\times\id_X)^*V$.
Combining the projection formula for the map $\id_X\times m\times\id_X$ with the
flat base change formula for $P\boxtimes P\in D_{qc}(G\times G)$ 
and the cartesian diagram
\begin{diagram}
X\times G\times G\times X&\rTo{p_{23}}& G\times G\\
\dTo{\id\times m\times\id}&&\dTo{m}\\
X\times G\times X&\rTo{p_2}& G
\end{diagram}
we get an isomorphism
$$(\id_X\times m\times\id_X)_*((\id_X\times m\times\id_X)^*V\otimes p_2^*P\otimes p_3^*P)\simeq
V\otimes p_2^*m_*(P\boxtimes P).$$
Hence, we get the required morphism
$$\mu_{V,P}:(\id_X\times m\times\id_X)_*[(V\circ_X V)\otimes p_2^*P\otimes p_3^*P]\to
V\otimes p_2^*m_*(P\boxtimes P)\to
V\otimes p_2^*P$$
over $X\times G\times X$, where the second arrow is induced by the convolution algebra structure on $P$.
The associativity is checked as follows. Let us consider the following diagram with cartesian squares
\begin{diagram}
X\times G\times X\times G\times X\times G\times X&\rTo{}&X\times G\times X\times G\times G\times X 
&\rTo{}& X\times G\times X\times G\times X\\
\dTo{}&&\dTo{}&&\dTo{}\\
X\times G\times G\times X\times G\times X &\rTo{}& X\times G\times G\times G\times X
&\rTo{\id_{X\times G}\times m\times\id_X}& X\times G\times G\times X\\
\dTo{}&&\dTo{\id_X\times m\times\id_{G\times X}}&&\dTo{}\\
X\times G\times X\times G\times X&\rTo{}&X\times G\times G\times X&\rTo{}&X\times G\times X
\end{diagram}
where all the arrows are either projections omitting one of the factors $X$, or maps induced by the product in two consecutive factors $G$. For each vertex $v$ in this diagram 
let $X(v)=X\times G^{n_1}\times X\times G^{n_2}\times\ldots$ denote the corresponding scheme. 
We have a natural object $W(v)\in D(X(v))$ defined as
follows: let $m(v):X(v)\to X\times G\times X\times G\ldots \times X$ 
denote the map induced by the product
maps $G^{n_1}\to G$, $G^{n_2}\to G$, etc. Then $W(v)$ is given by the (derived) tensor product
of pull-backs of $P$ from all the factors of $G$ with the pull-back by $m(v)$ of
the object $p_{123}^*V\otimes p_{345}^*V\otimes\ldots\in D(X\times G\times X\times G\ldots\times X)$.
Next, for every arrow $e:v\to v'$ let $f(e):X(v)\to X(v')$ be the corresponding morphism.
Then we have natural map $f(e)_*W(v)\to W(v')$ induced either by $\mu$ (in the case when
$f(e)$ is the projection omitting one of the factors $X$), or by the convolution product on $P$ (when
$f(e)$ is induced by $m$). Furthermore, we claim that for each of the four small cartesian squares
of our diagram 
\begin{diagram}
v_1&\rTo{e_{12}}&v_2\\
\dTo{e_{13}}&&\dTo{e_{24}}\\
v_3&\rTo{e_{34}}&v_4
\end{diagram}
the induced morphisms
$$f(e_{24}\circ e_{12})_*X(v_1)\to f(e_{24})_*f(e_{12})_*X(v_1)\to f(e_{24})_*X(v_2)\to X(v_4)$$
and
$$f(e_{34}\circ e_{13})_*X(v_1)\to f(e_{34})_*f(e_{13})_*X(v_1)\to f(e_{34})_*X(v_3)\to X(v_4)$$
are the same. Indeed, for the upper left (resp., lower right) square this amounts to the associativity
of $\mu$ (resp., of the convolution product on $P$). For the remaining two squares this can be deduced
from the base change formula. Now the required associativity of $\mu_P$ follows from the similar
compatibility applied to the ambient square in the big diagram above.

The unit for $V\otimes p_2^*P$ is given by the composition
$$(\id_X\times e\times\id_X)_*\De_*\OO_X\stackrel{u}{\to}
(\id_X\times e\times\id_X)_*(\id_X\times e\times\id_X)^*V\simeq
V\ot(\id_X\times e\times\id_X)_*\OO_{X\times X}\to V\otimes p_2^*P,$$
where the second arrow is induced by the unit for $P$.
\ed

Assume that $G$ is a finite group scheme (resp., locally nicely ind-finite formal group scheme) 
over $S$.
Let us consider the dualizing sheaf $\om_{G/S}:=\pi^!\OO_S$, where $\pi:G\to S$ is the projection
(in the formal case we use the functor $\pi^!$ from Proposition \ref{formal-finite-duality-prop}(ii)).
It is easy to see that $\om_{G/S}$ has a natural structure of a convolution algebra on $G$ (see
Theorem \ref{equiv-sh-thm} below for a more general construction).
Thus, in this case we can associate with any quasi-kernel representation $V$ of $G$ the
kernel algebra $\AA(V,\om_{G/S})$ over $X$.

For us the most useful property of kernel (quasi-)representation is that they
can be carried over under derived equivalences in a way compatible with the corresponding
operation on kernel algebras. This follows from the following more general statement.

\begin{prop}\label{ker-equiv-prop}
Let $X$ and $Y$ be  $S$-schemes, 
and let $(\PP,\QQ,\a,\b)$ be adjoint kernel data, where
$\PP\in D_{qc}(X\times Y)$ and $\QQ\in D_{qc}(Y\times X)$
(see section \ref{equiv-ker-sec}).
Let $G$ be a group scheme over $S$, and let $V$ be a kernel quasi-representation
of $G$ over $X$.
Then $\QQ\circ_X V\circ_X\PP$ has a natural structure of
a kernel quasi-representation of $G$ over $Y$. 
This correspondence is compatible with the restriction under a homomorphism of
group schemes $G'\to G$.
Furthermore, if $P$ is a convolution algebra on $G$
then we have an isomorphism of kernel algebras over $Y$:
\begin{equation}\label{ker-equiv-eq}
\AA((\QQ\circ_X V\circ_X\PP), P)\simeq\QQ\circ_X\AA(V, P)\circ_X\PP.
\end{equation}
Similar assertions hold if $G$ is a formal group scheme over $S$, $X$ and $Y$ are proper over $S$,
$V\in D_{G,\clim}(X\times G\times X)$ (in the analog of the last assertion we should require $G$ to
be {\it ldu-pft} over $S$ in order to consider a convolution algebra over it).
\end{prop}

\Pf . Using the associativity of the $\circ$-product, the map $\a:\PP\circ_Y\QQ\to\De_*\OO_X$ and the 
product map $\mu$ for $V$ we obtain a morphism
$$(\QQ\circ_X V\circ_X\PP)\circ_Y(\QQ\circ_X V\circ_X\PP)\to
\QQ\circ_X V\circ_X V\circ_X \PP\to\QQ\circ_X(\id_X\times m\times\id_X)^*V\circ_X\PP.$$
By Lemma \ref{circle-push-pull-lem}(ii), the target of this map is naturally isomorphic to
$(\id_Y\times m\times\id_Y)^*(\QQ\circ_X V\circ_X\PP)$.
This gives a product map for $\QQ\circ_X V\circ_X\PP$.
On the other hand, by Lemma \ref{circle-push-pull-lem}(iii), we have
$$(\id_Y\times e\times\id_Y)^*(\QQ\circ_X V\circ_X\PP)\simeq \QQ\circ_X 
(\id_X\times e\times\id_X)^*V\circ_X\PP.$$
Hence, the unit for $V$ together with the map $\b:\De_*\OO_Y\to\QQ\circ_X\PP$ induce a unit
for $\QQ\circ_X V\circ_X\PP$.
The isomorphism \eqref{ker-equiv-eq} is obtained using the associativity of the $\circ$-product:
\begin{align*}
&\AA((\QQ\circ_X V\circ_X\PP), P)=(\QQ\circ_X V\circ_X\PP)\circ_{G\times X} (P\Box\De_*\OO_X)\simeq\\
&\QQ\circ_X V\circ_{G\times X}(\De_*\OO_G\Box\PP)\circ_{G\times X} (P\Box\De_*\OO_X)\simeq
\QQ\circ_X V\circ_{G\times X}(P\Box\De_*\OO_X)\circ_{G\times X}(\De_*\OO_G\Box\PP)\simeq\\
&\QQ\circ_X (V\circ_{G\times X}(P\Box\De_*\OO_X))\circ_X \PP.
\end{align*}
\ed

\subsection{Equivariant sheaves as modules over kernel algebras}

We are going to show that in the situation when a finite group scheme 
(resp., locally nicely ind-finite formal group scheme)
acts on a scheme $X$ there is a natural way to construct a kernel algebra on $X$ that captures
the corresponding category of equivariant sheaves. In fact, we can deal with a more general situation
of a groupoid over $X$.

\begin{thm}\label{equiv-sh-thm} 
(i) Let $(s,t):M\to X\times X$ be a groupoid over $X$. 
Assume that the maps $s$ and $t$ are finite and flat. Then the relative 
dualizing sheaf  $\om_s=s^!\OO_X$
has a natural structure of a convolution algebra on $M$. Its push-forward to $X\times X$ is
a pure kernel algebra $\AA_M=\AA_M(\om_s)$ over $X$. 
The category of $\AA_M$-modules is equivalent to the category of $M$-equivariant quasicoherent
sheaves on $X$. By definition, these are quasicoherent sheaves $F$ on $X$ equipped with an
isomorphism $\a:t^*F\to s^*F$ such that 
\begin{equation}\label{equiv-rel-eq}
m^*\a=p_1^*\a\circ p_2^*\a, \ \ e^*\a=\id_F,
\end{equation}
where
$m, p_1, p_2:M\times_X M\to M$ are the product and the projections, $e:X\to M$ is the unit map.

\noindent (ii) The above statements hold also in the case when $M$ is a formal groupoid over $X$ 
with the following changes.
We assume that 
locally over $S$ each connected component of $M$ is affine over $X$ with respect
to both $s$ and $t$ and can be presented as the inductive limit of a system
of closed subschemes $M_n\sub M_{n+1}\sub\ldots\sub M$ such that each $M_n$ is finite and flat
over $X$ with respect to both $s$ and $t$. We set $\om_s=s^!\OO_X$
(see Proposition \ref{formal-finite-duality-prop}(ii)). 
The definition of an $M$-equivariant structure on $F$ should be modified as follows: 
it is given by a map $\a\in\widehat{\Hom}(t^*F, s^*F)$,
and the relations \eqref{equiv-rel-eq} are imposed
using $\widehat{\Hom}$-spaces (see \eqref{hat-Hom-eq}).
\end{thm}

\Pf . Part (i) is a particular case of part (ii), so it is enough to consider the case of a formal groupoid.
Our assumptions imply that both morphisms $s,t:M\to X$
(and hence, $p_1,p_2:M\times_X M\to M$) are nicely ind-finite.
Since $M$ is a groupoid, the map $m:M\times_X M\to M$ differs from the projection $p_1$
by an automorphism of $m$. Therefore, $m$ is also nicely ind-finite.
Thus, we can use the functor $m^!$ defined in Proposition \ref{formal-finite-duality-prop}. 
In particular, we have the following isomorphism on $M\times_X M$:
\begin{equation}\label{m-om-eq}
m^!\om_s=m^!s^!\OO_S\simeq (s\circ m)^!\OO_S=(s\circ p_1)^!\OO_S\simeq p_1^!\om_s.
\end{equation} 
Now we define the convolution algebra structure on $\om_s$ as the map
corresponding by adjunction to the following map on $M\times_X M$:
\begin{equation}\label{conv-om-eq}
p_1^*\om_s\otimes p_2^*\om_s\to p_1^!\om_s \simeq m^!\om_s,
\end{equation}
where the arrow is given by the isomorphism \eqref{dual-XY-eq}.
The associativity of this product reduces to the commutativity of the diagram in
Corollary \ref{formal-duality-cor} for $M\times_X M\times_X M$ with $F=\om_s$,
along with some easier compatibilities.
The unit $e_*\OO_S\to \om_s=s^!\OO_S$ corresponds by adjunction to 
the isomorphism $s_!e_*\OO_S\simeq \OO_S$.

Using Lemma \ref{aff-monoid-lem} we see that $\AA_M$ is pure. 
Indeed, both $s$ and $t$ are {\it ldu}-affine, so we only need to check flatness
of $\om_s$ over $X$ with respect to both $s$ and $t$. But this follows immediately
from the construction of Proposition \ref{formal-finite-duality-prop}, since
$\om_s$ is the inductive limit of line bundles over $M_n$.

By Proposition \ref{formal-finite-duality-prop}(iv), an element
$\a\in\widehat{\Hom}(t^*F,s^*F)$ required in the definition of
an $M$-equivariant sheaf $F$ on $X$ is the same as a morphism
\begin{equation}\label{st-F-eq}
t^*F\otimes \om_s\to s^*F\otimes \om_s \wt{\ra} s^!F,
\end{equation}
where the last arrow is an isomorphism by Proposition \ref{formal-finite-duality-prop}(ii). 
By adjunction this is the same as a morphism
$$s_!(t^*F\otimes\om_s)\to F.$$
By the projection formula, the source of this map is isomorphic to $\AA_M\circ_X F$.
Conversely, given a morphism $\AA_M\circ_X F\to F$, reversing the above procedure we
get a map $t^*F\otimes\om_s\to s^*F\otimes\om_s$, or equivalently, an element
$\a\in\widehat{\Hom}(t^*F,s^*F)$.

We have to check that the module axiom for $F$ is equivalent to the
conditions \eqref{equiv-rel-eq}. Using the projection formula and the base change
formula we get an isomorphism
$$\AA_M\circ_X\AA_M\circ_X F\simeq (sp_1)_!((tp_2)^*F\ot p_1^*\om_s\ot p_2^*\om_s),$$
where $p_1,p_2:M\times_X M\to M$ are the projections.
Furthermore, the map 
$$\AA_M\circ_X\AA_M\circ_X F\to \AA_M\circ_X F\to F$$
obtained from two action maps $\AA_M\circ_X F\to F$ corresponds by adjunction to the map
\begin{equation}\label{tpsp-eq}
\begin{array}{l}
(tp_2)^*F\ot p_1^*\om_s\ot p_2^*\om_s\to p_2^*(s^*F\ot\om_s)\ot p_1^*\om_s\to 
p_1^*(t^*F\ot \om_s)\ot p_1^!\OO_M
\to p_1^!(t^*F\ot\om_s)\to \\
p_1^!s^!F\simeq (sp_1)^!F,
\end{array}
\end{equation}
where the first and the last arrows are induced by \eqref{st-F-eq}, the second arrow is
induced by the isomorphism $p_2^*\om_s\simeq p_1^!\OO_M$ (using Proposition 
\ref{formal-finite-duality-prop}(iii)). Note that we have
$$(sp_1)^!\OO_X\simeq p_1^*\om_s\ot p_1^!\OO_M\simeq p_1^*\om_s\ot p_2^*\om_s,$$
hence, $(sp_1)^!F\simeq (sp_1)^*F\ot p_1^*\om_s\ot p_2^*\om_s$.
Now the first relation in \eqref{equiv-rel-eq} is equivalent to the condition that under this isomorphism
the map \eqref{tpsp-eq} coincides with the map
\begin{equation}\label{tsm*-eq}
(tm)^*F\ot p_1^*\om_s\ot p_2^*\om_s\to (sm)^*F\ot p_1^*\om_s\ot p_2^*\om_s
\end{equation}
induced by $\a$. Indeed, this can be deduced easily from the commutativity of the diagram in Proposition \ref{formal-finite-duality-prop}(iii) for the composition of $s$ with $p_1$ 
and from the isomorphism
of Proposition \ref{formal-finite-duality-prop}(iv) for the composed map $sp_1$.
On the other hand, the module axiom requires that \eqref{tpsp-eq} coincides with the map
\begin{equation}\label{tsm-eq}
(tm)^*F\ot p_1^*\om_s\ot p_2^*\om_s
\to m^*t^*F\ot m^!\om_s\to m^*(t^*F\ot\om_s)\ot m^!\OO_M\to m^!(t^*F\ot\om_s)\to m^!s^!F\simeq (sm)^!F,
\end{equation}
where the first arrow is induced by \eqref{conv-om-eq}, the second by the isomorphism
$m^!\om_s\simeq m^*\om_s\ot m^!\OO_M$, the third by the isomorphism \eqref{dual-formal-eq},
and the last by \eqref{st-F-eq}. Using the isomorphism 
$(sm)^!F\simeq (sm)^*F\ot p_1^*\om_s\ot p_2^*\om_s$ and the 
commutativity of the diagram in Proposition
\ref{formal-finite-duality-prop}(iii) for the composition of $s$ with $m$ one can check that the maps
\eqref{tsm-eq} and \eqref{tsm*-eq} are the same, which finishes the proof.
\ed

\begin{ex} 
In the case when $S=\Spec(k)$, where $k$ is a field of characteristic zero, and $M$ is the infinitesimal 
groupoid corresponding to a Lie algebroid $L$ over $X$, the kernel algebra $\AA_M$ is isomorphic
to the universal enveloping algebra of $L$ (viewed as a D-algebra, see \cite{BB}). Thus, the above
Theorem reduces in this case to the interpretation of $L$-modules as $M$-equivariant sheaves on $X$.
\end{ex}

Assume that a finite group scheme $G$ acts on $X$. Then we have the corresponding action groupoid
$$\Ga_G=\{(gx,g,x)\ |\ g\in G, x\in X\}\sub X\times G\times X,$$
where the maps $s,t:\Ga_G\to X$ are induced by $p_1$ and $p_3$ (projections from 
$X\times G\times X$). Hence, by the above Theorem we
have the corresponding pure kernel algebra over $X$
$$\AA^G_X:=\AA_{\Ga_G}=p_{13*}(p_2^*\om_{G/S}).$$

Now let $G$ be a formal group scheme, locally nicely ind-finite over $S$.
We define a $G$-equivariant structure on a quasicoherent sheaf $F$ on $X$ as
an isomorphism in the sense of $\widehat{\Hom}$-spaces
of two pull-backs of $F$ to $G\times X$ 
(via the projection and the action map $G_n\times X\to X$), satisfying the usual compatibility
conditions (still in $\widehat{\Hom}$-spaces).
Then we can apply Theorem \ref{equiv-sh-thm}(ii) to the formal
groupoid $\Ga_G$ to define the pure kernel algebra $\AA^G_X$ by the same formula as above 
(but with $p_{13*}$ replaced by $p_{13!}$). Theorem \ref{equiv-sh-thm}
implies the following result.

\begin{cor}\label{action-ker-cor} 
Let $G$ be finite group scheme (resp., locally nicely ind-finite
formal group scheme) acting on a scheme $X$. 
Then the category of $G$-equivariant quasicoherent sheaves on $X$ is equivalent
to the category of $\AA^G_X$-modules. 
\end{cor}

Note that in the above situation the coherent sheaf
$\OO_{\Ga_G}$ on $X\times G\times X$ has a natural structure of
a kernel representation of $G$ over $X$. Hence, using the construction
of Proposition \ref{ker-G-rep-prop} we can associate with this kernel representation the kernel algebra
$\AA(\OO_{\Ga_G},\om_{G/S})$. It is easy to see that it is isomorphic to $\AA^G_X$
(see Lemma \ref{quasi-1-coc-alg-lem} for a more general statement).

\subsection{Twisting geometric actions by $1$-cocycles}
\label{ker-twist-sec}

Recall that if $\XX$ is an algebraic stack then
a cohomology class $e\in H^2(\XX,\G_m)$ (i.e., a $\G_m$-gerbe over $\XX$)
allows to define a category of twisted quasicoherent
sheaves on $\XX$ (see e.g., section 2 of \cite{Lie}). 
We will show that some of the twisted categories of sheaves 
for the global quotient stack $\XX=[X/G]$ are related to kernel representations of $G$ over $X$
of a special kind. It is convenient to describe $\G_m$-gerbes over a stack using certain data over
its groupoid presentation. The relevant notion of a $1$-cocycle on a groupoid generalizes the
well known description of gerbes using open coverings and line bundles on pairwise intersections 
(see section 1 of \cite{Bryl}, and also \cite{Hitchin}, \cite{Ben-Bassat}). 

\begin{defi}
Let $X$ be a scheme, $(s,t):M\to X\times X$ a groupoid (resp., formal groupoid) over $X$. 
A {\it $1$-cocycle on $M$} is a line bundle $\LL$ on $M$ equipped with an isomorphism
\begin{equation}\label{twisting-def}
p_1^*\LL\otimes p_2^*\LL\to m^*\LL
\end{equation}
on $M\times_X M$, where $m,p_1,p_2:M\times_X M\to M$ are the product and the projection maps,
such that the following diagram on $M\times_X M\times_X M$) is commutative:
\begin{diagram}\label{1-coc-comm-eq}
p_1^*\LL\otimes p_2^*\LL\otimes p_3^*\LL&\rTo{}&
p_1^*\LL\otimes p_{23}^*m^*\LL\\
\dTo{}&&\dTo{}\\
p_{12}^*m^*\LL\otimes p_3^*\LL&\rTo{}&(mm)^*\LL
\end{diagram}
where $mm:M\times_X M\times_X M\to M$ is the composed product map.

In the case when $M=\Ga_G=G\times X$ is the action groupoid associated with an action of a group
scheme (resp. formal group scheme) $G$ on $X$ we will also say that $\LL$ is a
{\it $1$-cocycle of $G$ with values in $\PPic(X)$}.
\end{defi}

Let $e:X\to M$ be the unit.
Pulling back \eqref{twisting-def} via $e\times e: X\to M\times_X M$ we obtain an isomorphism
\begin{equation}\label{unit-1-coc}
\OO_X\to e^*\LL
\end{equation}
on $X$. Using the compatibility diagram above one can check that this isomorphism is 
compatible with the restrictions of maps \eqref{twisting-def} to $M$ via the maps
$e\times\id_M:M\to M\times_X M$ and $\id_M\times e:M\to M\times_X M$.

\begin{exs}
1. In the case when 
$X\to Y$ is an fppf map of schemes (or a presentation of an algebraic stack $Y$) and $M=X\times_Y X$
is the corresponding groupoid,
the notion of a $1$-cocycle on $M$ is an algebraic version of {\it gerbe data} of \cite{Bryl}
(resp., a {\it presentation of a gerbe} of \cite{Ben-Bassat}). 
Let us show how to construct a $\G_m$-gerbe $\GG_{\LL}$ on $Y$ 
corresponding to a $1$-cocycle $\LL$ on $M=X\times_Y X$. For an open $U\to Y$ (in flat topology) let us consider $E=X\times_Y U$, and let $\LL_U$ denote the pull-back of $\LL$ to $E\times_U E$. 
Then $\GG_{\LL}(U)$ is the category of line bundles $\xi$ over $E$ equipped with an isomorphism
$$p_1^*\xi\simeq p_2^*\xi\otimes\LL_U$$
over $E\times_U E$, satisfying the obvious compatibility on $E\times_U E\times_U E$. 
Note that the pull-back of this gerbe to $X$ is trivialized, namely, the $1$-cocycle $\LL$ itself 
may be viewed as an object of $\GG_{\LL}(X)$. Conversely, a $\G_m$-gerbe on $Y$ equipped
with a trivialization of its pull-back to $X$ gives rise to a $1$-cocycle on $X\times_Y X$: one
simply looks at the difference of two trivializations of our gerbe to $X\times_Y X$ induced by
two projections to $X$. It is easy to see in this way one gets an equivalence
between the category of $1$-cocycles on $X\times_Y X$ and the category of $\G_m$-gerbes
on $Y$ with trivialized pull-back to $X$. Let us observe also that a $1$-cocycle $\LL$ on $X\times_Y X$
gives rise to an element of the relative Picard group $\Pic(X/Y)$ defined by Grothendieck
in \cite{FGA}, and in the situation considered in \cite{Gr-BIII} the corresponding element
of the cohomological Brauer group is the class of the gerbe $\GG_{\LL}$.

\noindent
2. In the case when $G$ is discrete a $1$-cocycle of $G$ with values in $\PPic(X)$ induces
a $1$-cocycle of $G$ with values in the abelian group $\Pic(X)$ in the usual sense, where $\Pic(X)$ is
viewed as a $G$-module. The obstacle to lifting a usual $1$-cocycle of $G$ with values in
$\Pic(X)$ to $\PPic(X)$ is a certain cohomology class in $H^3(G,H^0(X,\OO^*))$ that measures
the defect for commutativity of the diagram \eqref{1-coc-comm-eq}.
Note also that to give a $1$-cocycle of $G$ with values in $\Pic(X)$ is equivalent to lifting
the homomorphism $G\to\Aut(X)$ to a group homomorphism $G\to\Aut(X)\ltimes\Pic(X)$.
On the other hand, we will show below that a $1$-cocycle with values in $\PPic(X)$ gives rise
to a kernel representation of $G$ over $X$.
Thus, if we think of the notion of kernel representation as a refinement of a homomorphism
$G\to\Auteq(D_{qc}(X))$ then $1$-cocycles with values in $\PPic(X)$ correspond
to homomorphisms that factor through $\Aut(X)\ltimes\Pic(X)$.

\noindent
3. The holomorphic analog of the above notion can be used to study holomorphic
gerbes on quotients by actions of discrete groups. Taking
$X$ to be a complex vector space and $G$ a lattice acting on $X$ by translations, one
gets a construction of gerbes on complex tori. It is easy to see that
all gerbes on complex tori appear in this way. Furthermore, due to triviality of the relevant
line bundles, these gerbes are in fact described by $2$-cocycles of $G$ with values in
the group of invertible functions on $X$. An explicit cocycle representative for each
equivalence class of gerbes was constructed in \cite{Ben-Bassat-gerbes} (these representatives
are similar to Appell-Humbert $1$-cocycles describing holomorphic
line bundles on complex tori).
\end{exs}

Given the connection of $1$-cocycles with gerbes outlined above we should expect that they
define twisted versions of the categories of equivariant sheaves.

\begin{defi} Let $(s,t):M\to X\times X$ be a groupoid over $X$, $\LL$ a $1$-cocycle on $M$.
An {\it $\LL$-twisted $M$-equivariant quasicoherent sheaf on $X$} is a quasicoherent sheaf 
$F$ on $X$ equipped with an
isomorphism $\a:\LL\otimes t^*F\to s^*F$ such that equations 
\eqref{equiv-rel-eq} are satisfied (these equations make sense because of the isomorphism
$m^*\LL\simeq p_1^*\LL\ot p_2^*\LL$ and the trivialization of $e^*\LL$).
Morphisms between $\LL$-twisted $M$-equivariant sheaves are defined in a natural way.

The above definition can also be applied with appropriate changes to the case when
$M$ is a formal groupoid over $X$, satisfying the conditions of Theorem \ref{equiv-sh-thm}(ii):
the morphism $\a$ and the compatibilities should be formulated using $\widehat{\Hom}$-spaces
(see \eqref{hat-Hom-eq}).
\end{defi}

If $\LL$ and $\LL'$ are $1$-cocycles on a groupoid $(s,t):M\to X\times X$ then
there is a natural structure of a $1$-cocycle on the tensor product of line bundles
$\LL\otimes\LL'$.

\begin{defi}. With a line bundle $L$ on $X$ we associate the line bundle 
$\de L$ on $M$ by $\de L=s^*L\otimes t^*L^{-1}$.
Note that $\de L$ has a natural structure of a $1$-cocycle on $M$ that we call
a {\it coboundary $1$-cocycle}.
If $\LL$ and $\LL'$ are $1$-cocycles on $M$
equipped with an isomorphism $\LL'\simeq \LL\otimes\de L$ of $1$-cocycles
for some line bundle $L$ on $X$,
then we say that the $1$-cocycles $\LL$ and $\LL'$ are {\it cohomologous}.
\end{defi}

In the situation of Example 1 above
cohomologous $1$-cocycles lead to equivalent gerbes. In general, corresponding
categories of twisted $M$-equivalent sheaves on $X$ will be equivalent.
Indeed, if $F$ is an $\LL$-twisted $M$-equivariant sheaf on $X$ then
$F\otimes L$ has a natural structure of an $\LL\otimes\de L$-twisted $M$-equivariant sheaf .

We have the following analog of Theorem \ref{equiv-sh-thm} in the twisted case.

\begin{thm}\label{tw-equiv-sh-thm} 
Let $(s,t):M\to X$ be a formal groupoid over $X$ satisfying the conditions of
Theorem \ref{equiv-sh-thm}(ii), and let $\LL$ be a $1$-cocycle on $M$. Then
the sheaf $\LL\otimes\om_s$
has a natural structure of a convolution algebra on $M$. Its push-forward to $X\times X$ is
a pure kernel algebra $\AA_M(\LL\otimes\om_s)$ over $X$. 
The category of $\AA_M(\LL\otimes\om_s)$-modules is equivalent to the category of $\LL$-twisted
$M$-equivariant quasicoherent sheaves on $X$.
\end{thm}

The proof is similar to that of Theorem \ref{equiv-sh-thm}. Note that the convolution
algebra structure on $\LL\otimes\om_s$ is induced by the maps \eqref{conv-om-eq} and
\eqref{twisting-def} along with the natural map
$m^*\LL\otimes m^!\om_s\to m^!(\LL\times\om_s)$ (see Lemma \ref{tw-conv-lem} below
for a more general construction).

The following generalization of the notion of $1$-cocycle leaves only the features of this notion that are necessary to get a convolution algebra on $M$.

\begin{defi} Let $(s,t):M\to X\times X$ be a formal groupoid over a scheme $X$.
A {\it quasi-$1$-cocycle on $M$} is an object $\FF\in D_{qc}(M)$ equipped with a morphism
$\a:p_1^*M\ot p_2^*M\to m^*M$ 
on $M\times_X M$ as in the definition of $1$-cocycle (however, this morphism is not required to be an isomorphism) and subject
to the same constraint on $G\times G\times G\times X$.
In addition, $\FF$ should be equipped
with a map $\OO_X\to e^*\FF$ (which we call a {\it unit for }$\FF$), compatible with the 
pull-backs of $\a$ to $M$ via $e\times\id_M$ and $\id_M\times e$.

In the case when $M=\Ga_G=G\times X$ is an action groupoid associated with an action of a group
scheme (resp. formal group scheme) $G$ on $X$ we will also say that $\FF$ is a
a {\it quasi-$1$-cocycle of $G$ with values in $D_{qc}(X)$}.
\end{defi}

\begin{lem}\label{tw-conv-lem} 
Let $(s,t):M\to X$ be a formal groupoid over $X$ satisfying the conditions of
Theorem \ref{equiv-sh-thm}(ii), and let
$\FF\in D_{qc}(M)$ be a quasi-$1$-cocycle.
Then $\FF\otimes \om_s\in D_{qct}(M)$ 
has a natural structure of a convolution algebra over $M$.
\end{lem}

\Pf . As the convolution product on $\FF\otimes \om_s$ we
take the following map on $M$:
$$m_*(p_1^*\FF\otimes p_2^*\FF\otimes p_1^*\om_s\otimes p_2^*\om_s)\to
m_*(m^*\FF\otimes p_1^*\om_s\otimes p_2^*\om_s)\to
\FF\otimes m_*(p_1^*\om_s\otimes p_2^*\om_s)\to\FF\otimes\om_s,$$
where we use consecutively the quasi-$1$-cocycle structure on $\FF$, the projection formula,
and the convolution product on $\om_s$ constructed in Theorem \ref{equiv-sh-thm}. 
The associativity axiom follows from the associativity
of the convolution product on $\om_{s}$ and the compatibility on $M\times_X M\times_X M$ 
in the definition of quasi-$1$-cocycle. 
The unit on $\FF\otimes \om_{s}$ is given by the composition
$$e_*\OO_X\to \FF\otimes e_*\OO_X\to\FF\otimes \om_s,$$
where the maps are induced by the units for $\FF$ and for $\om_s$. 
\ed

We are mostly interested in the case when $M$ is the action groupoid $\Ga_G$ associated with an action of a group scheme (resp., formal group scheme) $G$ on $X$. In this case a $1$-cocycle
$\LL$ of $G$ with values in $\PPic(X)$ can be viewed
as a line bundle on $G\times X$, and
an $\LL$-twisted $G$-sheaf on $X$ is a quasicoherent sheaf $F$ on $X$ equipped with
an isomorphism 
\begin{equation}\label{equiv-coc-def}
F_{gx}\to \LL_{g,x}\otimes F_x
\end{equation}
over $G\times X$, where $F_{gx}$ denotes the pull-back of $F$ under the map
$G\times X\to X:(g,x)\mapsto gx$, etc. This morphism should reduce to the identity over
$e\times X$ (recall that $\LL_{e,x}$ is trivialized), and
the natural diagram on $G\times G\times X$ should be commutative.
We denote the category of $\LL$-twisted $G$-sheaves on $X$ by $\Qcoh^{\LL}_G(X)$.
In the case when $G$ is finite we also denote by
$\Coh^{\LL}_G(X)\sub\Qcoh^{\LL}_G(X)$ the full subcategory
consisting of coherent sheaves equipped with the above data.
Note that for a line bundle $L$ on $X$ the coboundary $\de L$ is given by 
$\de L_{g,x}=L_{gx}\otimes L_x^{-1}$.


In the case when $G$ is a group scheme the category $\Qcoh^{\LL}_G(X)$ should be viewed as a category twisted quasicoherent sheaves on the global quotient stack $[X/G]$.
The next lemma states the corresponding invariance of this category under certain changes of a presentation (and extends it to the case of formal group schemes).

\begin{lem}\label{gerb-lem} Assume that $G$ is finite group scheme (resp., locally nicely ind-finite
formal group scheme) acting on $X$, and let $(X',G')$ be another data of the same kind.
Assume that we have a surjective homomorphism of groups $G'\to G$ with the kernel $K$ which is
finite group scheme (flat over $S$), and a 
$G'$-equivariant morphism $X'\to X$. Assume also that $X'\to X$ is a $K$-torsor.
Let $\LL$ be a $1$-cocycle of $G$ with values in $\PPic(X)$, and let $\LL'$ be its pull-back to 
$G'\times X'$. Then $\Qcoh^{\LL}_G(X)$ is equivalent to $\Qcoh^{\LL'}_{G'}(X')$.
\end{lem}

The proof is a straightforward application of the flat descent for quasicoherent sheaves.

With a $1$-cocycle $\LL$ of $G$ with values in $\PPic(X)$
one can associate a kernel representation of $G$ over $X$ by setting
$$V_{\LL}=\ga_{G*}(\LL),$$
where $\ga_G:G\times X\to X\times G\times X: (g,x)\mapsto (gx,g,x)$.
It is easy to see that for 
the trivial $1$-cocycle we get precisely the  kernel representation $V^G_X$
associated with the action of $G$ on $X$.

More generally, if $\FF$ is a quasi-$1$-cocycle of $G$ with values in $D_{qc}(X)$ then
one we associate with $\FF$ a kernel quasi-representation of $G$ over $X$ by 
\begin{equation}\label{quasi-1-coc-quasi-rep}
V_{\FF}:=\ga_{G*}(\FF)\in D_{qc}(X\times G\times X).
\end{equation}
This construction also works in the case when $G$ is a locally nicely ind-finite formal group scheme
provided $\FF\in D_{G,\clim}(G\times X)$ (i.e.,
locally over $G$ the cohomology of $\FF$ are in $A_{\clim}$).
In this case we have the corresponding kernel algebra
$$\AA^G_X(\FF):=\AA(V_{\FF},\om_{G/S})=(gx,x)_*(p_1^*\om_{G/S}\otimes \FF)$$
(see Proposition \ref{ker-G-rep-prop}).
One can easily check that this is the same kernel algebra as the one obtained from
the convolution algebra structure on $\FF\ot p_G^*\om_{G/S}$ constructed in Lemma
\ref{tw-conv-lem}, where $p_G:\Ga_G=G\times X\to G$ is the projection. More precisely,
one has the following statement (the proof is straightforward).

\begin{lem}\label{quasi-1-coc-alg-lem} Let $G$ be a finite group scheme (resp., locally nicely ind-finite
formal group scheme)
acting on $X$, and let $\FF\in D_{qc}(G\times X)$ (resp., $\FF\in D_{G,\clim}(G\times X)$)
be a quasi-$1$-cocycle $\FF$ of $G$ with values in $D_{qc}(X)$.
Then the push-forward of the convolution algebra $\FF\otimes p_G^*\om_{G/S}$
with respect to the homomorphism of groupoids
$\ga_G:\Ga_G\to X\times G\times X$ can be identified with
the convolution algebra $V_{\FF}\ot p_2^*\om_{G/S}$ on $X\times G\times X$
associated with the quasi-representation $V_{\FF}$ (see Proposition \ref{ker-G-rep-prop}). 
Hence, we have an isomorphism of kernel algebras over $X$:
$$\AA^G_X(\FF)\simeq\AA_{\Ga_G}(\FF\otimes p_G^*\om_{G/S}).$$
\end{lem}

We have the following analog of Corollary \ref{action-ker-cor} for twisted $G$-equivariant sheaves
that can be deduced from Theorem \ref{tw-equiv-sh-thm}.

\begin{cor}\label{tw-action-ker-cor}
Let $G$ be finite group scheme (resp., locally nicely ind-finite
formal group scheme) acting on a scheme $X$, and let
$\LL$ be a $1$-cocycle of $G$ with values in $\PPic(X)$.
Then
the category $\Qcoh^{\LL}_G(X)$ is equivalent to $\AA^G_X(\LL)-\mod$. In the case when
$G$ is finite, we also have $\Coh^{\LL}_G(X)\simeq\AA^G_X(\LL)-\mod^c$.
\end{cor}

The dependencies between the notions introduced above can be roughly summarized by 
the following picture

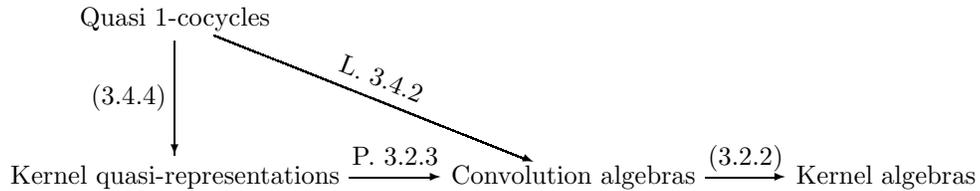
\begin{diagram}
\text{Quasi 1-cocycles} &&&&\\
\dTo{\eqref{quasi-1-coc-quasi-rep}} & \rdTo{\text{L. \ref{tw-conv-lem}}}&&&\\
\text{Kernel quasi-representations} &\rTo{\text{P. \ref{ker-G-rep-prop}}}&\text{Convolution algebras}&
\rTo{\eqref{conv-alg-ker-alg}}&\text{Kernel algebras}
\end{diagram}

Quasi-$1$-cocycles appear naturally via the following construction.

\begin{lem}\label{push-quasi-1-coc-lem} 
Let $G$ be a finite group scheme (resp., locally nicely ind-finite formal group scheme).
Let $f:X\to Y$ be an affine 
$G$-equivariant morphism between $S$-schemes equipped with a $G$-action, 
and let $\FF$ be a quasi-$1$-cocycle of $G$ with values in $D_{qc}(X)$ (in the formal case
we assume that $\FF\in D_{G,\clim}(G\times X)$). Then
$(\id_G\times f)_*\FF$ has a natural structure of a quasi-$1$-cocycle of $G$ with values in $D_{qc}(Y)$.
Furthermore, we have an isomorphism of kernel algebras
\begin{equation}\label{push-quasi-1-coc-isom}
(f\times f)_*\AA^G_X(\FF)\simeq\AA^G_Y((\id_G\times f)_*\FF).
\end{equation}
\end{lem}

\Pf . Set $\GG=(\id_G\times f)_*\FF$. 
Note that in the formal case we have 
$\GG\in D_{qc}(G\times Y)$ by 
Proposition \ref{formal-pull-push-prop}(vii), since $f$ is an affine morphism.
The unit $\OO_X\to(e\times\id_X)^*\FF$ for $\FF$ 
induces a unit map 
$$\OO_X\to (e\times\id_Y)^*\GG\simeq f_*(e\times\id_X)^*\FF,$$
where the latter isomorphism is given by the base change formula. Also,
we have a natural morphism on $G\times G\times Y$:
$$\GG_{g_1,g_2x}\ot \GG_{g_2,x}\to (\id_{G\times G}\times f)_*(\FF_{g_1,g_2x}\ot \FF_{g_2,x}).$$
Hence, the quasi-$1$-cocycle structure on $\FF$ induces a similar structure on $\GG$.

Recall that $\AA^G_X(\FF)$ comes from the convolution algebra
$\FF\otimes p_G^*\om_{G/S}$ over the action groupoid $\Ga_G(X)=G\times X$.
Similarly, $\AA^G_Y(\GG)$ comes from the convolution algebra
$\GG\otimes p_G^*\om_{G/S}$ over $\Ga_G(Y)=G\times Y$.
The projection formula and the base change formula give an isomorphism
$$\GG\otimes p_G^*\om_{G/S}\simeq (\id_G\times f)_*(\FF\otimes p_G^*\om_{G/S})$$
over $\Ga_G(Y)$ that induces the required isomorphism \eqref{push-quasi-1-coc-isom}.
The compatibility of the product maps and of units is easy to check.
\ed

\begin{prop}\label{ker-rep-twist-prop} 
Let $G$ be a finite group scheme (resp., locally nicely ind-finite formal group) 
acting on $X$, and let $\FF$ be a flat quasicoherent
sheaf on $G\times X$ equipped with a quasi-$1$-cocycle structure. 
Then the kernel algebra
$\AA^G_X(\FF)$ is pure. If $G$ is a finite group scheme and 
$\FF$ is a vector bundle then the pure kernel algebra $\AA^G_X(\FF)$ is finite.
\end{prop}

\Pf . By Lemma \ref{quasi-1-coc-alg-lem}, we can realize $\AA^G_X(\FF)$ as coming from
the convolution algebra $\FF\otimes p_G^*\om_{G/S}$ on $\Ga_G$. Now the first assertion
follows easily from Lemma \ref{aff-monoid-lem}. The second assertion is clear since
in this case $\AA^G_X(\FF)$ is the push-forward of $\FF\otimes p_G^*\om_{G/S}$ under a finite map.
\ed

\begin{rem} It is easy to see that if we change a $1$-cocycle $\LL$ to a cohomologous
cocycle $\LL\otimes\de L$ then the corresponding kernel representations of $G$ over $X$ 
are related via
the autoequivalence of $D_{qc}(X)$ given by tensoring with $L$ (as in Proposition
\ref{ker-equiv-prop}):
$$V_{\LL\otimes\de L}\simeq \De_*L\circ_X V_{\LL}\circ_X \De_*L^{-1}.$$
Hence, the corresponding kernel algebras are also related in a similar way.
\end{rem}

\begin{ex} Assume that the action of $G$ on $X$ is trivial. Then
a $1$-cocycle of $G$ with values in $\PPic(X)$ is the same as 
a monoidal functor (over $S$) of stacks of groupoids
$G\to\PPic(X)$ over $S$, i.e., in the terminology of \cite{De-sm}, 
a {\it central extension} of $G$ by $\PPic(X)$). By definition, such a central extension 
is given by a line bundle $\LL$ on $G\times X$ equipped with
an isomorphism
\begin{equation}\label{hom-pic-eq}
\LL_{g_1g_2,x}\simeq \LL_{g_1,x}\otimes\LL_{g_2,x}
\end{equation}
satisfying the natural cocycle condition.
In this case the morphism of groupoids $\Ga_G\to X\times X$ is a composition of the natural
projection $p_2:\Ga_G\simeq G\times X\to X$ with the diagonal embedding $\De:X\to X\times X$.
Hence,
$$\AA^G_X(\LL)\simeq \De_*\OO_X[G]^{\LL},$$
where
\begin{equation}\label{twist-group-eq}
\OO_X[G]^{\LL}:=p_{2*}(p_G^*\om_G\otimes\LL)
\end{equation}
is an algebra in $D_{qc}(X)$ associated with $\LL$. 
Here the product on the right-hand side is induced by the
convolution algebra structure on $p_G^*\om_G\otimes\LL$, where $G\times X$ is
viewed as a group scheme (resp., formal group scheme) over $X$.
Note that the algebra $\OO_X[G]^{\LL}$ is not always commutative even if $G$ is commutative. 
This corresponds to the fact that our monoidal functor $G\to\PPic(X)$ does not have to respect
the commutativity constraints (i.e., it is not necesserily a homofunctor between the Picard categories). For example, if $X=S$ (the base scheme) then a $1$-cocycle of $G$ with values in 
$\PPic(S)$ over $S$ is the same as a central extension of group schemes
$$1\to \G_m\to \wt{G}\to G\to 1.$$
The corresponding algebra is the twisted group algebra corresponding to this extension.
\end{ex}

Now assume that $H$ is a group scheme (or a formal group scheme)
acting on $X$ (and $G$ acts trivially on $X$).
We say that {\it an $H$-equivariant structure} on a homofunctor $\LL$
from $G$ to $\PPic(X)$ is given if $\LL$ is equipped with an $H$-equivariant structure with
respect to the natural action of $H$ on $G\times X$, and the isomorphism \eqref{hom-pic-eq} is
compatible with the $H$-action. In this situation we will also write that $\LL$ gives a homofunctor
$G\to\PPic^H(X)$. Recall that for a commutative group scheme $A$ we denote by $\Eext(A,\G_m)$
the Picard stack of (commutative) extensions of $A$ by $\G_m$. In the case when we have
a homomorphism $H\to A$ we denote by $\Eext^H(A,\G_m)$ the Picard stack of extensions of
$A$ by $\G_m$ equipped with a trivialization over $H$.

\begin{lem}\label{prod-cocycle-lem} 
(i) In the above situation the category of
monoidal functors $G\to\PPic^H(X)$ over $S$ is equivalent to the category of
$1$-cocycles of $H\times G$ with values in $\PPic(X)$, trivialized over $H$.

\noindent
(ii) Now assume that $X$ is a commutative group scheme and the action of $H$ on $X$
is given by translations (so we have a homomorphism $H\to X$). Then we can replace $\PPic(X)$ and
$\PPic^H(X)$ by the Picard stacks
$\Eext(X,\G_m)$ and $\Eext^H(X,\G_m)$, respectively. Assume also that $G$ is a commutative group
scheme or a formal group scheme. Then the following two categories are equivalent:

\noindent
(a) biextensions of $G\times X$ by $\G_m$ trivialized over $G\times H$;

\noindent
(b) homofunctors $G\to\Eext^H(X,\G_m)$.

\noindent
This category is also equivalent to a full subcategory in the category of 
$1$-cocycles of $H\times G$ with values in $\Eext(X,\G_m)$, trivialized over $H$.
\end{lem}

\Pf . (i) To a  monoidal functor $\LL:G\to\PPic^H(X)$ over $S$ we associate the line bundle
$p_{23}^*\LL$ on $H\times G\times X$. 
It is easy to see that it is equipped with the required
structures. The inverse functor is given by the restriction to $e_H\times G\times X$,
where $e_H\in H$ is the unit.

\noindent
(ii) The equivalence of (a) and (b) follows from Proposition \ref{biext-prop1}. To realize the
category in (b) in terms of $1$-cocycles of $H\times G$, trivialized over $H$, use (i).
\ed

\subsection{$1$-cocycles and Cartier duality}\label{Cartier-sec}

Henceforward, whenever formal $k$-groups are mentioned it is assumed that $S=\Spec(k)$,
where $k$ is a field of characteristic zero.

Let $X$ be an $S$-scheme, $G$ a finite commutative group scheme (flat) over $S$ 
(resp., formal $k$-group), 
and let $G^*$ denote the Cartier dual group scheme (resp., affine commutative
algebraic group over $k$).

Given a $G^*$-torsor $\phi:E\to X$ consider its pull-back to $G\times X$: 
$(\id_G\times\phi): G\times E\to G\times X$. Viewing the universal
bicharacter $b:G\times G^*\to\G_m$ as a homomorphism of group schemes $G^*\to\G_m$
over $G\times X$ we can associate with the above $G^*$-torsor 
a line bundle $\LL_E$ over $G\times X$.
It is easy to see that $\LL_E$ has a natural structure of a homofunctor $G\to\PPic(X)$.
Note that this construction is local over $X$. Hence, the obtained homofunctor is
{\it locally trivial}, i.e., there exists an open covering $\UU\to X$ such that
the induced homofunctor $G\to\PPic(\UU)$ is trivial.

\begin{prop}\label{torsor-dual-prop} 
(i) The above construction gives an equivalences between the
category of $G^*$-torsors over $X$ and the category of locally trivial
homofunctors $G\to\PPic(X)$. Furthermore, for a $G^*$-torsor $\phi:E\to X$ we have
an isomorphism of sheaves of $\OO_X$-algebras
$$\phi_*\OO_E\simeq \OO_X[G]^{\LL_E},$$
where the algebra structure on the right is defined by \eqref{twist-group-eq}.

\noindent
(ii) Assume now that $X$ is a formal scheme over $k$ and $G$ is a formal $k$-group.
As above, we have an equivalence between the category of $G$-torsors over $X$ and
the category of locally trivial homofunctors $G^*\to\PPic(X)$. For a $G$-torsor
$\psi:E\to X$ and a torsion quasicoherent sheaf $F$ on $X$ we have a natural isomorphism of 
sheaves on $X$:
$$\psi_!\psi^!F\simeq F\otimes p_{2*}\LL_E,$$
where $\LL_E$ is a line bundle on $G^*\times X$ associated with $E$.
\end{prop}

\Pf . (i) The quasi-inverse functor is constructed as follows. Given a locally trivial homofunctor
$\LL:G\to\PPic(X)$, consider the functor of trivializations of $\LL$ on the category of
$X$-schemes. By definition, it associates with an $X$-scheme $S$ the set of isomorphisms
of the pull-back of $\LL$ to $G\to\PPic(S)$ with the trivial homofunctor. It is easy to see
that this functor is represented by a $G^*$-torsor.

To construct an isomorphism of $\OO_X$-algebras we consider their pull-backs to $E$.
For the left-hand side we have a $G$-equivariant isomorphism of algebras 
$$\phi^*\phi_*\OO_E\simeq \pi_E^*(\pi_{G^*})_*\OO_{G^*},$$
where for every $S$ scheme $Y$ we denote by $\pi_Y:Y\to S$ the projection to $S$.
On the other hand, the canonical trivialization of $\phi^*\LL_E$ gives rise to an algebra isomorphism
$$\phi^*\OO_X[G]^{\LL_E}\simeq \pi_E^*(\pi_{G})_*\om_{G},$$
where $(\pi_{G})_*\om_{G}$ is equipped with the convolution product.
It remains to use the $G^*$-equivariant isomorphism 
\begin{equation}\label{formal-group-dual-isom}
(\pi_{G^*})_*\OO_{G^*}\simeq(\pi_{G})_!\om_{G}.
\end{equation}

\noindent
(ii) It is enough to construct such an isomorphism for the trivial $G$-torsor and to check its functoriality
with respect to automorphisms of the trivial $G$-torsor. In this case we have $E=G\times X$ and
$\LL_E$ is trivial; $\psi:G\times X\to X$ is the projection. 
Hence, we have to construct an isomorphism 
\begin{equation}\label{psi!-dual-isom}
\psi_!\psi^!F\simeq F\otimes \pi_X^*(\pi_{G^*})_*\OO_{G^*}, 
\end{equation}
compatible with automorphisms of the trivial $G$-torsor over $X$.
Note that the projection $\psi:G\times X\to X$ has a nicely ind-finite structure
(since $G$ is a union of finite $k$-schemes). Hence, 
by Proposition \ref{formal-finite-duality-prop}(iii),
we have $\psi^!F\simeq\psi^*F\otimes \psi^!\OO_X$.
By the projection formula, this leads to an isomorphism
$$\psi_!\psi^!F\simeq F\otimes \psi_!\psi^!\OO_X.$$
Using the compatibility of the formation of $f^!\OO_X$ with
base changes (see Proposition \ref{formal-finite-duality-prop}(iii)) and the base change
formula (Theorem \ref{formal-base-change}(ii),(iii)) we get
$$F\otimes\psi_!\psi^!\OO_X\simeq F\otimes\psi_!p_G^*\om_G\simeq F\otimes\pi_X^*(\pi_G)_!\om_G,$$
where $p_G:G\times X\to G$ is the projection. 
It remains to use the $G$-equivariant isomorphism 
\eqref{formal-group-dual-isom}. The compatibility of \eqref{psi!-dual-isom}
with automorphisms of both sides induced by automorphisms of the trivial $G$-torsor over $X$
follows easily from the compatibility of the isomorphism \eqref{formal-group-dual-isom}
with the natural $G$-action on both sides.
\ed

It follows that in the situation of Proposition \ref{torsor-dual-prop}(ii) we have an isomorphism of 
kernel algebras over $X$ 
\begin{equation}\label{torsor-alg-eq}
\De_*\phi_*\OO_E\simeq \AA^G_X(\LL_E).
\end{equation}

Now assume that $H$ is a group scheme (resp., formal $k$-group) acting on $X$.
Since the equivalence of the above proposition is compatible with the pull-back functors
(with respect to morphisms $X'\to X$), we deduce the following corollary.

\begin{cor} Let $\phi:E\to X$ be a $G^*$-torsor and let $\LL_E$ be the corresponding homofunctor
$G\to\PPic(X)$. A lifting of the $H$-action on $X$ to an action on $E$ (commuting with $G^*$)
is the same as an $H$-equivariant structure on $\LL_E$ (compatible with the homofunctor structure).
\end{cor}

Now let us consider a special case when $X$ is a commutative group scheme (resp., formal $k$-group).
Then instead of considering all $G^*$-torsors we can restrict our attention to extensions $E$ of $X$
by $G^*$ in the category of commutative groups. Here is the corresponding specialization
of Proposition \ref{torsor-dual-prop}.

\begin{lem}\label{ext-Car-lem} 
(i) The construction of Proposition \ref{torsor-dual-prop} induces an equivalence between
the category of extensions of $X$ by $G^*$ (resp., $G$) in the category of sheaves of commutative groups and
the category of biextensions of $G\times X$ (resp., $G^*\times X$) by $\G_m$.

\noindent
(ii) Let $H\to X$ be a homomorphism. Then the category of Lemma \ref{prod-cocycle-lem}(ii)
is equivalent to the category of extensions $E$ of $X$ by $G^*$ in the category of sheaves of commutative groups, equipped with a splitting $H\to E$ over $H$. 
\end{lem}

\Pf .
(i) This follows from Proposition \ref{biext-prop2} because
of the vanishings 
$$\underline{\Ext}^1(G,\G_m)=\underline{\Ext}^1(G^*,\G_m)=0$$ 
(see Lemma \ref{finite-ext-van-lem}; in the case 
when $G$ is a formal $k$-group we use Lemma \ref{formal-ext-van-lem}(i)).

\noindent
(ii) This follows immediately from part (i) using the definition (a) of Lemma \ref{prod-cocycle-lem}(ii).
\ed

Now let us go back to the situation when $X$ does not have a group structure.
Assume that $H$ is a finite commutative group scheme (resp., formal $k$-group) acting on 
$X$, $G$ is a finite commutative group scheme (resp., formal $k$-group), 
and $\phi:E\to X$ is a $G^*$-torsor equipped with an $H$-action that commutes with the $G^*$-action
and is compatible with the $H$-action on $X$. 
In addition, assume that we have an extension of commutative group schemes 
(resp., formal $k$-groups)
$$0\to G\to\wt{H}\to H\to 0.$$
Let $\BB$ be the corresponding biextension of $H\times G^*$ by $\G_m$ 
(see Lemma \ref{ext-Car-lem}(i)). 
We would like to identify the kernel algebra 
$(\phi\times\phi)_*\AA^H_E(\MM)$, where $\MM$ is a $1$-cocycle of $H$ with values in
$\PPic(E)$, with some kernel algebra corresponding to the action of $\wt{H}$ on $X$
(and a $1$-cocycle). For this we need to impose the following compatibility of $\MM$ with
the biextension $\BB$.

\begin{defi} Let $\MM$ be a $1$-cocycle of $H$ with values in $\PPic(E)$. We say that
$\MM$ is equipped with $(G^*,\BB)$-equivariant structure if an isomorphism
$$\MM_{h,g^*e}\simeq\BB_{h,g^*}\otimes\MM_{h,e}$$
is given, and the following diagrams on $H\times H\times G^*\times E$ and 
$H\times G^*\times G^*\times E$ are commutative:
\begin{diagram}
\MM_{h_1,h_2g^*e}\otimes\MM_{h_2,g^*e}&\rTo{}&\BB_{h_1,g^*}\otimes\BB_{h_2,g^*}\otimes
\MM_{h_1,h_2e}\otimes\MM_{h_2,e}\\
\dTo{}&&\dTo{}\\
\MM_{h_1h_2,g^*e}&\rTo{}&\BB_{h_1h_2,g^*}\otimes\MM_{h_1h_2,e}
\end{diagram}
\begin{diagram}
\MM_{h,g_1^*g_2^*e}&\rTo{}&\BB_{h,g_1^*}\otimes\MM_{h,g_2^*e}\\
\dTo{}&&\dTo{}\\
\BB_{h,g_1^*g_2^*}\otimes\MM_{h,e}&\rTo{}&\BB_{h,g_1^*}\otimes\BB_{h,g_2^*}\otimes\MM_{h,e}
\end{diagram}
\end{defi}

Since the pull-back of $\BB$ to $\wt{H}\times G^*$ is trivial, the pull-back of a
$(G^*,\BB)$-equivariant $1$-cocycle $\MM$ of $H$ with values in $\PPic(E)$ 
to $\wt{H}\times E$ is equipped
with a $G^*$-equivariant structure. Hence, it descends to a $1$-cocycle $\wt{\MM}$ of $\wt{H}$ with
values in $\PPic(X)$. It is easy to see that the restriction of $\wt{\MM}$ to $G$
is naturally isomorphic to $\LL_E:G\to\PPic(X)$, the homofunctor associated with $E$.

\begin{prop}\label{equiv-ker-prop} 
The above construction $\MM\mapsto\wt{\MM}$
gives an equivalences between the categories of
$(G^*,\BB)$-equivariant $1$-cocycles $\MM$ of $H$ with values in $\PPic(E)$ and
$1$-cocycles $\wt{\MM}$ of $\wt{H}$ with values in $\PPic(X)$ extending $\LL_E$.
Furthermore, we have an isomorphism of convolution algebras over $H\times X$ (viewed
as a groupoid over $X$):
\begin{equation}\label{quasi-coc-isom}
p_1^*\om_H\otimes(\id_H\times\phi)_*\MM\simeq (\pi\times\id_X)_!(p_1^*\om_{\wt{H}}\otimes\wt{\MM}),
\end{equation}
where $\pi:\wt{H}\to H$ is the projection. Here the left-hand side is the convolution algebra
associated with the quasi-$1$-cocycle structure on $(\id_H\times\phi)_*\MM$ (see Lemma 
\ref{push-quasi-1-coc-lem}), and
the right-hand side is the push-forward of the convolution algebra associated with $\wt{\MM}$
with respect to the map $\pi\times\id:\wt{H}\times X\to H\times X$ of groupoids over $X$
(see Lemma \ref{push-convolution-lem}).
Hence, we have an isomorphism of kernel algebras over $X$
\begin{equation}\label{main-ker-alg-isom}
(\phi\times\phi)_*\AA^H_E(\MM)\simeq \AA^{\wt{H}}_X(\wt{\MM}).
\end{equation}
\end{prop}

\Pf . Let us show how to go back from a $1$-cocycle $\wt{\MM}$ of $\wt{H}$ with values in
$\PPic(X)$ extending $\LL_E$ to a $(G^*,\BB)$-equivariant $1$-cocycle of $H$ with values in 
$\PPic(E)$. Let $\MM'$ denote the pull-back of $\wt{\MM}$ under the map 
$\id\times\phi:\wt{H}\times E\to\wt{H}\times X$. Then $\MM'$ is a $1$-cocycle of $\wt{H}$
with values in $\PPic(X)$. Note that the pull-back of $\LL_E$ under the map
$\id_G\times\phi:G\times E\to G\times X$ is naturally trivialized, and the $G^*$-equivariant structure
on it (where $G^*$ acts on $E$) is induced by the universal bicharacter $G\times G^*\to\G_m$.
It follows that $\MM'$ restricts to the trivial $1$-cocycle of $G$, so it descends to a $1$-cocycle 
$\MM$ of $\wt{H}/G\simeq H$ with values in $\PPic(E)$. The $(G^*,\BB)$-equivariant structure
on $\MM$ comes from the $G^*$-equivariance of $\MM'$ (recall that the pull-back of $\BB$
to $\wt{H}\times G^*$ is naturally trivialized but $G$ acts on it via the universal bicharacter
$G\times G^*\to\G_m$).

To construct an isomorphism \eqref{quasi-coc-isom} let us first consider the pull-backs of both sides
under the map $\id_H\times\phi:H\times E\to H\times X$.
By the definition, we have an isomorphism
$$(\pi\times\id_E)^*\MM\simeq (\id_{\wt{H}}\times\phi)^*\wt{\MM}$$
over $\wt{H}\times E$, compatible with $1$-cocycle structures (for $\wt{H}$).
Applying the flat base change formula for the cartesian diagram 
\begin{diagram}
\wt{H}\times E &\rTo{\id_{\wt{H}}\times\phi}& \wt{H}\times X\\
\dTo{\pi\times\id_E}&&\dTo{\pi\times\id_X}\\
H\times E&\rTo{\id_H\times \phi}& H\times X
\end{diagram}
and the projection formula for the morphism $\pi\times\id_E$
we derive an isomorphism 
$$(\id_H\times\phi)^*(\pi\times\id_X)_!(p_1^*\om_{\wt{H}}\otimes\wt{\MM})\simeq
(\pi\times\id_E)_!(p_1^*\om_{\wt{H}})\otimes\MM$$
of convolution algebras over $H\times E$ (viewed as a groupoid over $E$). 
Another use of the base change formula gives an isomorphism of
this with $p_1^*(\pi_!\om_{\wt{H}})\otimes \MM$.
Note that we can view $\pi:\wt{H}\to H$ as a $G$-torsor over $H$ 
(see Lemma \ref{formal-ext-van-lem}(ii)). Hence,
using the isomorphism of Proposition \ref{torsor-dual-prop}(ii)
we get
$$\pi_!\om_{\wt{H}}\simeq\pi_!\pi^!\om_H\simeq \om_H\otimes p_*\BB.$$
Thus, we derive an isomorphism
\begin{equation}\label{phiMHB-eq}
(\id_H\times\phi)^*(\pi\times\id_X)_!(p_1^*\om_{\wt{H}}\otimes\wt{\MM})\simeq
p_1^*(\om_H\otimes p_*\BB)\otimes \MM
\end{equation}
On the other hand, using the $(G_*,\BB)$-equivariance of $\MM$ we get 
$$(\id_H\times\phi)^*(\id_H\times\phi)_*\MM\simeq p_{13*}(p_{12}^*\BB)\otimes \MM\simeq
p_1^*(p_*\BB)\otimes \MM,$$
where $p_{ij}$ are projections from $H\times G^*\times E$, $p_1:H\times E\to H$ and
$p:H\times G^*\to H$ are also projections. Hence,
$$(\id_H\times\phi)^*(p_1^*\om_H\otimes(\id_H\times\phi)_*\MM)\simeq 
p_1^*(\om_H\otimes p_*\BB)\otimes \MM.$$
Comparing this with \eqref{phiMHB-eq} we get an isomorphism of the pull-backs of both sides
of \eqref{quasi-coc-isom} under $\id_H\times\phi$, compatible with convolution algebra structures.
It is easy to check that it is also compatible with $G^*$-actions, so it descends to 
an isomorphism over $H\times X$.

The isomorphism \eqref{main-ker-alg-isom} follows from \eqref{quasi-coc-isom} by passing
to associated kernel algebras, taking into account the isomorphism
$$(\phi\times\phi)_*\AA^H_E(\MM)\simeq\AA^H_X((\id\times\phi)_*\MM)$$
that follows from Lemma \ref{push-quasi-1-coc-lem}.
\ed

By Lemma \ref{prod-cocycle-lem}(i), an $H$-equivariant structure on a homofunctor 
$\LL:G\to\PPic(X)$ gives rise to an extension of $\LL$ to a $1$-cocycle $\wt{\LL}$ of
$H\times G$ with values in $\PPic(X)$, trivial over $H$.
Hence, we have
the corresponding pure kernel algebra 
$\AA^{H\times G}_X(\wt{\LL})$ over $X$.

On the other hand, a lifting of the $H$-action on $X$ to an action on $E$ gives
rise to a pure kernel algebra $\AA^H_E$ over $E$. 
Hence, we get a pure kernel algebra 
$(\phi\times\phi)_*\AA^H_E$ over $X$. 

\begin{cor}\label{equiv-ker-cor} 
Assume that $G$ and $H$ are finite commutative group schemes (resp., formal $k$-groups),
$X$ is a scheme with $H$-action, 
$\phi: E\to X$ is a $G^*$-torsor equipped with an $H$-action,
$\LL_E: G\to\PPic(X)$ is the corresponding $H$-equivariant homofunctor, and
$\wt{\LL}_E$ is the corresponding $1$-cocycle of $H\times G$ with values in $\PPic(X)$. 
Then there is an isomorphism of derived kernel algebras over $X$:
$$(\phi\times\phi)_*\AA^H_E\simeq \AA^{H\times G}_X(\wt{\LL}_E).$$
\end{cor}

\subsection{Compatibility with open coverings}

We are going to give a criterion for a kernel algebra coming from a convolution algebra
to be of affine type (see section \ref{Cech-res-sec}).

\begin{defi} Let $X$ be a scheme, and let $M$ be a formal groupoid over $X$. We say that
$M$ is compatible with an open covering $\UU=\sqcup_i U_i\to X$ (in flat topology) if
an isomorphism $\a:M_{\UU\times X}\wt{\to} M_{X\times\UU}$ over $M$ is given, where 
for a scheme $Y$ over $X\times X$ we denote by $M_Y$ the fibered product 
$M\times_{X\times X} Y$.
We require $\a$ to satisfy the following compatibilities.
First, the diagram
\begin{diagram}
M_{\UU\times X}\times_X M &\rTo{}& M\times_X M_{X\times\UU}\\
\dTo{} &&\dTo{}\\
M_{\UU\times X} &\rTo{\a}& M_{X\times\UU}
\end{diagram}
should be commutative, where the vertical arrows are obtained from the product map 
$m:M\times_X M\to M$ by the base change, and the upper arrow is given by the composition 
$$M_{\UU\times X}\times_X M\simeq (M\times_X M)_{\UU\times X\times X}\stackrel{\a_1}{\to}
(M\times_X M)_{X\times \UU\times X}\stackrel{\a_2}{\to} (M\times_X M)_{X\times X\times\UU}\simeq
M\times_X M_{X\times\UU}$$
with $\a_1$ and $\a_2$ induced by $\a$.
Second, the composition of $\a$ with the map $u_1:\UU\to M_{\UU\times X}$ should
coincide with the map $u_2:\UU\to M_{X\times\UU}$, where $u_1$ and $u_2$
obtained by the base change from the unit map $u:X\to M$.
\end{defi}

\begin{prop}\label{aff-groupoid-prop} 
Let $M$ be a groupoid over $X$ (resp., formal groupoid over $X$ 
such that the structure maps $s,t:M\to X$ are flat
and of {\it ldu}-pseudo-finite type). Assume that
$M$ is compatible
with an open covering $j:\UU=\sqcup_i U_i\to X$, where $j_i:U_i\to X$ are affine morphisms, 
Then for every convolution algebra $\VV$ on $M$
the kernel algebra $\AA_M(\VV)$ over $X$ is compatible with this open covering.
\end{prop}

\Pf . By definition, $\AA_M(\VV)$ is the push-forward of $\VV$ with respect to the structure
morphism $\pi:M\to X\times X$. Let us denote by $j_1:M_{\UU\times X}\to M$ and
$j_2:M_{X\times\UU}\to M$ the natural morphisms. Then
$$\AA_M(\VV)\circ_X\De(\UU)\simeq (\id\times j)_*(\id\times j)^*\AA_M(\VV)$$
can be identified with the push-forward to $X\times X$ of $j_{2*}j_2^*\VV\in D_{qc}(M)$.
Similarly, $\De(\UU)\circ_X\AA_M(\VV)$ can be identified with the push-forward of
$j_{1*}j_1^*\VV\in D_{qc}(M)$. Since we have an isomorphism of $M_{X\times\UU}$
and $M_{\UU\times X}$ over $M$, this leads to the required isomorphism
$$\AA_M(\VV)\circ_X\De(\UU)\simeq\De(\UU)\circ_X\AA_M(\VV).$$ 
Compatibility with the product (resp., unit) on $\AA_M(\VV)$ follows from the commutativity of the
diagram (resp., compatibility of $\a$ with $u_1$ and $u_2$)
in the above definition. 
\ed

For a formal group scheme $G$, {\it ldu-pf} ({\it ldu}-pseudo-finite) over a field $k$, we denote by 
$G_0$ its connected component of $1$, so that $G_0$ is infinitesimal, and
$G/G_0$ is an \'etale formal group corresponding to the $\Gal(\ov{k}/k)$-module $G(\ov{k})$.

\begin{prop}\label{aff-ker-prop} 
Let $G$ be a finite group scheme acting on $X$. 
Then for every
flat quasicoherent sheaf $\FF$ on $G\times X$ equipped with a quasi-$1$-cocycle structure,
the kernel algebra $\AA^G_X(\FF)$ is pure and of affine type. 
Similar result holds if $S=\Spec(k)$, where $k$ is a field, and
$G$ is a formal group scheme, {\it ldu-pf} over $k$, such that the action of $\Gal(\ov{k}/k)$ on $G(\ov{k})$
factors through $\Gal(k'/k)$ for some finite field extension $k\sub k'$.
\end{prop}

\Pf . Recall that $\AA^G_X(\FF)$ is the kernel algebra associated with a convolution algebra structure
on the sheaf $\FF\otimes p_G^*\om_{G/S}$ over the action groupoid $\Ga_G$ (see Lemma
\ref{quasi-1-coc-alg-lem}). It is pure by Proposition \ref{ker-rep-twist-prop}.
To check that it is of affine type, by Proposition \ref{aff-groupoid-prop}, 
it suffices to find an open affine covering $\UU\to X$, such that the action groupoid $\Ga_G$ is
compatible with $\UU$. 

Consider first the case when $G$ is a finite group scheme over $S$ (recall that we assume
it to be flat over $S$).
We can assume that $S$ is affine and pick a finite open affine covering $\wt{\UU}\to X$.
Then setting $\UU=G\times\wt{\UU}$ we obtain a $G$-equivariant
open affine covering of $X$ in flat topology, where the map $\UU\to X$
is the composition
$$\UU=G\times\wt{\UU}\to G\times X\to X$$
(the last arrow is given by the $G$-action on $X$). 

Now let us consider the case when $S=\Spec(k)$ and $G$ is a formal group scheme, {\it ldu-pf} over $k$.
We start by choosing an arbitrary Zariski open affine covering 
$\wt{\UU}=(\wt{U}_i)$ of $X$. Since the action of $G_0$ preserves each open subset $\wt{U}_i\sub X$,
for every $g\in G(\ov{k})$ there is a well-defined affine open subset $g(\wt{U}_i)\sub X_{\ov{k}}$, where
$X_{\ov{k}}$ is obtained from $X$ by extending scalars to $\ov{k}$. 
If $O\sub G(\ov{k})$ is a $\Gal(\ov{k}/k)$-orbit (necessarily finite by our assumption) then
the affine scheme $\sqcup_{g\in O}g(\wt{U}_i)$ comes from an \'etale open $U_{O,i}\to X$.
These form an \'etale covering of $X$. Furthermore, it is easy to see that there is an action of $G$ on
$(U_{O,i})$ compatible with its action on $X$, as required.
\ed

\subsection{Projective kernel representations and the Fourier-Mukai transform}
\label{projFM-sec}

The formalism of the Fourier-Mukai transform for abelian schemes
requires to introduce a projective version
of the notion of kernel representation (similar to the notion of projective representation
of a group).

\begin{defi} (i) Let $G$ be a group scheme (resp., formal group scheme) over $S$.
{\it A $2$-cocycle of $G$ with values in $\PPic(S)$} is a line bundle
$\LL$ over $G\times G$ equipped with an isomorphism
$$\LL_{g_1,g_2}\otimes\LL_{g_1g_2,g_3}\wt{\to}\LL_{g_1,g_2g_3}\otimes\LL_{g_2,g_3}$$
on $G^3$ such that the following diagram on $G^4$ is commutative:
\begin{diagram}
\LL_{g_1,g_2}\otimes\LL_{g_1g_2,g_3}\otimes\LL_{g_1g_2g_3,g_4}
&\rTo{}&&&\LL_{g_1,g_2g_3}\otimes\LL_{g_2,g_3}\otimes\LL_{g_1g_2g_3,g_4}\\
\dTo{}&&&&\dTo{}\\
\LL_{g_1,g_2}\otimes\LL_{g_1g_2,g_3g_4}\otimes\LL_{g_3,g_4}
&&&&\LL_{g_1,g_2g_3g_4}\otimes\LL_{g_2,g_3}\otimes\LL_{g_2g_3,g_4}\\
&\rdTo{}&&\ldTo{}&\\
&&\LL_{g_1,g_2g_3g_4}\otimes\LL_{g_2,g_3g_4}\otimes\LL_{g_3,g_4}&&
\end{diagram}
In addition, we assume that the bundle $(e\times e)^*\LL=\LL_{e,e}$ on $S$ is trivialized,
where $e:S\to G$ is a neutral element.
This induces a trivialization of the bundles $\LL_{e,g}$ and $\LL_{g,e}$ on $G$ (using
the cocycle isomorphism).

\noindent
(ii) Given a $2$-cocycle $\LL$ of a group scheme $G$ with values in $\PPic(S)$, 
and a scheme $X$ over $S$,
a {\it projective kernel representation of $G$ with the cocycle $\LL$} over $X$ is 
an object $V\in D_{qc}(X\times G\times X)$ equipped with an isomorphism
$$\mu:V\circ_X V\wt{\to} (\id_X\times m\times\id_X)^*V\ot p_{23}^*\LL$$
over $X\times G\times G\times X$ (where $m:G\times G\to G$ is the product), and an isomorphism
$$u:\De_*\OO_X\wt{\to}(\id_X\times e\times\id_X)^*V$$
subject to the following conditions.
The morphism $\mu$ induces the following morphism on $X\times G^3\times X$:
$$(V\circ_X V)\circ_X V\to \left((\id_X\times m\times\id_X)^*V\otimes p_{23}^*\LL\right)\circ_X V
\simeq \LL_{g_1,g_2}\otimes (\id_X\times m\times\id_G\times\id_X)^*(V\circ_X V).$$
Using $\mu$ again we get a map (an isomorphism)
$$(V\circ_X V)\circ_X V\to\LL_{g_1,g_2}\otimes \LL_{g_1g_2,g_3}\otimes 
(\id_X\times mm\times \id_X)^*V,$$
where $mm:G^3\to G$ sends $(g_1,g_2,g_3)$ to $g_1g_2g_3$.
Similarly, we construct a map
$$V\circ_X (V\circ_X V)\to\LL_{g_1,g_2g_3}\otimes \LL_{g_2,g_3}\otimes 
(\id_X\times mm\times \id_X)^*V,$$
and we require these maps to be the same.
The maps $u$ and $\mu$ should be compatible as follows: the composition
$$V\simeq \De_*\OO_X\circ_X V\stackrel{u\circ_X\id}{\to} (\id_X\times e\times\id_X)^*V\circ_X V\to V,$$
where the second arrow is induced by $(\id_X\times e\times\id_G\times\id_X)^*\mu$ and by
the trivialization of $(e\times\id_X)^*\LL$
(resp., similar composition starting with $\id\circ_X u$) should be equal to the identity map.

In the case when $G$ is a formal group scheme the above definition still make sense provided
$V\in D_{G,\clim}(X\times G\times X)$ and the support of $V$ is proper over $X\times G$ and over
$G\times X$.

\noindent
(iii) Given a $2$-cocycle $\LL$ of $G$ as above, and an $S$-scheme $X$ equipped with
an action of $G$, a
{\it coboundary for $\LL$ with values in $\PPic(X)$}
is a line bundle $\MM$ over $G\times X$ together with an isomorphism
$$\a(g_1,g_2;x):\MM_{g_1,g_2x}\ot\MM_{g_2,x}\wt{\to}\MM_{g_1g_2,x}\ot\LL_{g_1,g_2}$$
on $G\times G\times X$ and a trivialization of $(e\times\id_X)^*\MM=\MM_{e,x}$ compatible
with $\a(e,e;x)$,
such that the following diagram on $G^3\times X$
is commutative:
\begin{diagram}
\MM_{g_1,g_2g_3x}\MM_{g_2,g_3x}\MM_{g_3,x}&\rTo{\a(g_2,g_3;x)}&
\MM_{g_1,g_2g_3x}\MM_{g_2g_3,x}\LL_{g_2,g_3}&&\\
\dTo{\a(g_1,g_2;g_3x)}&&&\rdTo{\a(g_1,g_2g_3;x)}&\\
\MM_{g_1g_2,g_3x}\MM_{g_3,x}\LL_{g_1,g_2}&\rTo{\a(g_1g_2,g_3;x)}&
\MM_{g_1g_2g_3,x}\LL_{g_1g_2,g_3}\LL_{g_1,g_2}&\rTo{\b}&
\MM_{g_1g_2g_3,x}\LL_{g_1,g_2g_3}\LL_{g_2,g_3}
\end{diagram}
where $\b$ is induced by the structure of the $2$-cocycle on $\LL$ (we skipped the signs of
tensor product for brevity).
\end{defi}

In the case when the $2$-cocycle $\LL$ is trivial 
the definition (ii) (resp., (iii)) above reduces to the notion of a 
kernel representation (resp., a $1$-cocycle with values in $\PPic(X)$) considered in 
section \ref{ker-twist-sec}.

\begin{exs} 
1. Let $A$ be a commutative group scheme, $\LL$ a biextension of $A\times A$ by $\G_m$.
Then $\LL$ has a natural structure of a $2$-cocycle of $A$ with values in $\PPic(S)$.

\noindent
2. If $f:G\to G'$ is a homomorphism of group schemes over $S$ and $\LL'$ is a $2$-cocycle
of $G'$ with values in $\PPic(S)$ then $(f\times f)^*\LL'$ has a natural structure of a $2$-cocycle of $G$.
In particular, by the previous example if we have a homomorphism $f$ from $G$ to a commutative group scheme $A$ and a biextension $\LL$ of $A\times A$ by $\G_m$ then $(f\times f)^*\LL$ has the 
structure of a $2$-cocycle of $G$. More generally, if $A$ and $B$ are commutative group schemes
and $\LL$ is a biextension of $A\times B$ by $\G_m$ then for a pair of homomorphisms
$f:G\to A$, $g:G\to B$ the pull-back $(f\times g)^*\LL$ has a natural structure of a $2$-cocycle of $G$.
\end{exs} 

Note that $2$-cocycles of $G$ with values in $\PPic(S)$ form a commutative Picard category with
respect to the tensor product operation. 
For any line bundle $\MM\in\PPic(G)$ we can form a $2$-cocycle of $G$ with values in $\PPic(S)$
$$\La(\MM)_{g_1,g_2}:=\MM^{-1}_{g_1}\ot\MM^{-1}_{g_2}\ot\MM_{g_1g_2}.$$
A coboundary with values in $\PPic(S)$ for a $2$-cocycle $\LL$ is given by a line bundle $\MM$ 
together with an isomorphism $\LL\simeq\La(\MM)^{-1}$ of $2$-cocycles. 

For a $2$-cocycle $\LL$ of $G$ with values in $\PPic(S)$ and a line bundle $\MM$ on $G$
we have an equivalence from the category of projective kernel representations of $G$ over $X$
with the $2$-cocycle $\LL\ot\La(\MM)$ to the similar category for $\LL$ 
sending $V$ to $V\ot p_G^*\MM$,
where $p_G:X\times G\times X\to G$ is the projection.

Now assume that we have an action of $G$ on $X$.
Generalizing the construction of section \ref{ker-twist-sec}, given a $2$-cocycle $\LL$ of $G$ with
values in $\PPic(S)$ and a coboundary $\MM$ for $\LL$ with values in $\PPic(X)$,
we can define a projective kernel representation of $G$ over $X$ by setting
$$V_{\MM}=\ga_{G*}(\MM),$$
where $\ga_G(g,x)=(gx,g,x)$. For $\NN\in\PPic(G)$ the line bundle
$\MM\ot p_1^*\NN\in\PPic(G\times X)$ acquires a natural structure of
a coboundary for $\LL\ot\La(\NN)^{-1}$ with values in $\PPic(X)$.
The corresponding projective kernel representation is
$$V_{\MM\ot p_1^*\NN}\simeq V_{\MM}\ot p_G^*\NN.$$

\begin{ex} Let $A$ be an abelian variety over a field $k$, $\hat{A}$ the dual abelian variety.
Set $X_A=A\times \hat{A}$.
There is a natural homomorphism from the group $X_A(k)$
to the group of autoequivalences of $D^b(A)$ (viewed up to an isomorphism)
such that $A(k)$ acts by translations and 
$\hat{A}(k)$ acts by tensoring with line bundles. Let us give a 
categorified version of this picture that also works in the relative setting 
(cf. \cite{P-as}, sec.~3). Assume now that $A$ is an abelian scheme over $S$,
$\hat{A}$ is the dual abelian scheme, and $X_A=A\times_S\hat{A}$.
Then we are going to construct a projective kernel representation of $X_A$ over
$A$ with the $2$-cocycle $\LL_A=p_{32}^*\PP$, where $\PP$ is the 
Poincar\'e line bundle on $X_A$ and $p_{32}:X_A^2=(A\times\hat{A})^2\to A\times\hat{A}$
sends $(x_1,\xi_1,x_2,\xi_2)$ to $(x_2,\xi_1)$. 
Namely, if we let $X_A$ act on $A$ by translations (with $\hat{A}$
acting trivially), then we have a natural coboundary $\MM$ for $\LL_A$ with values in 
$\PPic(A)$ given by 
\begin{equation}\label{cob-L-A}
\MM_{(x,\xi;y)}=\PP_{y,\xi}, 
\end{equation}
where $(x,\xi;y)\in (A\times\hat{A})\times A$.
Indeed, the biextension structure on $\PP$ gives rise to an isomorphism
$$\MM_{x_1,\xi_1;x_2+y}\ot\MM_{x_2,\xi_2;y}=\PP_{x_2+y,\xi_1}\ot\PP_{y,\xi_2}\simeq
\PP_{y,\xi_1+\xi_2}\ot \PP_{x_2,\xi_1}=\MM_{x_1+x_2,\xi_1+\xi_2;y}\ot\PP_{x_2,\xi_1}.$$
Therefore, we have the corresponding projective kernel representation of $X_A$ over $A$ given by
\begin{equation}\label{V-A}
V(A):=V_{\MM}=(x+y,x,\xi,y)_*\MM\in\Coh(A\times X_A\times A),
\end{equation}
where $(x,\xi,y)\in A\times\hat{A}\times A$.
\end{ex}

We have the following version of Proposition \ref{ker-equiv-prop} for projective representations. The proof is analogous, so we skip it.

\begin{prop}
Let $X$ and $Y$ be schemes over $S$, 
and let $(\PP,\QQ,\a,\b)$ be adjoint kernel data, where
$\PP\in D_{qc}(X\times Y)$ and $\QQ\in D_{qc}(Y\times X)$.
Let also $G$ be a group scheme, and let $V$ be a kernel projective
representation of $G$ over $X$ with a $2$-cocycle $\LL$.
Then $\QQ\circ_X V\circ_X\PP$ has a natural structure of
a kernel projective representation of $G$ over $Y$ with the same $2$-cocycle. 
This correspondence is compatible with restriction under homomorphisms of group schemes 
and with the functors
of the form $V\mapsto V\ot p_G^*\MM$ for $\MM\in\PPic(G)$
(where $\LL$ gets replaced by $\LL\ot\de\MM$).
Similar assertions hold if $G$ is a formal group scheme, $X$ and $Y$ are proper over $S$, and
$V\in D_{G,\clim}(X\times G\times X)$.
\end{prop}

Recall that the Fourier-Mukai transform $D_{qc}(A)\simeq D_{qc}(\hat{A})$ corresponds
to taking as $\PP$ the Poincar\'e line bundle on $A\times\hat{A}$ and 
$\QQ=\PP^{-1}\ot p_{\hat{A}}^*\om_{\hat{A}/S}[g]$ on $\hat{A}\times A$, where $p_{\hat{A}}$
is the projection to $\hat{A}$, $g$ is the relative dimension of $A/S$.
Using the fact that this transform exchanges the operation of translation on $A$ and of
tensoring with the corresponding line bundle on $\hat{A}$, 
we can easily calculate the Fourier-Mukai dual of the projective kernel representation of
$X_A=A\times\hat{A}$ over $A$ considered above. 

\begin{prop}\label{Heis-Four-prop} 
The Fourier-Mukai dual of the projective representation $V(A)$ of 
the group $X_A$ is given by
\begin{equation}\label{FM-dual-proj-eq}
\QQ\circ_A V(A)\circ_A\PP\simeq p_X^*\PP^{-1}\ot 
(\id_{\hat{A}}\times\kappa\times\id_{\hat{A}})^*V(\hat{A}),
\end{equation}
where $\kappa$ is the isomorphism
$$\kappa:X_A=A\times\hat{A}\to X_{\hat{A}}=\hat{A}\times A:(x,\xi)\mapsto(\xi,-x),$$
and $p_X:\hat{A}\times X_A\times\hat{A}\to X_A$ is the projection.
The isomorphism \eqref{FM-dual-proj-eq} is compatible with the kernel projective representation
structures on $V(A)$ and $V(\hat{A})$ via the isomorphism $\kappa:X_A\wt{\to}X_{\hat{A}}$ and
the isomorphism of $2$-cocycles of $X_A$
$$\La(\PP)\ot(\kappa\times\kappa)^*\LL_{\hat{A}}\simeq\LL_A.$$
\end{prop}

\Pf . Let us consider two simpler kernel representations related to $V(A)$: the one corresponding
to the action of translations, and the one corresponding to tensoring with line bundles associated
with points in $\hat{A}$.
The first is the kernel representation of $A$ (viewed as a group scheme) over $A$ given by
$$V_t(A):=(x+y,x,y)_*\OO_{A\times A}\in\Coh(A\times A\times A),$$
where $(x,y)\in A\times A$.
The second is the kernel representation of $\hat{A}$ over $A$ given by
$$V_{\otimes}(A)=(x,\xi,x)_*\PP\in\Coh(A\times\hat{A}\times A),$$
where $x\in A$, $\xi\in\hat{A}$.
We have natural isomorphisms
\begin{equation}\label{VA-eq1}
V_t(A)\circ_A V_{\otimes}(A)\simeq V(A),
\end{equation}
\begin{equation}\label{VA-eq2}
V_{\otimes}(A)\circ_A  V_t(A)\simeq p_X^*\PP\otimes V(A)
\end{equation}
that are compatible with the structures of kernel representations on $V_t(A)$ and $V_{\otimes}(A)$
and with the structure of a projective kernel representation on $V(A)$ in the following sense.
Combining \eqref{VA-eq1} and \eqref{VA-eq2} we get the following ``commutation" relation:
\begin{equation}\label{VA-eq3}
V_{\otimes}(A)\circ_A V_t(A)\simeq p_{23}^*\PP\otimes V_t(A)\circ_A V_{\otimes}(A).
\end{equation}
Now the compatibility asserts that the projective kernel representation structure on $V(A)$ is
given by the following composition:
\begin{align*}
&V(A)\circ_A V(A)\simeq V_t(A)\circ_A V_{\otimes}(A)\circ_A V_t(A)\circ_A V_{\otimes}(A)\simeq
p_{23}^*\PP\otimes V_t(A)\circ_A V_t(A)\circ_A V_{\otimes}(A)\circ_A V_{\otimes}(A)\to\\
&p_{23}^*\PP\otimes V_t(A)\circ_A V_{\otimes}(A)\simeq p_{23}^*\PP\circ V(A),
\end{align*}
where the first and the last isomorphisms are induced by \eqref{VA-eq1}, the second isomorphism is
induced by \eqref{VA-eq3}, and the arrow is given by the kernel representation structures on
$V_t(A)$ and on $V_{\otimes}(A)$. 
Thus, it is enough to compute the Fourier duals of the kernel
representations $V_t(A)$ and $V_{\otimes}(A)$. Using the relation 
$\QQ\circ_A\PP\simeq\De_*{\hat{A}}$ and the biextension structure of $\PP$, one gets
the following isomorphisms of kernel representations:
$$\QQ\circ_A V_{\otimes}(A)\circ_A\PP\simeq V_t(\hat{A}),$$
$$\QQ\circ_A V_t(A)\circ_A\PP\simeq (\id_{\hat{A}}\times [-1]_A\times\id_{\hat{A}})^*V_{\ot}(\hat{A}).$$
One can check that this indeed  leads to \eqref{FM-dual-proj-eq}.
\ed

We are going to use the above result to calculate the Fourier-Mukai duals of some kernel
representations obtained from $V(A)$ by restricting to subgroups over which the 
corresponding $2$-cocycle becomes a coboundary.

\begin{defi} Let $G$ be a group scheme (resp., a formal group scheme) over $S$,
$A$ an abelian scheme over $S$.
A {\it $G$-twisting data} $T=(f,f',\a,\iota)$ for $A$ consists of homomorphisms
$f:G\to A$, $f':G\to \hat{A}$, and of a line bundle $\a$ over $G$ equipped with an isomorphism
of $2$-cocycles of $G$ with values in $\PPic(S)$
$$\iota:\La(\a)\simeq (f\times f')^*\PP,$$ 
where $\PP$ is the Poincar\'e biextension of $A\times\hat{A}$. 
If $T=(f,f',\a,\iota)$ is a $G$-twisting data for $A$ then the {\it dual $G$-twisting data} for $\hat{A}$
is $\hat{T}=(f',-f,\hat{\a},\iota')$, where 
\begin{equation}\label{dual-twist-eq}
\hat{\a}=\a\ot (f,f')^*\PP^{-1},
\end{equation}
and $\iota'$ is induced by $\iota$.
\end{defi}

Note that the $2$-cocycle $(f\times f')^*\PP$ of $G$ is the pull-back under
$(f,f'):G\to X_A$ of the natural $2$-cocycle $\LL_A=p_{32}^*\PP$ of $X_A$.
Hence, the coboundary \eqref{cob-L-A} for $\LL_A$ with values in
$\PPic(A)$ together with the isomorphism $\iota$ induce a $1$-cocycle of $G$ with values in $\PPic(A)$. The underlying line bundle on $G\times A$ is
$$\LL(T):=p_1^*\a\otimes (\id\times f')^*\PP.$$
It is easy to see that a $1$-cocycle $\LL$ of $G$ with values in $\PPic(A)$ appears in this
way iff $\LL|_{g\times A_s}\in\Pic^0(A_s)$ for every geometric point $g$ of $G$ (where $s$ is the image
of $g$ in $S$). Indeed, in this case we can write $\LL$ in the form $p_1^*\a\otimes(\id\times f')^*\PP$
for some line bundle $\a$ over $G$ and some
morphism $f':G\to \hat{A}$. Unraveling the cocycle condition we get that $f'$ should be a homomorphism
and $\a$ should satisfy $\La(\a)\simeq (f\times f')^*\PP$.

We set $\Qcoh^T_G(A)=\Qcoh^{\LL(T)}_G(A)$ (resp., $\Coh^T_G(A)=\Coh^{\LL(T)}_G(A)$).
Note that by Corollary \ref{tw-action-ker-cor}, we have
$\Qcoh^T_G(A)\simeq\AA^G_A(\LL(T))-\mod$ (resp., $\Coh^T_G(A)\simeq\AA^G_A(\LL(T))-\mod^c$).
The following result generalizes Theorem 15.2 of \cite{P-ab}.

\begin{thm}\label{twist-eq-thm} Assume $S$ is semi-separated.
The kernel representations of $G$ over $A$ and $\hat{A}$ associated with $\LL(T)$ and 
$\LL(\hat{T})$ are Fourier-Mukai dual to each other.
Hence, if $G$ is a finite flat group scheme over $S$ then
the corresponding pure kernel algebras $\AA^G_A(\LL(T))$ and $\AA^G_{\hat{A}}(\LL(\hat{T}))$ are
also Fourier-Mukai dual, and we get exact equivalences 
$$D(\Qcoh^T_G(A))\simeq D(\Qcoh^{\hat{T}}_G(\hat{A})),$$ 
$$D^b(\Coh^T_G(A))\simeq D^b(\Coh^{\hat{T}}_G(\hat{A})).$$
The first of these equivalences also holds if $S=\Spec(k)$, where $k$ is a field, and
$G$ is a formal group scheme, {\it ldu-pf} over $k$, such that the action of $\Gal(\ov{k}/k)$ on $G(\ov{k})$
factors through $\Gal(k'/k)$ for some finite field extension $k\sub k'$.
\end{thm}

\Pf . Let us denote by $V(A)|_G\in \Coh(A\times G\times A)$ (resp.,
$V(\hat{A})|_G\in \Coh(\hat{A}\times G\times\hat{A})$) the restriction of the projective
kernel representation $V(A)$ of $X_A$ over $A$, via the homomorphism
$(f,f'):G\to X_A$ (resp., $(f',-f):G\to X_{\hat{A}}$). Then we have natural isomorphisms
of kernel representations of $G$
$$V_{\LL(T)}\simeq p_G^*\a\ot V(A)|_G, \ \ 
V_{\LL(\hat{T})}\simeq p_G^*\hat{\a}\ot V(\hat{A})|_G.$$
Hence, by Proposition \ref{Heis-Four-prop},
$$\QQ\circ_A V_{\LL(T)}\circ_A\PP\simeq p_G^*(\a\ot(f,f')^*\PP^{-1})\ot V(\hat{A})|_G\simeq
V_{\LL(\hat{T})}.$$
By Proposition \ref{ker-equiv-prop}, this implies the isomorphism of kernel algebras
$$\QQ\circ_A \AA^G_A(\LL(T))\circ_A\PP\simeq \AA^G_{\hat{A}}(\LL(\hat{T})).$$
Now the first equivalence of categories follows from Theorem \ref{PR-thm}(i).
Note that we can apply it because the kernel algebras $\AA^G_A(\LL(T))$ and 
$\AA^G_{\hat{A}}(\LL(\hat{T}))$ are pure and of affine type by Proposition \ref{aff-ker-prop}.
To deduce the second equivalence of categories we use Theorem \ref{PR-thm}(iii)
(recall that if $G$ is finite then
the kernel algebra $\AA^G_A(\LL(T))$ is finite).
\ed

\begin{rem} It would be interesting to try to describe all 
exact equivalences between derived categories of twisted 
coherent sheaves over abelian varieties by generalizing the picture of \cite{Orlov-ab}.
The above theorem provides many examples of such equivalences. 
\end{rem}

\begin{ex} Let us explain how the equivalences between modules over algebras of twisted differential operators (tdo's) considered in \cite{PR} fit the above theorem.
Let $A$ be an abelian variety over a field $k$ of characteristic zero, and let
$G$ be a formal group scheme over $k$, 
isomorphic as a formal scheme to $\Spf k[[t_1,\ldots,t_n]]$. Then 
the kernel algebra $\AA=\AA^G_A(\LL(T))$ over $A$ 
associated with a $G$-twisting data $T=(f,f',\a,\iota)$ for $A$
is a D-algebra (see \cite{BB}, and Example 1 of section \ref{ker-alg-sec} above).
Furthermore, the natural filtration of $\om_G$ (given by the infinitesimal neighborhoods
of zero) induces an exhaustive algebra filtration 
$\AA_0\sub\AA_1\sub\ldots$ such that 
the associated graded algebra is isomorphic to the commutative $\OO_A$-algebra
$S(L_G)\otimes\OO_A$, where $L_G$ is the Lie algebra of $G$, $S(L_G)$ is the symmetric algebra
on $L_G$ .
According to Lemma 5.1 of \cite{PR}, such a $D$-algebra $\AA$ is the universal enveloping 
algebra $U^0(\wt{L})$ corresponding to
a structure of a Lie algebroid on $L_G\otimes\OO_A$ and a central
extension $\wt{L}$ of this Lie algebroid by $\OO_A$. The Fourier-Mukai duality of such algebras
(and derived categories of modules over them) developed in section 7 of \cite{PR} matches the
duality picture obtained from Theorem \ref{twist-eq-thm}.
Note that
the algebra $\AA^G_A(\LL(T))$ is a tdo iff the tangent map
$df:L_G\to L_A$ is an isomorphism, where $L_A$ is the Lie algebra of $A$
(see section 5 of \cite{PR}).
In this case $f$ gives an isomorphism of $G$ with the formal group obtained from $A$. 
\end{ex}

Here is an example of dual twisting data that will be relevant for the next section.

\begin{cor}\label{GH-cor} 
Let $A$ be a abelian scheme over a semi-separated scheme $S$, and let $G$ and $H$ be finite
flat commutative group schemes over $S$
equipped with homomorphisms $f:G\to A$, $f':H\to \hat{A}$ along with a trivialization
of the biextension $(f\times f')^*\PP$ of $G\times H$, where $\PP$ is the Poincar\'e
biextension of $A\times\hat{A}$. We can view the line bundle $\LL(f')=(p_3,f'p_2)^*\PP$ on
$G\times H\times A$ as a $1$-cocycle of $G\times H$ with values in
$\PPic(A)$ (where $H$ acts trivially on $A$), and the line bundle $\LL(f)=(fp_1,p_3)^*\PP$ on 
$G\times H\times\hat{A}$ as a $1$-cocycle of $G\times H$ with values in $\PPic(\hat{A})$
(where $G$ acts trivially on $\hat{A}$). Then the kernel algebras
$\AA^{G\times H}_A(\LL(f'))$ and $\AA^{G\times H}_{\hat{A}}(\LL(f))$ are Fourier-Mukai
dual to each other, so we have exact equivalences of categories
$$D(\Qcoh^{\LL(f')}_{G\times H}(A))\simeq D(\Qcoh^{\LL(f)}_{G\times H}(\hat{A})),$$ 
$$D^b(\Coh^{\LL(f')}_{G\times H}(A))\simeq D^b(\Coh^{\LL(f)}_{G\times H}(\hat{A})).$$ 
The first of these equivalences also holds if $S=\Spec(k)$, where $k$ is a field, and
$G$ (resp., $H$) is a formal group scheme, {\it ldu-pf} over $k$, such that the action of $\Gal(\ov{k}/k)$ on $G(\ov{k})$ (resp., $H(\ov{k})$)
factors through $\Gal(k'/k)$ for some finite field extension $k\sub k'$.
\end{cor}

\Pf . By definition, the $1$-cocycle $\LL(f')$ of $G\times H$ with values in $\PPic(A)$ 
has form $\LL(T)$ for the $G\times H$-twisting data for $A$
$$T=(fp_1,f'p_2,\OO_{G\times H},\iota),$$
where $\iota$ is induced by the trivialization of $(f\times f')^*\PP$ on $G\times H$.
Similarly, $\LL(f)=\LL(\wt{T})$, for the $G\times H$-twisting data for $\hat{A}$
$$\wt{T}=(f'p_2,fp_1,\OO_{G\times H},\wt{\iota})$$
with $\wt{\iota}$ induced by the above trivialization.
Note that $\wt{T}$ differs from the twisting data $\hat{T}$ by the automorphism
$\id\times[-1]$ of $G\times H$. Hence, the kernel algebras on $\hat{A}$ 
associated with $\wt{T}$ and $\hat{T}$ are isomorphic.
It remains to use Theorem \ref{twist-eq-thm}.
\ed

\section{Fourier-Mukai duality for orbi-abelian schemes and generalized 1-motives}
\label{FM-sec}


In this section we will apply the techniques of kernel algebras to get versions of the Fourier-Mukai
equivalences in two situations. The first situation arises when we have a homomorphism
$G\to E$, where $G$ is a finite flat commutative
group scheme, $E$ is an extension of an abelian scheme by a finite
flat commutative group scheme, and we consider the category of $G$-equivariant sheaves on $E$. 
It turns out that one can construct a dual data of this kind and an analog of the Fourier-Mukai functor.
The second situation is the generalization of the one considered by Laumon in \cite{Lau}:
here we work over a field $k$ of characteristic zero, $G$ is a formal $k$-group (see
beginning of section \ref{kerrep-sec} for our conventions on formal $k$-groups), 
and $E$ is a commutative algebraic group over $k$.

We start by studying in sections \ref{proper-sec} and \ref{orbi-sec}
the duality in the particular case of the first situation when the base scheme is
a field, since in this case a nice duality functor exists on certain derived category containing
both abelian varieties and finite commutative group schemes. Then we will
proceed to more general setups mentioned above.

\subsection{Duality functor on the derived category of commutative proper group schemes over a field}
\label{proper-sec}

Let $k$ be a field.
For a commutative algebraic group $X$ of finite type over $k$ 
we denote by $X_0$ its connected component of $1$ taken
with reduced scheme structure (this is denoted as $CR(X)$ in \cite{Oort}).
Also for an integer $N\neq 0$ we denote by $X_N$ the kernel of the map $[N]_X:X\to X:x\mapsto Nx$.

\begin{lem}\label{subgr-sch-lem} 
The following conditions for a commutative algebraic group $X$ over $k$ are equivalent:

\noindent
(a) $X$ is proper;

\noindent
(b) $X_0$ is an abelian variety;

\noindent
(c) $X$ is isomorphic to a subgroup scheme of an
abelian variety.
\end{lem}

\Pf . The equivalence of (a) and (b) is well known.
It is clear that (c) implies (a). To prove (b)$\implies$(c)
we use the fact that every finite group scheme can be embedded into an abelian variety as
a subgroup scheme (see \cite{Oort} (15.4)). Let $X$ be any commutative proper group scheme. Then there exists an integer $N>0$ such that $X/X_0$ is annihilated by $N$. Hence, $X_0=NX_0\sub NX\sub X_0$, so that $NX=X_0$. It follows that 
$X/X_N\simeq X_0$ is an abelian variety. Now pick an embedding $X_N\sub A$, where $A$
is an abelian variety. Then $X$ is a subgroup scheme of the induced extension of $X_0$ by
$A$, which is an abelian variety.
\ed

Let $\GG_k^{pr}$ denote the category of commutative proper group schemes over $k$. 
This is a full subcategory of the category
$\GG_k$ of all commutative algebraic groups over $k$. It is well known that
$\GG_k$ is an abelian category (see \cite{SGA3}, exp. $VI_A$, 5.4). 
A nice exposition of the properties of this category can be found in
chapter II of \cite{Oort}.

\begin{lem} The subcategory $\GG_k^{pr}\sub\GG_k$ is a Serre subcategory.
\end{lem}

\Pf . It is clear that $\GG_k^{pr}$ is closed under quotients and subgroup schemes. Assume that 
$$0\to X\to Y\to Z\to 0$$
is an exact sequence in $\GG_k$ with $X,Z\in\GG_k^{pr}$. 
Then the morphism $Y\to Z$ is proper and $Z$ is proper, hence, $Y$ is proper.
\ed

For the theory of torsion pairs and tilting the reader can consult \cite{HRS} (a concise
exposition can be found also in section 5.4 of \cite{BVdB}).

\begin{lem}
Let $\GG_k^f$ and $\AV_k$ be the subcategories of finite group schemes and abelian varieties in $\GG_k^{pr}$,
respectively. Then $(\AV_k,\GG_k^f)$ is a tilting torsion pair in $\GG_k^{pr}$.
\end{lem}

\Pf . All morphisms from an abelian variety to a finite group scheme are trivial.
Also, by Lemma \ref{subgr-sch-lem}, for $X\in\GG_k^{pr}$ we have $X_0\in\AV_k$ and 
$\pi_0 X:=X/X_0\in\GG_k^f$, so $(\AV_k,\GG_k^f)$ is a torsion pair.
The fact that it is tilting follows from Lemma \ref{subgr-sch-lem}.
\ed

\begin{thm} There exists an exact functor 
$\Du:D^b(\GG_k^{pr})^{opp}\to D^b(\GG_k^{pr})$
such that 

\noindent
(i) $\Du\circ \Du\simeq\Id$;

\noindent
(ii) $\Du(A)\simeq \hat{A}$ for an abelian variety $A$;

\noindent
(iii) $\Du(G)\simeq G^*[1]$ for a finite group scheme $G$.
\end{thm}

\Pf . The idea is to use resolutions in terms of abelian varieties.
More precisely, the fact that $(\AV_k,\GG_k^f)$ is a tilting torsion pair in $\GG_k^{\pr}$ gives an
equivalence $D^b(\GG_k^{pr})\simeq D^b(\AV_k)$, where $\AV_k$ is viewed as an exact category 
(see \cite{BVdB}, Lemma 5.4.2; see also \cite{BBD}, Ex.~1.3.23 (iii)). Since the duality for abelian varieties preserves short exact sequences, 
it extends to an exact functor $\Du:D^b(\AV_k)\to D^b(\AV_k)$.
This proves (i) and (ii). Embedding a finite group scheme $G$ into an abelian variety and
using the Isogeny Theorem (see \cite{Oort}, III.19.1) gives (iii).
\ed

\begin{rem} It is easy to see that $\Du$ exchanges the standard $t$-structure on $D^b(\GG_k^{pr})$
with the tilted $t$-structure associated with the torsion pair $(\AV_k,\GG_k^f)$. 
More precisely, for every $X\in\GG_k^{pr}$ the duality sends
the canonical exact sequence 
$$0\to X_0\to X\to \pi_0X\to 0,$$
where $X_0$ is an abelian variety, to
the exact triangle
$$
(\pi_0X)^*[1]\to \Du(X)\to \hat{X_0}\to (\pi_0X)^*[2],
$$
so that $(\pi_0X)^*[1]\simeq \tau_{\le -1}\Du(X)$ and $\hat{X_0}\simeq \tau_{\ge 0}\Du(X)$.
Thus, we have
\begin{equation}\label{H-dual-eq}
H^{-1}\Du(X)\simeq(\pi_0X)^*,
\end{equation}
\begin{equation}\label{H-dual-eq2}
H^0\Du(X)\simeq\hat{X_0}.
\end{equation}
\end{rem}


\subsection{Orbi-abelian varieties}\label{orbi-sec}

Let $k$ be a field.

\begin{defi} An {\it orbi-abelian variety over $k$} is an object $K$ of $D^b(\GG_k^{pr})$ 
such that $H^{-1}K$ is a finite group scheme and $H^iK=0$ for $i\not\in\{-1,0\}$.
\end{defi}

\begin{lem}\label{orbi-ab-lem} 
(i) For any complex $K$ in $\GG_k^{pr}$ 
concentrated in degrees $-1$ and $0$, such that $H^{-1}K$ is finite,
there exists a quasiisomorphic subcomplex of the form $[G\to X]$,
where $G$ is a finite group scheme.
In particular, every orbi-abelian variety 
can be presented by a complex $[G\to X]$ with $G$ finite.

\noindent (ii)
An object $K\in D^b(\GG_k^{pr})$ is an orbi-abelian variety iff $\Du(K)$ is.
Furthermore, we have natural isomorphisms
\begin{equation}\label{orbi-ab-dual-eq}
(H^0\Du(K))_0\simeq \widehat{(H^0K)_0},
\end{equation}
\begin{equation}\label{H-orbi-dual-eq}
H^{-1}\Du(K)\simeq (\pi_0H^0 K)^*.
\end{equation}
\end{lem}

\Pf . (i) 
Let $K=[i:Y\to X]$, where $\ker(i)$ is finite. We claim that there exists a homomorphism $p:X\to Y_0$ such that $p\circ i=[N]_Y$ for some $N>0$ (where $[N]_YY\sub Y_0)$. Indeed,
since $Y_0$ and $X_0$ are abelian varieties the assertion is true
for the induced homomorphism $i_0:Y_0\to X_0$, so we can find $p_0:X_0\to Y_0$ such that
$p_0\circ i_0=[N_1]_{Y_0}$ for some $N_1>0$. If $N_2=|X/X_0|$ then
the multiplication by $N_2$ maps $X$ to $X_0$, so we can view $\wt{p}=N_2p_0$ as a homomorphism
$X\to Y_0$. Finally, the homomorphisms $\wt{p}\circ i$ and $[N_1N_2]_Y$ agree on $Y_0$,
hence, 
$$N_3\wt{p}\circ i=[N_1N_2N_3]_Y,$$ 
where $N_3=|Y/Y_0|$. Thus, we can set $p=N_3\wt{p}$, $N=N_1N_2N_3$, and our claim follows.
In this situation $K$ is quasi-isomorphic to the subcomplex $[i':Y_N\to\ker(p)]$, where $i'$ is induced
by $i$. 

\noindent  
(ii) Let $K=[G\to X]$ be an orbi-abelian variety, where $G$ is a finite group scheme. Then we have an exact triangle
$$X\to K\to G[1]\to X[1].$$ 
The dual of this triangle is
$$G^*\to \Du(K)\to \Du(X)\to G^*[1].$$
The corresponding long exact sequence of cohomology shows that
$H^i\Du(K)=0$ for $i\not\in\{-1,0\}$ and that $H^{-1}\Du(K)$ 
is a subgroup scheme in $H^{-1}\Du(X)\simeq (\pi_0X)^*$ (see \eqref{H-dual-eq}). Hence,
$H^{-1}\Du(K)$ is finite, so $\Du(K)$ is an orbi-abelian variety. The same exact sequence
shows that
$$H^{-1}\Du(K)\simeq\ker((\pi_0X)^*\to G^*)\simeq\coker(G\to\pi_0 X)^*\simeq(\pi_0H^0K)^*.$$
On the other hand, by \eqref{H-dual-eq2}, the abelian variety $H^0\Du(K)_0$ is dual to
$H^0\Du H^0\Du(K)$. Dualizing the standard exact triangle
$$H^{-1}\Du(K)[1]\to\Du(K)\to H^0\Du(K)\to\ldots$$
and passing to cohomology we get an exact sequence
$$0\to H^0\Du H^0\Du(K)\to H^0K\to (H^{-1}\Du(K))^*\to\ldots$$
which induces an isomorphism $H^0\Du H^0\Du(K)\simeq (H^0K)_0$.
\ed

Here is a simple way to realize the dual $\Du(K)$ to an orbi-abelian variety $K$ concretely.
Let $K=[G\to X]$ as in Lemma \ref{orbi-ab-lem}(i).
We can pick a finite subgroup $H\sub X$ such that
$A=X/H$ is an abelian variety (e.g., one can take $H=X_N$, where $X/X_0$ is annihilated
by $N$).
Note that the dual of $L=[G\to A]$ has only one cohomology due to an exact triangle
$$G^*\to \Du(L)\to\hat{A}\to G^*[1].$$
On the other hand, dualizing the exact triangle 
$$K\to L\to H[1]\to K[1]$$
we get
$$H^*\to\Du(L)\to\Du(K)\to H^*[1],$$
hence, the complex $[H^*\to\Du(L)]$ represents $\Du(K)$.

\begin{rem} If $k$ has characteristic zero then the category $\GG_k^{pr}$ has homological dimension
$1$ (see \cite{Oort}, II.14), so every object of $D^b(\GG_k^{pr})$ is isomorphic to the direct sum of
its cohomologies. By duality, this implies that in this case every object of $\GG_k^{pr}$
is a direct sum of an abelian variety and a finite group scheme.
Thus, in characteristic zero the duality for orbi-abelian varieties takes form 
$$\Du(G[1]\oplus (H\oplus A))\simeq H^*[1]\oplus (G^*\oplus \hat{A}),$$
where $A$ is an abelian variety, $H$ and $G$ are finite commutative group schemes.
Now assume that the characteristic of $k$ is $p>0$.
Let $\GG^{pr}_{k,0}\sub\GG_k^{pr}$ denote the full subcategory consisting of
proper groups $G$ such that $\Hom(\a_p,G)=0$.
We claim that $\GG_{k,0}^{pr}$ is a Serre subcategory in $\GG_k^{pr}$. 
Indeed, it is clear that $\GG_{k,0}^{pr}$ is closed under extensions and passing to subobjects. To check that it is closed under quotients it suffices to consider quotients by elementary groups. But for such a group $G$ one has $\Hom(\a_p,G)=0$ iff
$\Ext^1(\a_p,G)=0$ (see \cite{Oort}, II.14), and our claim follows. Note that the category $\GG^{pr}_{k,0}$
also has homological dimension $1$ (see \cite{Oort}, II.14), so the above remarks about orbi-abelian
varieties in characteristic zero apply as well to orbi-abelian varieties in $D^b(\GG^{pr}_{k,0})$.
\end{rem}

\subsection{Abstract duality setup}\label{abs-dual-sec}

Let $\SS$ be a site. We denote by $\Sh_{\SS}$ the category of sheaves of abelian groups on $\SS$.
We assume that we have fixed a certain sheaf of abelian groups $\bG$. We consider the associated
duality functor 
$$\Du:D^b(\Sh_{\SS})\to D^+(\Sh_{\SS}): K\mapsto R\und{\Hom}(K,\bG)[1].$$
Assume also that we have three full subcategories stable under extensions
$\For, \Aff, \Ab\sub \Sh_{\SS}$ with the following properties:

\noindent
(i) $\Du(\Ab)\sub \Ab$, $\Du(\For)[-1]\sub\Aff$, $\Du(\Aff)[-1]\sub\For$;

\noindent
(ii) for $K$ in one of the subcategories $\For$, $\Aff$ or $\Ab$, the natural morphism $K\to\Du\Du(K)$
is an isomorphism.

In this situation for $A\in\Ab$ we set $\hat{A}:=\Du(A)$, while for $G$ either in $\For$ or in $\Aff$ we set $G^*=\Du(G)[-1]$.
Note that (i) and (ii) imply that for $A\in\Ab$ and for $G$ either in $\For$ or in $\Aff$ one has
$\Hom(A,G)=0$. Indeed, we have $\Du(A)\in\Sh_{\SS}$ and $\Du(G)\in\Sh_{\SS}[1]$.
Hence, $\Hom(\DD(G),\DD(A))=0$, and the required vanishing follows from (ii).

\begin{defi} Given the data 
$\De=(\bG,\For,\Aff,\Ab)$ as above, we define a {\it generalized $1$-motive of type $\De$} as
a complex $[G\to E]$ of sheaves of abelian groups on $\SS$ concentrated in degrees $-1$ and
$0$, such that $G\in\For$, and $E$ fits into an exact sequence
$$0\to L\to E\to A\to 0$$
with $L\in\Aff$ and $A\in\Ab$. Morphisms between generalized $1$-motives as above are
simply morphisms in the derived category $D^b(\Sh_{\SS})$.
\end{defi}

We will consider the following two examples of this situation.
The first is when $\SS$ is the fppf site of schemes over a given scheme $S$, $\bG=\G_m$
is the multiplicative group, $\Ab$ consists
of abelian schemes over $S$, while $\For=\Aff$ are finite flat commutative group schemes over $S$.
We denote the corresponding duality type $\De^{fin}_S$ and call generalized $1$-motives of type
$\De^{fin}_S$ {\it orbi-abelian schemes} over $S$.

In the second example $\SS$ is the fppf site of affine schemes over a field $k$ of characteristic
zero, $\bG=\G_m$,
$\Ab$ consists of abelian varieties over $k$, $\For$ are formal $k$-groups 
(see our conventions in the beginning of section \ref{kerrep-sec}), and
$\Aff$ are affine commutative algebraic groups over $k$.
We denote the corresponding duality type $\De^{for}_k$.
To get the usual generalized $1$-motives of Laumon \cite{Lau} one would have
to modify this type by taking $\For$ to be formal $k$-groups without torsion and 
considering only connected groups in $\Aff$.

Sometimes it is convenient to view generalized $1$-motives as a full subcategory in the category of
Picard stacks by associating with $K=[G\to E]$ the Picard stack $\ch(K)$
(see section \ref{biext-sec} and \cite{SGA4}, 1.4). 

\begin{prop} If $K$ is a generalized $1$-motive of type $\De$ then $\Du(K)$ is also a
generalized $1$-motive of type $\De$ and the natural map $K\to\Du\Du(K)$ is an isomorphism.
\end{prop}

\Pf . Let $K=[G\to E]$, where $G\in\For$ and $E$ is an extension of $A\in\Ab$ by $L\in\Aff$.
Set $\ov{K}=[G\to A]$.
Then $\Du(\ov{K})$ is an extension of $\hat{A}=\Du(A)$ by $G^*=\Du(G)[-1]$, and
$\Du(K)$ is represented by the complex $[L^*\to\Du(\ov{K})]$. The second assertion now follows from
the assumption (ii) on the data $\De$.
\ed

The duality of generalized $1$-motives has a convenient interpretation in terms of biextensions.
Recall that
if $[Y\to X]$ and $[Y'\to X']$ are complexes over $\Sh_{\SS}$ concentrated in degrees $-1$ and $0$ then
a {\it biextension of $[Y\to X]$ and $[Y'\to X']$ by $\bG$} is a biextension of $X\times X'$ by
$\bG$ equipped with trivialization of its pull-backs to $Y\times X'$ and $X\times Y'$
such that the induced trivialization of the pull-back to $Y\times Y'$ are the same 
(see \cite{De}, 10.2). Such biextensions form a (commutative) Picard category and
we denote by $\Biext^1([Y\to X], [Y'\to X'];\bG)$ (resp., 
$\Biext^0([Y\to X],[Y'\to X'],\bG)$) the group of isomorphism classes in this category
(resp, automorphism group of an object). 
These groups are isomorphic to $\Ext^i([Y\to X]\otimes[Y'\to X'],\bG)$, $i=0,1$, so they
depend only on isomorphism classes of $[Y\to X]$ and
$[Y'\to X']$ in the derived category $D^b(\Sh_{\SS})$ (see \cite{De}, 10.2.1). 

\begin{prop}\label{dual-prop} 
Let $K$ and $K'$ be generalized $1$-motives of type $\De$. 

\noindent
(i) One has a functorial isomorphism 
$$\Hom_{D^b(\Sh_{\SS})}(K',\Du(K))\simeq\Biext^1(K, K';\bG).$$

\noindent
(ii) For a presentation $K=[G\to E]$ 
let $\Eext([G\to E],\bG)$ denote the Picard stack of extensions of $[G\to E]$ by $\bG$
in the category of complexes over $\Sh_{\SS}$. Equivalently, $\Eext([G\to E],\bG)$
classifies extensions of $E$ by $\bG$ with trivialized 
pull-back to $G$. Then we have an isomorphism of Picard stacks
\begin{equation}\label{dual-ext-eq}
\ch\Du(K)\simeq\Eext(K,\bG).
\end{equation}

\noindent
(iii) Let $K=[G\to E]$, where $E$ is an extension of $A\in\Ab$ by $L\in\Aff$, so that
$\Du(K)$ is represented by $[L^*\to \Du(\ov{K})]$ with $\ov{K}=[G\to A]$. 
Then the biextension of $K\times\Du(K)$ corresponding to the identity map $\id_{\Du(K)}$ under
the isomorphism of (i) ({\em the Poincar\'e biextension})
is represented by the pull-back to $E\times\Du(\ov{K})$ of the
similar biextension of $A\times\hat{A}$, where $\hat{A}=\Du(A)$. This pull-back is 
equipped with natural trivializations along
$G\times \Du(\ov{K})$ and $E\times L^*$ that are compatible over $G\times L^*$. 
\end{prop}

\Pf . (i) We have 
$$\Hom(K',\Du(K))=\Hom(K',R\und{\Hom}(K,\bG[1]))\simeq\Hom(K'\otimes K,\bG[1])=
\Ext^1(K'\otimes K,\bG),$$
which is isomorphic to $\Biext^1(K,K';\bG)$ (see \cite{De}, 10.2.1).

\noindent
(ii) By Lemma \ref{Pic-ext-lem}, we have 
$$\Eext(K,\bG)\simeq HOM(\ch(K),\ch(\bG[1]))\simeq\ch(\tau_{\le 0}R\und{\Hom}(K,\bG[1])),$$
so the assertion follows from $R\und{\Hom}(K,\bG[1])=\Du(K)\in D^{\le 0}(\Sh_{\SS})$.

\noindent
(iii) Note that the Poincar\'e biextension of $K\times\Du(K)$ corresponds to the canonical
morphism $K\otimes\Du(K)\to\bG[1]$ via the isomorphism
$\Biext^1(K,\Du(K);\bG[1])\simeq\Ext^1(K\otimes\Du(K),\bG)$.
Applying this to $\ov{K}$ and using the commutative diagram
\begin{diagram}
A\otimes\Du(\ov{K}) &\rTo{}& A\otimes\hat{A}\\
\dTo{}&&\dTo{}\\
\ov{K}\otimes\Du(\ov{K})&\rTo{}&\bG[1]
\end{diagram}
we see that the Poincar\'e biextension of $\ov{K}\times\Du(\ov{K})$ is represented by
the pull-back to $A\times\Du(\ov{K})$ of the Poincar\'e biextension of $A\times\hat{A}$
(this pull-back is equipped with a trivialization over $G\times\Du(\ov{K})$).
Now the assertion follows in a similar way from the commutative diagram
\begin{diagram}
K\otimes\Du(\ov{K}) &\rTo{}& \ov{K}\otimes\Du(\ov{K})\\
\dTo{}&&\dTo{}\\
K\otimes\Du(K)&\rTo{}&\bG[1]
\end{diagram}
\ed

\subsection{Fourier-Mukai transform}

In this section we work with generalized $1$-motives of type $\De=(\G_m,\For,\Aff,\Ab)$, 
where $\De$ is either
$\De^{fin}_S$ (orbi-abelian schemes) or $\De^{for}_k$. 
Thus, our generalized $1$-motives are of the form $K=[G\to E]$,
where either (i) $E$ is an extension of an abelian scheme by a finite flat commutative group scheme
over $S$, $G$ is a finite flat commutative group scheme over $S$, or 
(ii) $E$ is a commutative algebraic group over a field 
$k$ of characteristic zero, $G$ is a formal $k$-group. 
Note that in case (ii) $E$ is an extension of an abelian variety by an affine commutative algebraic group.
In both cases by a quasicoherent sheaf on $K$ we mean a $G$-equivariant quasicoherent sheaf on
$E$ (where $G$ acts on $E$ by translations). We denote the category of quasicoherent sheaves
by $\Qcoh K$. In case (i) we can also consider the subcategory $\Coh K$ of coherent sheaves.
It is easy to check that  up to equivalence the category $\Qcoh K$ 
does not depend on a presentation $K=[G\to E]$. Indeed, this is clear in the case when we have surjective morphism of complexes $f:[G'\to E']\to [G\to E]$ (i.e., both maps $f_{-1}:G'\to G$ and
$f_0:E'\to E$ are surjective) such that $\ker(f_{-1})$ maps isomorphically to $\ker(f_0)$ (cf. Lemma \ref{gerb-lem}---note that $\ker(f_{-1})\simeq \ker(f_0)$ is a finite group scheme). The general
case follows because of the following simple result. 

\begin{lem}\label{mot-present-lem} 
Let $[G\to E]$ and $[G'\to E']$ be two presentations of the same generalized $1$-motive
$K$. Then there exists a third presentation $[G_0\to E_0]$ of $K$ equipped with surjective
maps (quasi-isomorphisms) to $[G\to E]$ and $[G'\to E']$.
\end{lem} 

\Pf . Let us define $E_0$ in the derived category of sheaves from the exact triangle
$$E_0\to E\oplus E'\to K\to E_0[1].$$
Using the octahedron axiom one can easily see that $E_0$ is an extension of $E$ by $G'$
(resp., of $E'$ by $G$), and that $[G\oplus G'\to E_0]$ will be the required third presentation of $K$.
\ed

Let $K=[G\stackrel{f}{\ra} E]$ be a generalized $1$-motive of one of the two types above.
By definition, $E$ is an extension of $A\in\Ab$ by $L\in\Aff$. Note that $A$ is an abelian scheme over
$S$ (where $S=\Spec(k)$ in case (ii)). 
Let $\pi:E\to A$ denote the projection. Then we associate with the presentation $K=[G\to E]$ the
kernel algebra 
$$\AA(G\to E)=(\pi\times\pi)_*\AA^G_E$$ 
over $A$. By Corollary \ref{action-ker-cor} and Lemma \ref{push-alg-lem},
the category $\Qcoh K$ is equivalent
to the category of $\AA(G\to E)$-modules on $A$.
In the case $\De=\De^{fin}_S$ the kernel algebra $\AA(G\to E)$ is finite and we have an equivalence
$\Coh K\simeq\AA(K)-\mod^c$.

As we have seen before, the dual $1$-motive is represented as $\Du(K)=[L^*\to E']$,
where $E'=\Du[G\to A]$ is an extension of $\hat{A}$ by $G^*$.
Thus, denoting by $\rho:E'\to\hat{A}$ the natural map we also have the corresponding kernel algebra
$$\AA(L^*\to E')=(\rho\times\rho)_*\AA^{L^*}_{E'}$$
over $\hat{A}$, such that $\Qcoh \Du(K)$ is equivalent to the category of modules over $\AA(L^*\to E')$. 

\begin{thm}\label{FM-thm} 
The kernel algebras $\AA(G\to E)$ and 
$\AA(L^*\to E')$ are Fourier-Mukai dual to each other.
Hence, we get an exact equivalence
$$D(\Qcoh K)\simeq D(\Qcoh \Du(K))$$ 
that also induces an equivalence between the bounded
derived categories of coherent sheaves in the case $\De=\De^{fin}_S$. 
\end{thm}

\Pf . Set $H=L^*$.
By construction, the extension $E\to A$ by $L$ is dual to the homomorphism $f':H\to\hat{A}$
(in the sense of Lemma \ref{ext-Car-lem}), while the extension $E'\to\hat{A}$ by $G^*$ is dual to the homomorphism $f:G\to A$.
Applying Corollary \ref{equiv-ker-cor} to the $G$-equivariant $L$-torsor $E\to A$
(resp., $H$-equivariant $G^*$-torsor $E'\to\hat{A}$), we get isomorphisms of kernel algebras 
on $A$ and on $\hat{A}$:
$$\AA(G\to E)\simeq \AA^{G\times H}_A(\LL(f')),$$
$$\AA(H\to E')\simeq \AA^{G\times H}_{\hat{A}}(\LL(f)),$$
where we use the notation of Corollary \ref{GH-cor}. It remains to apply this Corollary.
\ed

\subsection{Fourier-Mukai duality for twisted sheaves}

As in the previous section we work with generalized $1$-motives of type 
$\De\in\{\De^{fin}_S,\De^{for}_k\}$.
Let $K$ be such a generalized $1$-motive.
Let also $G$ be an object of $\For$, so $G$ is a finite flat commutative group scheme over $S$ if 
$\De=\De^{fin}_S$, and $G$ is a $k$-formal group if $\De=\De^{for}_k$. 
The definition of twisting data from section \ref{projFM-sec} 
has an obvious extension to this situation.

\begin{defi} A {\it $G$-twisting data} $T=(f,f',\a,\iota)$ for $K$ consists of homomorphisms
$f:G\to K$, $f':G\to \Du(K)$, and of a line bundle $\a$ over $G$ equipped with an isomorphism
of $2$-cocycles 
$\iota:\La(\a)\simeq (f\times f')^*\PP$, where $\PP$ is the Poincar\'e biextension. 
If $T=(f,f',\a,\iota)$ is a $G$-twisting data for $K$ then the {\it dual $G$-twisting data} for $\Du(K)$
is $\Du(T)=(f',-f,\hat{\a},\iota')$, where $\hat{\a}$ is given by \eqref{dual-twist-eq} 
and $\iota'$ is induced by $\iota$.
\end{defi}

We will use $G$-twisting data as above to define twisted versions of the category of $G$-equivariant
sheaves on $K$. By Proposition \ref{dual-prop}(i),  
twisting data $T=(f,f',\a,\iota)$ can be equivalently described by the data $(f,\BB,\a,\iota)$, 
where $f:G\to K$ is a homomorphism, $\BB$ is a biextension
of $G$ and $K$ by $\G_m$, $\a$ is a line bundle over $G$ and 
$\iota:\La(\a)\simeq (f\times\id)^*\BB$ is an isomorphism of $2$-cocycles.

Given a representation $K=[H\to E]$ let us consider the
composition of $f:G\to K$ with the corresponding morphism $K\to H[1]$.
This will give an extension 
$$0\to H\to\wt{H}(f)\to G\to 0$$
such that the composition $\wt{H}(f)\to G\to K\to H[1]$ is zero.
We claim that to a morphism $f$ one can canonically 
(up to an automorphism of $\wt{H}(f)$, compatible with the extension structure)
associate a morphism of exact triangles 
\begin{diagram}
H&\rTo{}&\wt{H}(f)&\rTo{}&G &\rTo{}& H[1]\\
\dTo{\id}&&\dTo{\wt{f}}&&\dTo{f} &&\dTo{\id}\\
H&\rTo{}&E&\rTo{}& K &\rTo{}& H[1]
\end{diagram}
Indeed, a morphism in the derived category $f: G\to K$ can be represented
by a map $[G_1\to G_0]\to [H\to E]$
in the homotopy category of complexes, where $G_1\sub G_0$ and
$G_0/G_1\simeq G$. To such a map of complexes we associate the
natural homomorphism $(G_0\oplus H)/G_1\to E$, identical on $H$.
It remains to observe that one has a natural isomorphism $(G_0\oplus H)/G_1\simeq\wt{H}(f)$
of extensions of $G$ by $H$, and that our construction is compatible with the homotopy and
with changing the complex $[G_1\to G_0]$. 
Mimicking the definition for the abelian varieties, let us consider the line bundle
\begin{equation}\label{LT-eq}
\LL(T)=p_1^*\a\ot\BB
\end{equation} 
on $G\times E$ (where $\BB$ is the biextension coming from our
twisting data). Its pull-back $\wt{\LL}(T)$ to $\wt{H}(f)\times E$ has the natural structure of 
a $1$-cocycle of $\wt{H}(f)$ with values in $\PPic(E)$, where the action of $\wt{H}(f)$ on $E$
is induced by the homomorphism $\wt{f}$.
Now we define the twisted category of sheaves on $[K/G]$ associated with $T$ to be
\begin{equation}\label{orbi-twisted-eq}
\Qcoh^T_G(K)=\AA^{\wt{H}(f)}_E(\wt{\LL}(T))-\mod.
\end{equation}
If we have another representation $K=[H'\to E']$ such that $f:G\to K$ factors through $E'\to K$ then
we can construct another presentation $K=[H_0\to E_0]$ using Lemma
\ref{mot-present-lem}, so that we have an exact triangle
$$E_0\to E\oplus E'\to K\to E_0[1].$$
This easily implies that $f:G\to K$ factors also through $E_0\to K$. Also, we have morphisms
from $[H_0\to E_0]$ to the original two presentations of $K$.
In this situation the corresponding $1$-cocycle $\wt{\LL}(T_0)$
of $\wt{H}_0(f)$ with values in $\PPic(E_0)$
is isomorphic to the pull-back of $\wt{\LL}(T)$ with respect to the natural
morphisms $\wt{H}_0(f)\to \wt{H}(f)$ and
$E_0\to E$. Hence, we get an equivalence
$$\AA^{\wt{H}_0(f)}_{E_0}(\wt{\LL}(T_0))-\mod\simeq\AA^{\wt{H}(f)}_E(\wt{\LL}(T))-\mod
$$
(see Lemma \ref{gerb-lem}), 
which shows that the right-hand side of \eqref{orbi-twisted-eq} does not depend
on the choice of a representation $K=[H\to E]$.

Now we will prove a generalization of Theorem \ref{twist-eq-thm} in this situation by constructing
an equivalence of the derived categories of twisted sheaves on $[K/G]$ and $[\Du(K)/G]$ 
associated with $T$ and $\Du(T)$, respectively. 

\begin{thm}\label{orbi-twist-eq-thm} 
In the above situation one has an exact equivalence
$$D(\Qcoh^T_G(K))\simeq D(\Qcoh^{\Du(T)}_G(\Du(K)).$$
In the situation of orbi-abelian schemes (i.e., $\De=\De^{fin}_S$), this equivalence
induces an equivalence between the bounded derived categories of coherent sheaves.
\end{thm}

\Pf . We start by choosing dual representations for $K$ and $\Du(K)$.
Namely, let $K=[H_1\to E]$ be a representation of $K$, where $H_1\in\For$ and
$E$ is an extension
$$0\to L_1\to E\to A\to 0,$$
where $A\in\Ab$ and $L_1\in\Aff$. 
Consider the corresponding representation for the dual $1$-motive
$\Du(K)=[H_2\to E']$, where $E'=\Du[H_1\to A]$, and $H_2=L_1^*$,
so that we have an exact sequence
$$0\to L_2\to E'\to\hat{A}\to 0$$
with $L_2=H_1^*$.

We want to express everything in terms of kernel algebras over $A$ and $\hat{A}$, 
so that the dual sides enter into the picture in a symmetric way (as in Theorem \ref{twist-eq-thm}).
Recall that the two categories we want to compare are
$$\Qcoh^T_G(K)=\AA^{\wt{H}_1(f)}_E(\wt{\LL}(T))-\mod,$$
$$\Qcoh^{\Du(T)}_G(\Du(K))=\AA^{\wt{H}_2(f')}_{E'}(\wt{\LL}(\Du(T))-\mod),$$
where $\wt{H}_1(f)$ and $\wt{H}_2(f')$ are extensions of $G$ corresponding to the composed morphisms
$G\stackrel{f}{\to} K\to H_1[1]$ and $G\stackrel{f'}{\to}\Du(K)\to H_2[1]$, respectively. Recall that to define
the relevant $1$-cocycles we also use
the morphisms $\wt{f}:\wt{H}_1(f)\to E$ and $\wt{f'}:\wt{H}_2(f')\to E'$ lifting $f:G\to K$ and 
$f':G\to\Du(K)$, respectively.

Let $\wt{G}(f,f')$ denote the extension of $G$ by $H_1\oplus H_2$ such that the corresponding
class in $\Ext^1(G,H_1\oplus H_2)\simeq\Ext^1(G,H_1)\oplus\Ext^1(G,H_2)$ has
components represented by the extensions $\wt{H}_1(f)$ and $\wt{H}_2(f')$.
Then we have an exact sequence
$$0\to H_2\to \wt{G}(f,f') \to \wt{H}_1(f)\to 0.$$
Recall that $H_2$ is Cartier dual to $L_1$, so 
this extension is dual (in the sense of Lemma \ref{ext-Car-lem}(i)) 
to some biextension $\BB_1$ of $\wt{H}_1(f)\times L_1$ by $\G_m$.
It is easy to check that $\BB_1$ is isomorphic to the pull-back of 
the biextension $\BB$ of $G\times E$ (corresponding to the homomorphism $f':G\to \Du(K)$)
under the natural homomorphism $\wt{H}_1(f)\times L_1\to G\times E$.
It follows from the definition that $\wt{\LL}(T)$ is $(L_1,\BB_1)$-equivariant. Hence,
we are in the situation of Proposition \ref{equiv-ker-prop}, where we view $E$ as an $L_1$-torsor
over $A$. 
Using this Proposition we obtain an isomorphism of kernel algebras over $A$
\begin{equation}\label{twist-phi-isom}
(\phi\times\phi)_*\AA^{\wt{H}_1(f)}_{E}(\wt{\LL}(T))\simeq\AA^{\wt{G}(f,f')}_A(\wt{\LL}_1),
\end{equation}
where $\phi:E\to A$ is the projection, and
$\wt{\LL}_1$ is the $1$-cocycle of $\wt{G}(f,f')$ with values in $\PPic(A)$ 
associated with $\wt{\LL}(T)$ by the construction of Proposition \ref{equiv-ker-prop}.
Thus, by Lemma \ref{push-alg-lem}, we get an equivalence
$$\Qcoh^T_G(K)\simeq\AA^{\wt{G}(f,f')}_A(\wt{\LL}_1)-\mod.$$ 
It is not hard to check that in fact $\wt{\LL}_1$ 
comes from a natural $\wt{G}(f,f')$-twisting data $\wt{T}_1$ for $A$.
More precisely, this twisting data consists of homomorphisms
$$\bff:\wt{G}(f,f')\to\wt{H}_1(f)\stackrel{\wt{f}}{\to}E\to A,$$
$$\bff':\wt{G}(f,f')\to\wt{H}_2(f')\stackrel{\wt{f'}}{\to}E'\to \hat{A},$$
and of the line bundle $\pi^*\a$ on $\wt{G}(f,f')$, where $\pi:\wt{G}(f,f')\to G$ is the projection,
equipped with the isomorphism $\La(\pi^*\a)\simeq(\bff\times\bff')^*\PP$ induced by $\iota$.
To see that the corresponding $1$-cocycle of $\wt{G}(f,f')$ is isomorphic to $\wt{\LL}_1$
one should look at their pull-backs to $\wt{G}(f,f')\times E$ and observe that
these pull-backs are both isomorphic to the pull-back of the line
bundle $\LL(T)$ on $G\times E$ given by \eqref{LT-eq}.
Indeed, this follows from the fact that the pull-back of $\BB$ to $\wt{G}(f,f')\times E$
coincides with the pull-back
under the homomorphism $\wt{G}(f,f')\times E\to E'\times A$ of 
the natural biextension of $E'\times A$.
It is easy to check that the above isomorphism on $\wt{G}(f,f')\times E$ is
compatible with the $L_1$-action and with the $1$-cocycle structures.

Next, we should repeat the above procedure for the dual data $(\Du(K),\Du(T))$.
Note that the extension $\wt{G}(f,f')$ will get replaced by $\wt{G}(f',-f)$ (an extension
of $G$ by $H_2\oplus H_1$) that maps
to $\wt{H}_2(f')$ (an extension of $G$ by $H_2$) and to $\wt{H}_1(-f)$ (an extension of $G$
by $H_1$).
The result will be an isomorphism similar to \eqref{twist-phi-isom}
$$(\phi'\times\phi')_*\AA^{\wt{H}_2(f')}_E(\wt{\LL}(T))\simeq\AA^{\wt{G}(f',-f)}_{\hat{A}}(\wt{\LL}_2),$$
of kernel algebras over $\hat{A}$,
where $\phi':E'\to\hat{A}$ is the projection, and $\wt{\LL}_2$ is the $1$-cocycle of
$\wt{G}(f',-f)$ with values in $\PPic(\hat{A})$ coming from a $\wt{G}(f',-f)$-twisting data
$\wt{T}_2$ for $\hat{A}$ that is defined similarly to $\wt{T}_1$. Again, by Lemma \ref{push-alg-lem},
this gives an equivalence
$$\Qcoh^{\Du(T)}_G(\Du(K))\simeq\AA^{\wt{G}(f',-f)}_{\hat{A}}(\wt{\LL}_2)-\mod.$$ 

We have a natural isomorphism 
$\si:\wt{H}_1(f)\to\wt{H}_1(-f)$ inducing
$[-1]$ on $H_1$ and identity on $G$.
Since $\wt{G}(f',-f)$ is built from the extensions $\wt{H}_2(f')$ and $\wt{H}_1(-f)$, we get the induced
isomorphism
$$\tau:\wt{G}(f,f')\to\wt{G}(f',-f)$$
compatible with the projections to $\wt{H}_2(f')$ and fitting into the commutative square
\begin{diagram}
\wt{G}(f,f')&\rTo{\tau}&\wt{G}(f',-f)\\
\dTo{}&&\dTo{}\\
\wt{H}_1(f)&\rTo{\si}&\wt{H}_1(-f)
\end{diagram}
The map $\wt{G}(f',-f)\to A$ forming a part of the twisting data $\wt{T}_2$ is given by the composition
\begin{equation}\label{T2-map}
\wt{G}(f',-f)\to\wt{H}_1(-f)\stackrel{h}{\to} E\to A,
\end{equation}
where $h$ fits into the morphism of exact triangles
\begin{diagram}
H_1&\rTo{}&\wt{H}_1(-f)&\rTo{}&G &\rTo{}& H_1[1]\\
\dTo{\id}&&\dTo{h}&&\dTo{-f} &&\dTo{\id}\\
H_1&\rTo{}&E&\rTo{}& K &\rTo{}& H_1[1]
\end{diagram}
It follows that the composition $h\circ\si:\wt{H}_1(f)\to E$ differs from $-\wt{f}$ by an automorphism
of $\wt{H}_1(f)$, inducing identity on $H_1$ and $G$. Adjusting $\si$ (and hence $\tau$)
by this automorphism we can assume that $h\circ\si=-\wt{f}$.
This implies that under the isomorphism $\tau:\wt{G}(f,f')\to\wt{G}(f',-f)$ the map
\eqref{T2-map} gets identified with $-\bff$. On the other hand, since $\tau$ is compatible
with the projections to $\wt{H}_2(f')$, the map $\wt{G}(f',-f)\to\hat{A}$ forming a part of
the twisting data $\wt{T}_2$ gets identified under $\tau$ with $\bff'$. It is easy to see from this
that under the isomorphism $\tau$ the data $\wt{T}_2$ gets identified with the dual twisting
data to $\wt{T}_1$.
Now the required equivalence follows from Theorem \ref{twist-eq-thm}.
\ed

\begin{rem} One can generalize the notion of a $G$-twisting data to allow $G$
to be noncommutative (as we did in section \ref{projFM-sec}). 
Namely, let $G$ be a finite flat group scheme over 
$S$ (resp., a formal group scheme, {\it ldu-pf} over $k$, 
such that the action of $\Gal(\ov{k}/k)$ on $G(\ov{k})$ factors through $\Gal(k'/k)$ for some finite field extension $k\sub k'$) equipped with a homomorphism $\pi:G\to G_0$ with $G_0\in\For$.
Then we can consider twisting data consisting of morphisms
$f:G_0\to K$, $f':G_0\to \Du(K)$, and of a line bundle $\a$ over $G$ equipped with an isomorphism
of $2$-cocycles $\iota:\La(\a)\simeq (f\pi\times f'\pi)^*\PP$ over $G\times G$.
One can still define the dual twists of  the categories of $G$-equivariant sheaves on $K$ and on 
$\Du(K)$  and prove the corresponding equivalence of derived categories.
\end{rem}

\end{document}